\documentclass[11pt]{article}
\usepackage{amsmath,amssymb,amsthm}
\usepackage{fullpage}
\usepackage{mathtools}
\usepackage{tikz,pgfplots}
\usepackage{color}
\usepackage{combelow}

\usepackage[breaklinks,colorlinks=true,linkcolor=blue,citecolor=red,backref=page]{hyperref}

\def\proscal<f,g>{\langle\#1,\#2\rangle}
\DeclareMathAlphabet{\itbf}{OML}{cmm}{b}{it}
 \DeclareMathAlphabet\mathbfcal{OMS}{cmsy}{b}{n}

\newcommand{\EE}{\mathbb{E}}
\newcommand{\PP}{\mathbb{P}}
\newcommand{\RR}{\mathbb{R}}
\newcommand{\cS}{\mathcal{S}}
\def\eps{\varepsilon}

\def\bZ{{\itbf Z}}
\def\bX{{\itbf X}}
\def\bx{{\itbf x}}
\def\by{{\itbf y}}
\def\bz{{\itbf z}}
\def\bU{{\itbf X}^0}
\def\bW{{\itbf W}}

\def\bb{{\itbf b}}
\def\bg{{\itbf g}}

\def\bsigma{\boldsymbol{\sigma}}
\def\boeta{\boldsymbol{\eta}}
\def\bGamma{\boldsymbol{\Gamma}}
\def\bpsi{\boldsymbol{\psi}}

\newtheorem{thm}{Theorem}[section]

\newtheorem{lem}[thm]{Lemma}
\newtheorem{prop}[thm]{Proposition}
\newtheorem{example}[thm]{Example}
\newtheorem{definition}[thm]{Definition}

\newcommand{\white}[1]{{\textcolor{white}{#1}}}

\usepackage{slashbox}
\usepackage{amsfonts,amssymb,amsmath,amsthm}
\usepackage{algorithmicx, algorithm}
\usepackage{algcompatible}
\usepackage{color,ulem}
\usepackage{fullpage}
\usepackage{graphicx}%
\usepackage[active]{srcltx}
\usepackage{mathabx}
\usepackage{verbatim}
\usepackage{cancel}
\usepackage{mathtools}
 \usepackage{pgfplots}
\usepackage{tikz,pgfplots}
\usetikzlibrary{snakes}
\usepackage{cite}
\usepackage{hyperref}
\hypersetup{
  linkcolor  = green!60!black,
  citecolor  = blue!60!black,
  urlcolor   = red!60!black,
  colorlinks = true,
}
 \usepackage{tikz}
 \usetikzlibrary{calc,patterns,
                 decorations.pathmorphing,
                 decorations.markings}
                 
\begin{document}
 
\title{A control variate method driven by diffusion approximation}

\author{Josselin Garnier and Laurent Mertz}

\maketitle

\begin{abstract}
In this paper we examine a control variate estimator
for a quantity that can be expressed as the expectation of a function of a random process, that is itself the solution
of a differential equation driven by fast mean-reverting ergodic forces.
The control variate is the same function for the limit diffusion process that approximates
the original process when the mean-reversion time goes to zero.
To get an efficient control variate estimator, we propose a coupling method for the original process and the limit diffusion process.
We show that the correlation between the two processes indeed goes to one when the mean reversion time goes to zero
and we quantify the convergence rate, which makes it possible to characterize the variance reduction of the proposed control variate method.
The efficiency of the method is illustrated on a few examples.
\end{abstract}

\section{Introduction}

In this paper we consider a system driven by external time-dependent random forces and we aim to compute a quantity of interest that is the expectation of a function of the system. 
The system state is the solution of an ordinary differential equation (or a system of ordinary differential equations) driven by external forces which are modeled as stationary random processes. The driving processes may have complicated spectra that have to be taken into account to compute the quantity of interest. This happens for instance in seismic probabilistic risk assessment studies or in the analysis of the structural performance of installations under seismic excitations \cite{patil} or under other loading sources such as wind or waves \cite{kameshwar,quilligan}.
For instance, the reliability of complex systems such as fixed or floating offshore wind turbines depends on its resistance against fatigue damage.  Fatigue damage can be assessed by time-domain simulations in which the structure is subjected to wind, wave and current loads \cite{cordle}.
The different loads can be described by (locally) stationary Gaussian processes with tabulated power spectral densities (such as the JONSWAP spectrum \cite{hassemann}). We may then wish to estimate the mean cumulative fatigue damage or a probability of failure which corresponds to the exceedance of a threshold value.

Monte Carlo simulations are standard to estimate the quantities of interest but they may be very time consuming. We look for an efficient variance reduction technique in this framework. It is known from the diffusion approximation theory \cite{papa74,e05,book} that the driving forces can often be approximated by white noises and the responses of the system can then be modeled by stochastic differential equations. 
This makes it possible to implement a partial differential equation approach to compute the quantity of interest. 
However, the bias due to the approximation of the original driving force by a white noise may be significant and  difficult to assess. 
To compensate for this bias, one may think at a control variate method \cite{glasserman}. 
Such a strategy has already been implemented in a Markov chain Monte Carlo context,  where the goal was to sample
 from a complex invariant probability distribution of a Markov chain for which an approximate distribution has a known expression.
The expectation of the approximate distribution then provides an initial guess, which can be corrected by simulating the two coupled processes to estimate the difference (in expected values) between the true distribution and the approximate distribution \cite{goodman}.
The implementation of a control variate method in our framework requires to be able to simulate the system driven by the original driving force with its complicated spectrum and the limit system driven by the white noise in such a way that both systems are strongly correlated. 
Unfortunately, most diffusion approximation results are established in a weak sense \cite{e05,book}. Some strong results have been obtained 
but only when the drift is a term of order one \cite{kifer01,givon06,liu10}, not when it is a zero-mean large term as we deal with in this paper.
In this paper we build an efficient coupling between the original and limit systems,
we establish a strong convergence result by quantifying the mean square distance between the original and limit processes, 
and we characterize the variance reduction of the control variate method.
We show by our theoretical results and numerical simulations that the variance reduction can be dramatic.

Our method is relevant when the quality of the approximation of the driving forces by a white noise is moderate. 
If it is very accurate,
then the quantity of interest can be estimated (up to a very small and negligible bias) by 
resolution of a Kolmogorov equation based on the limit diffusion system 
(or by a brute force Monte Carlo method applied to the limit system), so there is no need to apply a control variate method.
If it is very poor,
then the limit diffusion system is not correlated to the original system and the control variate method is not efficient.
If it is moderate, then the bias of the estimation method that consists in replacing the original system by the limit one is non-negligible, and the two systems are correlated, so that the control variate method turns out to be very efficient.

The paper is organized as follows.
In Section \ref{sec:ode} we introduce the random ordinary differential equations addressed in this paper and we state the main results of the paper. Note that, motivated by applications in engineering mechanics and physics such as the study of the risk analysis of failure for mechanical structures subjected to random vibrations \cite{MR3733758,MR3335058,MR2989391,siadspaper} or the modeling of the stochastic dynamics of fluid-structure interaction in turbulent thermal convection \cite{MR3853037}, we also consider the case of multivalued ordinary differential equations.
Sections \ref{sec:adif1}-\ref{sec:num} consider random ordinary differential equations.
In Section \ref{sec:adif1} we state the diffusion approximation theorem that gives the convergence in probability of the original process to the limit process. 
In Section \ref{sec:num} we apply the control variate method to a few examples.
The results are extended to the multivalued case in Sections \ref{sec:adif2}-\ref{sec:nummulti}.
In particular Section \ref{sec:nummulti} report numerical results for relevant engineering mechanics problems.
The concluding remark of Section \ref{sec:conclu} connects our findings to the multilevel Monte Carlo literature.

\section{Main results}
\label{sec:ode}
We consider the $\RR^n$-valued process $\bX^\eps = (\bX^\eps_t)_{t \in [0,T]}$ solution of the ordinary differential equation (ODE)\footnote{Throughout the paper,
symbols of scalar quantities are printed in italic type,
symbols of vectors are printed in bold italic type, and symbols of matrices are printed in bold type.}
\begin{equation}
\label{eq:sde1}
\frac{{\rm d}\bX^\eps}{{\rm d}t} = \bb(\bX^\eps) +\frac{1}{\eps}  \bsigma(\bX^\eps)\boeta^\eps  , \quad \bX_0^\eps=\bx_0,
\end{equation}
where $\bb(\bx)$ is a Lipschitz function from $\RR^n$ to $\RR^n$,
$\bsigma(\bx)$ is a function of class ${\cal C}^2$ with bounded derivatives from $\RR^n$ to ${\cal M}_{n,d}(\RR)$, 
and  $\boeta^\eps$ is a $\RR^d$-valued rapidly varying mean-reverting process, with a mean equal to zero, a unique invariant distribution, and a mean 
reversion time of the order of $\eps^2$.
More exactly, in this paper we address the case when $\boeta^\eps$ is a multivariate $d$-dimensional Ornstein-Uhlenbeck process
\begin{equation}
\label{eq:sde2}
 {\rm d}\boeta^\eps = \frac{{\bf K}}{\eps} {\rm d}\bW_t-\frac{{\bf A}}{\eps^2} \boeta^\eps {\rm d}t  ,
\end{equation}
where ${\bf A}$ is a $d\times d$ matrix, whose eigenvalues have positive real parts, 
${\bf K}$ is a $d\times {d'}$ matrix,
and $\bW$ is a ${d'}$-dimensional Brownian motion.
This model is classical. It can be encountered in earthquake engineering \cite{lin}
and also in finance \cite{vasicek}. It can model stationary Gaussian processes with very general spectra (see Section \ref{sec:adif1}).

Our main motivation is to estimate a quantity of the form
\begin{equation}
\label{def:Ieps}
I^\eps \triangleq \EE[ F(\bX^\eps)]
\end{equation}

for a fixed, small or moderate, parameter $\eps$, for a smooth real-valued function 
$F$ defined on the space of continuous functions over $[0,T]$.
We may think at $F(\bX)= f(\bX_T)$ where $f$
is smooth with  polynomial growth,
or $F(\bX)= \int_0^T h(\bX_s) {\rm d} s + f(\bX_T)$.
By the Feynman-Kac formula it is possible to get the value of $I^\eps$ for the model (\ref{eq:sde1}-\ref{eq:sde2}) by solving a parabolic equation, but this equation is formulated in a $d+n$-dimensional space
and it possesses large terms (of order $\eps^{-2}$) that give rapid fluctuations.
These rapid fluctuations need to be resolved by
the numerical scheme, which imposes to take a time step smaller than $\eps^2$.
The numerical resolution (with a finite difference method) is, therefore, challenging, if not impossible, and we look for other resolution methods.
It is also possible to estimate $I^\eps$ by a brute force Monte Carlo method. The Monte Carlo method, however, requires many simulations
to get an accurate estimation, and each simulation requires to resolve the  rapid fluctuations at the scale $\eps^2$, so we would like to propose an efficient variance reduction method.
The main idea is to find a limiting process $\bU$ that approximates $\bX^\eps$ in a strong sense 
when $\eps \to 0$ and for which the value 
\begin{equation}
I^0 = \EE[ F(\bU)]
\end{equation}
is known or can be estimated efficiently.
It is then possible to propose a control variate method to estimate $I^\eps$ for a fixed $\eps$.

We consider the limiting $\RR^n$-valued process $\bU$ solution of the stochastic differential equation (SDE)
\begin{equation}
\label{eq:sde3}
{\rm d}\bU = \widetilde{\bb}(\bU) {\rm d}t + \bGamma(\bU) {\rm d}\bW_t ,
\end{equation}
where $\bU$ share the same driving Brownian motion as $\boeta$, with the functions $\widetilde{\bb}(\bx)$ from $\RR^n$ to $\RR^n$ and 
$\bGamma(\bx)$ from $\RR^n$ to ${\cal M}_{n,{d'}}(\RR)$ given by
\begin{align}
\label{def:widetildeb}
\widetilde{b}_j(\bx)
&\triangleq
{b}_j(\bx) + \sum_{i=1}^n \big( (\partial_{x_i} \bsigma(\bx) ){\bf A}^{-1} {\bf C} \bsigma(\bx)^T \big)_{ji} ,\\
\bGamma(\bx)&\triangleq \bsigma(\bx) {\bf A}^{-1} {\bf K} ,
\label{def:bGamma}
\end{align}
and ${\bf C}$ is the $d \times d$ matrix defined by
\begin{equation}
\label{def:C}
{\bf C} \triangleq \int_0^\infty e^{-{\bf A} s }  {\bf K} {\bf K}^T e^{ -{\bf A}^T s} {\rm d} s.
\end{equation}
The matrix ${\bf C}$ is the covariance matrix of the 
stationary distribution of the process $\boeta^\eps$.
We show in Proposition \ref{prop:0} that the continuous process $(\bX^\eps-\bU)$ converges in probability to zero
as $\eps\to 0$.
The fact that the continuous process $\bX^\eps$ converges in distribution to $\bU$ is well-known \cite[Chapter 6]{book}, but here we get a
stronger result with a particular coupling between the two processes $\bX^\eps$ and $\bU$,
 that is needed to implement the control variate method that we have in mind.

The form of the limiting equation (\ref{eq:sde3}) is not surprising.
Indeed, by (\ref{eq:sde2}), we can anticipate that $\frac{1}{\eps} \boeta^\eps {\rm d}t \simeq  {\bf A}^{-1}{\bf K} {\rm d} \bW_t +$corrections, which explains the form (\ref{def:bGamma}) of the diffusion $\bGamma$. 
The form (\ref{def:widetildeb}) of the drift $\widetilde{\itbf b}$ is a manifestation of the It\^o-versus-Stratonovich problem \cite{pavliotis07}.
This problem is whether one should interpret the stochastic integral in the limiting equation  in It\^o sense, Stratonovich sense, or another sense.
The Wong-Zakai theory \cite{wongzakai} claims that the limiting diffusion
should be a Stratonovich equation when $d=1$. 
Indeed, Eqs.~(\ref{def:widetildeb}-\ref{def:C}) then reduce to
$\bGamma(\bx) = \frac{1}{A} \bsigma(\bx) {\bf K}$, 
${\bf C}= \frac{1}{2A} {\bf K}{\bf K}^T$,
 $$\widetilde{b}_j(\bx)-{b}_j(\bx) = 
\frac{1}{2A^2}  
 \sum_{i=1}^n \big( \partial_{x_i}( \bsigma(\bx) {\bf K} ) (\bsigma(\bx) {\bf K} )^T \big)_{ji}
=
\frac{1}{2}  
 \sum_{i=1}^n \big( \partial_{x_i} \bGamma(\bx) \bGamma(\bx)^T \big)_{ji}
 ,
 $$
 so that (\ref{eq:sde3}) can be written as
\begin{equation}
\label{eq:sde3s}
{\rm d}\bU = {\bb}(\bU) {\rm d}t + \bGamma(\bU) \circ {\rm d}\bW_t ,
\end{equation}
where $\circ$ stands for the Stratonovich integral,
because
\begin{align*}
\big( \bGamma(\bU) \circ {\rm d}\bW_t  \big)_j &=
\big( \bGamma(\bU)  {\rm d}\bW_t \big)_j +\frac{1}{2} 
\sum_{i=1}^n \sum_{j'=1}^{d'} \partial_{x_{i}} \Gamma_{jj'} (\bU) 
{\rm d} \left<  X_{i}^0 , W_{j'}\right>_t \\
&=
\big( \bGamma(\bU)  {\rm d}\bW_t \big)_j +\frac{1}{2} \sum_{i=1}^n \big( \partial_{x_{i} }  \bGamma (\bU)  \bGamma (\bU)^T\big)_{ji}  {\rm d} t .
\end{align*}
The form (\ref{eq:sde3s}) is valid when $d=1$ and it looks simpler than (\ref{eq:sde3}), 
but we have chosen to write the stochastic integral in (\ref{eq:sde3})
in It\^o's sense and to add the appropriate It\^o-Stratonovich drift correction $\widetilde{\itbf b}-{\itbf b}$, because it is a natural starting point for numerical schemes \cite{kloeden} and it is the appropriate form to express
the martingale problems used in the proofs (see Appendix). 
When $d > 1$ the difference between $\widetilde{\itbf b}$ and ${\itbf b}$ is an It\^o-Stratonovich correction that is more complex and the limiting equation (\ref{eq:sde3}) cannot be reduced to (\ref{eq:sde3s}). 

We can now introduce the Monte Carlo method for the estimation of $I^\eps$.
Let $\bW^{k}$, $k=1,\ldots,N$, be
$N$ independent and identically distributed ${d'}$-dimensional Brownian motions.
We consider
\textcolor{black}{three} Monte Carlo-type estimators of $I^\eps$:\\
1) The brute force Monte Carlo estimator is 
\begin{equation}
\label{def:JepsN}
\hat{I}^{\eps}_N \triangleq \frac{1}{N} \sum_{k=1}^N F( \bX^\eps(\bW^{k})),
\end{equation}
where $\bX^\eps(\bW^k)$ is the solution of (\ref{eq:sde1}-\ref{eq:sde2}) with $\bW^k$. The estimator $\hat{I}^{\eps}_N$ is unbiased and its variance is
\begin{equation}
\label{eq:expandJepsN}
{\rm Var}(\hat{I}^{\eps}_N ) =  \frac{1}{N} {\rm Var}(F( \bX^\eps)) .
\end{equation}
\textcolor{black}{
It is asymptotically normal as $N \to +\infty$:
\begin{equation}
\sqrt{N}\big(\hat{I}^{\eps}_N -I^\eps \big) \stackrel{dist.}{\longrightarrow} {\cal N}\big( 0 ,\sigma_{I^\eps}^2\big),
\end{equation}
with the asymptotic variance
\begin{equation}
\label{def:sigmaJeps}
\sigma_{I^\eps}^2 = {\rm Var}(F( \bX^\eps)) ,
\end{equation}
which has the following behavior as $\eps \to 0$ when $F$ is continuous and bounded (because $\bX^\eps$ weakly converges to $\bX^0$):
\begin{equation}
\sigma_{I^\eps}^2 = {\rm Var}(F( \bU )) +o(1) .
\end{equation}
}

2) The control variate estimator \cite{glasserman} is 
\begin{equation}
\label{def:IepsN}
\hat{J}^{\eps}_N \triangleq I^0 + \frac{1}{N} \sum_{k=1}^N 
F( \bX^\eps(\bW^{k})) -F( \bU(\bW^{k})), 
\end{equation}
where  $I^0=\EE[ F(\bU)]$
is supposed to be known exactly (or with high accuracy).
\textcolor{black}{The value $I^0$ can be obtained by solving a Kolmogorov equation in a $n$-dimensional framework and without large term;
if this is not possible (because $n$ is too large for instance), 
then the value $I^0$ can be obtained  by a brute force Monte Carlo method
which is easier than for $I^\eps$ because there is no large term
of order $\eps^{-2}$, so that a standard Euler scheme for stochastic differential equations can be used \cite{kloeden}}.
\textcolor{black}{
The control variate estimator $\hat{J}^{\eps}_N$ is unbiased and its variance is
\begin{equation}
{\rm Var}(\hat{J}^{\eps}_N ) =  \frac{1}{N} {\rm Var}(F( \bX^\eps) - F(\bU) )  .
\label{eq:VarIpesN}
\end{equation}
It is asymptotically normal as $N \to +\infty$:
\begin{equation}
\sqrt{N}\big(\hat{J}^{\eps}_N -I^\eps \big) \stackrel{dist.}{\longrightarrow} {\cal N}\big( 0 ,\sigma_{J^\eps}^2\big),
\end{equation}
with the asymptotic variance
\begin{equation}
\label{def:sigmaIeps}
\sigma_{J^\eps}^2 ={\rm Var}(F( \bX^\eps) - F(\bU) ).
\end{equation}
When $F$ is continuous and bounded, we have by Proposition \ref{prop:0} that $\sigma_{J^\eps}^2$ goes to zero as $\eps \to 0$.
More quantitatively, 
if $F(\bX)=f(\bX_T)$ for a smooth $f$ with bounded derivatives, then 
the asymptotic variance has the following behavior as $\eps \to 0$ (by Lemma \ref{prop:3}):
\begin{equation}
\label{eq:boundepsilon2}
 \sigma_{J^\eps}^2 \leq C \eps^2 .
\end{equation}
The order of magnitude $\eps^2$ of the asymptotic variance of $\hat{J}^\eps_N$ is confirmed by the numerical simulations that we report in Section \ref{sec:num}.
}

\noindent
\textcolor{black}{
3) The theoretical optimal control variate estimator is
\begin{equation}
\hat{O}^{\eps}_N \triangleq 
\rho^\eps I^0 + \frac{1}{N} \sum_{k=1}^N F( \bX^\eps(\bW^{k})) - \rho^\eps F( \bU(\bW^{k})),
\end{equation}
with  
\begin{equation}
\rho^\eps = {\rm Cov}( F( \bX^\eps), F(\bU))/{\rm Var}(F(\bU)).
\end{equation}
This estimator is unbiased and has the minimal variance
\begin{equation}
\label{eq:varOepsN}
{\rm Var} \big( \hat{O}^{\eps}_N  \big) =  \frac{1}{N} {\rm Var}(F( \bX^\eps) - \rho^\eps F(\bU) ) 
 ,
\end{equation}
amongst all control variate estimators of the form
$$
\rho I^0+ \frac{1}{N} \sum_{k=1}^N F( \bX^\eps(\bW^{k})) - \rho F( \bU(\bW^{k})) .
$$
Note that $\rho=0$ corresponds to the brute force Monte Carlo estimator $\hat{I}^\eps_N$, $\rho=1$ 
corresponds to the control variate estimator $\hat{J}_N^\eps$, and $\rho=\rho^\eps$ corresponds to the optimal control variate estimator $\hat{O}_N^\eps$.
The estimator $\hat{O}^\eps_N$ is asymptotically normal as $N \to +\infty$:
\begin{equation}
\sqrt{N}\big(\hat{O}^{\eps}_N -I^\eps \big) \stackrel{dist.}{\longrightarrow} {\cal N}\big( 0 ,\sigma_{O^\eps}^2\big),
\end{equation}
with the asymptotic variance
\begin{equation}
\sigma_{O^\eps}^2 ={\rm Var}(F( \bX^\eps) - \rho^\eps F(\bU) ).
\end{equation}
The estimator $\hat{O}^\eps_N$ is, however, not practical as it depends on $\rho^\eps$ which is unknown.
The practical optimal control variate estimator \cite{glasserman} is
\begin{equation}
\label{def:KepsN}
\hat{K}^{\eps}_N \triangleq  \hat{\rho}^\eps_N I^0+ \frac{1}{N} \sum_{k=1}^N F( \bX^\eps(\bW^{k})) - \hat\rho^\eps_N F( \bU(\bW^{k})), 
\end{equation}
where $\hat\rho^\eps_N$ is the empirical correlation 
\begin{equation}
\hat\rho^\eps_N = \frac{
 \sum_{k=1}^N ( F( \bX^\eps(\bW^{k})) - \hat{I}^\eps_N) (F( \bU(\bW^{k})) - \hat{I}^0_N)
 }
 {
  \sum_{k=1}^N (F( \bU(\bW^{k})) - \hat{I}^0_N )^2
},
\label{def:hatrhoNeps}
\end{equation}
with $ \hat{I}^\eps_N= \frac{1}{N}  \sum_{k=1}^N F( \bX^\eps(\bW^{k}))$ as in (\ref{def:JepsN}) 
and $\hat{I}^0_N = \frac{1}{N}  \sum_{k=1}^N F( \bU(\bW^{k})) $.
This estimator is a practical and approximate version of the theoretical optimal control variate estimator 
$\hat{O}^{\eps}_N$ in which the unknown correlation coefficient $\rho^\eps$ has been replaced by its empirical estimator
$\hat{\rho}_N^\eps$.
The estimator $\hat{K}^\eps_N$ may be slightly biased 
and may have a variance slightly larger than (\ref{eq:varOepsN}) because of the empirical estimation of $\rho^\eps$.
$ \hat{K}^{\eps}_N $ is, however, asymptotically normal with
an asymptotic variance that is the same one as that of the optimal estimator $ \hat{O}^{\eps}_N $,
as shown by the following proposition.
\begin{prop}
\label{prop:21}
As $N\to +\infty$,
\begin{equation}
\label{eq:asynorKeps}
\sqrt{N} \big(  \hat{K}^{\eps}_N - I^\eps\big) \stackrel{dist.}{\longrightarrow}
{\cal N}\big(0, \sigma_{K^\eps}^2\big),
\end{equation}
with
\begin{equation}
\label{def:sigmaKeps}
\sigma_{K^\eps}^2 =\sigma_{O^\eps}^2 =
 {\rm Var}(F( \bX^\eps) - \rho^\eps F(\bU) )  .
\end{equation}
Furthermore, 
if $F$ is continuous and bounded, then  $\sigma_{K^\eps}^2$ goes to zero as $\eps \to 0$.
If $F(\bX)=f(\bX_T)$, with $f$ with bounded derivatives, then there exists $C>0$ such that
\begin{equation}
\label{eq:estimsigmaKeps}
 \sigma^2_{K^\eps} \leq C \eps^2 ,\quad \quad 
0\leq \sigma^2_{J^\eps} -  \sigma^2_{K^\eps} \leq C \eps^4.
 \end{equation}
\end{prop}
\noindent
{\it Proof.}
By the law of large numbers, $\hat\rho^\eps_N$ converges to $\rho^\eps$ as $N \to +\infty$.
The convergence holds almost surely, hence in probability.
We have
$$
 \hat{K}^{\eps}_N - I^\eps= ( \hat{O}^\eps_N-I^\eps)  - (\hat{\rho}_N^\eps-\rho^\eps) (\hat{I}_N^0-I^0) ,
$$
so we get (\ref{eq:asynorKeps}-\ref{def:sigmaKeps}) from Slutsky's theorem.\\
Furthermore, we
have 
$$
\sigma^2_{J^\eps} - \sigma^2_{K^\eps}  = {\rm Var}( F(\bX^0)) (1 - \rho^\eps)^2 .
$$
If $F$ is continuous and bounded, then  $\rho^\eps$ goes to one 
and $\sigma^2_{J^\eps}$ goes to zero as $\eps \to 0$ by Proposition \ref{prop:0}.
If $F(\bX)=f(\bX_T)$, then, by Lemma \ref{prop:3},
$1-\rho^\eps$ and  $\sigma^2_{J^\eps}$  are of order $O(\eps^2)$ for small $\eps$.
This shows the desired result (\ref{eq:estimsigmaKeps}).
\qed
}

Proposition \ref{prop:21} shows that the asymptotic variances of the estimators $\hat{K}^\eps_N$ and $\hat{J}^\eps_N$ are equivalent
for vanishingly small $\eps$ and of the order of $O(\eps^2)$, and that  the asymptotic variance of the estimator $\hat{K}^\eps_N$
is slightly smaller than that of $\hat{J}^\eps_N$ for moderately small $\eps$.
These statements are confirmed by the numerical simulations that we report in Section \ref{sec:num}.

In addition, motivated by the examples that we address in Section \ref{sec:nummulti},
we consider the case where the $\RR^n$-valued process $\bX^\eps$ satisfies a multivalued ODE of the form 
\begin{equation}
\label{eq:mode1}
\frac{{\rm d}\bX^\eps}{{\rm d}t} + \partial \varphi (\bX^\eps) \ni \bb(\bX^\eps) +\frac{1}{\eps}  \bsigma \boeta^\eps  , \quad \bX_0^\eps=\bx_0,
\end{equation}
and the case where $\bX^\eps$ together with a $\RR^m$-valued process $\bZ^\eps$ satisfy the multivalued ODE
\begin{equation}
\label{eq:mode2}
\begin{dcases}
\frac{{\rm d}\bX^\eps}{{\rm d}t} + \partial \varphi (\bX^\eps) \ni \bb^X(\bX^\eps,\bZ^\eps) +\frac{1}{\eps}  \bsigma \boeta^\eps, \quad \bX_0^\eps=\bx_0 ,
\\
\frac{{\rm d}\bZ^\eps}{{\rm d}t} + \partial \psi (\bZ^\eps) \ni \bb^Z(\bX^\eps,\bZ^\eps), \quad \bZ_0^\eps=\bz_0
.
\end{dcases}
\end{equation}
Here $\bsigma \in {\cal M}_{n,d}(\RR)$ is constant, 
$\bb(\bx)$ from $\RR^n$ to $\RR^n$, 
$\bb^Z(\bx,\bz)$ from $\RR^{n+m}$ to $\RR^m$ and $\bb^X(\bx,\bz)$ from $\RR^{n+m}$ to $\RR^n$ are Lipschitz functions. 
The operators $\partial \varphi$ and $\partial \psi$ are the subdifferentials of some lower semi continuous (l.s.c.) convex functions 
$\varphi$ from $\RR^n$ to $[0,+\infty]$ and $\psi$ from $\RR^m$ to $[0,+\infty]$.
Stronger hypotheses will be assumed on $\varphi$ compared to $\psi$ as explained in Section \ref{sec:adif2}
and important examples motivate the two situations as shown in Section \ref{sec:nummulti}.
\textcolor{black}{It is important to observe that the multivalued operators that appear in the differential inclusions above are subdifferential of convex functions, therefore existence and uniqueness are guaranteed \cite[page 72]{MR0348562}. For the reader's convenience, proofs of existence and uniqueness are given in Appendix \ref{app:A}.
It is worth mentioning that there is an alternative formulation using the language of variational inequalities, that is equivalent to differential inclusions. 
Eq.~(\ref{eq:mode1}) is equivalent to
$$
\forall \boldsymbol{\xi} \in \RR^n, \: \forall t >0, \:  
\Big( \bb(\bX^\eps) +\frac{1}{\eps}  \bsigma \boeta^\eps - \frac{{\rm d}\bX^\eps}{{\rm d}t} \Big) \cdot \big( \boldsymbol{\xi} - \bX^\eps \big)
+ \varphi(\bX^\eps) \leq \varphi(\boldsymbol{\xi}) , \quad \bX_0^\eps=\bx_0 ,
$$
and Eq.~(\ref{eq:mode2}) is equivalent to
\begin{align*}
& \forall \boldsymbol{\xi} \in \RR^n, \: \forall \boldsymbol{\zeta} \in \RR^m, \: \forall t >0, \:\\  
& \Big( \bb^X(\bX^\eps,\bZ^\eps) +\frac{1}{\eps}  \bsigma \boeta^\eps - \frac{{\rm d}\bX^\eps}{{\rm d}t} \Big) \cdot \big( \boldsymbol{\xi} - \bX^\eps \big)
+ \varphi(\bX^\eps) \leq \varphi(\boldsymbol{\xi}) 
, \quad \bX_0^\eps=\bx_0,\\
& \Big( \bb^Z(\bX^\eps,\bZ^\eps) - \frac{{\rm d}\bZ^\eps}{{\rm d}t}\Big) \cdot \big( \boldsymbol{\zeta} - \bZ^\eps \big)
+ \psi(\bZ^\eps) \leq \psi(\boldsymbol{\zeta}) 
, \quad \bZ_0^\eps=\bz_0.
\end{align*}}

Propositions \ref{prop:1} and \ref{prop:2} show that the multi-valued process $\bX^\eps$ strongly converges to a limiting process
solution of a multivalued SDE.
Eqs.~(\ref{eq:prop5:1}) and (\ref{eq:prop5:2}) show that
the control variate estimators have asymptotic variances of order $\eps^2$ for (\ref{eq:mode1}) and $\eps$ for (\ref{eq:mode2}).

\begin{figure}[!h]
\center
\begin{tikzpicture}[scale=1.5]
\draw[color=gray] (-1.125,0.875) -- (-1,1);
\draw[color=gray] (-1.125,0.75) -- (-1,0.875);
\draw[color=gray] (-1.125,0.625) -- (-1,0.75);
\draw[color=gray] (-1.125,0.5) -- (-1,0.625);
\draw[color=gray] (-1.125,0.375) -- (-1,0.5);
\draw[color=gray] (-1.125,0.25) -- (-1,0.375);
\draw[color=gray] (-1.125,0.125) -- (-1,0.25);
\draw[color=gray] (-1.125,0) -- (-1,0.125);
\draw[color=gray] (-1.125,0.125) -- (-1,0.25);
\draw[color=gray] (-1.125,-0.125) -- (-1,0);
\draw[color=gray] (-1.125,-0.25) -- (-1,-0.125);
\draw[color=gray] (-1.125,-0.375) -- (-1,-0.25);
\draw[color=gray] (-1.125,-0.5) -- (-1,-0.375);
\draw[color=gray] (-1.125,-0.625) -- (-1,-0.5);
\draw[color=gray] (-1.125,-0.75) -- (-1,-0.625);
\draw[color=gray] (-1.125,-0.875) -- (-1,-0.75);
\draw[color=gray] (-1.125,-1) -- (-1,-0.875);
\draw[thick] (-1,-1) -- (-1,1) ;
\draw[color=gray] (-1.125,-1.125) -- (-1,-1);
\draw[color=gray] (-1,-1.125) -- (-0.875,-1);
\draw[color=gray] (-0.875,-1.125) -- (-0.750,-1);
\draw[color=gray] (-0.75,-1.125) -- (-0.625,-1);
\draw[color=gray] (-0.625,-1.125) -- (-0.5,-1);
\draw[color=gray] (-0.5,-1.125) -- (-0.375,-1);
\draw[color=gray] (-0.375,-1.125) -- (-0.250,-1);
\draw[color=gray] (-0.250,-1.125) -- (-0.125,-1);
\draw[color=gray] (-0.125,-1.125) -- (0,-1);
\draw[color=gray] (0,-1.125) -- (0.125,-1);
\draw[color=gray] (0.125,-1.125) -- (0.25,-1);
\draw[color=gray] (0.250,-1.125) -- (0.375,-1);
\draw[color=gray] (0.375,-1.125) -- (0.5,-1);
\draw[color=gray] (0.5,-1.125) -- (0.625,-1);
\draw[color=gray] (0.625,-1.125) -- (0.75,-1);
\draw[color=gray] (0.75,-1.125) -- (0.875,-1);
\draw[color=gray] (0.875,-1.125) -- (1,-1);
\draw[color=gray] (1,-1.125) -- (1.125,-1);
\draw[color=gray] (1.125,-1.125) -- (1.25,-1);
\draw[color=gray] (1.250,-1.125) -- (1.375,-1);
\draw[color=gray] (1.375,-1.125) -- (1.5,-1);
\draw[color=gray] (1.5,-1.125) -- (1.625,-1);
\draw[color=gray] (1.625,-1.125) -- (1.75,-1);
\draw[color=gray] (1.75,-1.125) -- (1.875,-1);
\draw[color=gray] (1.875,-1.125) -- (2,-1);
\draw[color=gray] (2,-1.125) -- (2.125,-1);
\draw[color=gray] (2.125,-1.125) -- (2.25,-1);
\draw[color=gray] (2.250,-1.125) -- (2.375,-1);
\draw[color=gray] (2.375,-1.125) -- (2.5,-1);
\draw[color=gray] (2.5,-1.125) -- (2.625,-1);
\draw[color=gray] (2.625,-1.125) -- (2.75,-1);
\draw[color=gray] (2.75,-1.125) -- (2.875,-1);
\draw[color=gray] (2.875,-1.125) -- (3,-1);

\draw[color=gray] (3,-1.125) -- (3.125,-1);
\draw[color=gray] (3.125,-1.125) -- (3.25,-1);
\draw[color=gray] (3.250,-1.125) -- (3.375,-1);
\draw[color=gray] (3.375,-1.125) -- (3.5,-1);
\draw[color=gray] (3.5,-1.125) -- (3.625,-1);
\draw[color=gray] (3.625,-1.125) -- (3.75,-1);
\draw[color=gray] (3.75,-1.125) -- (3.875,-1);
\draw[color=gray] (3.875,-1.125) -- (4,-1);

\draw[color=gray] (4,-1.125) -- (4.125,-1);
\draw[color=gray] (4.125,-1.125) -- (4.25,-1);
\draw[color=gray] (4.250,-1.125) -- (4.375,-1);
\draw[color=gray] (4.375,-1.125) -- (4.5,-1);
\draw[color=gray] (4.5,-1.125) -- (4.625,-1);
\draw[color=gray] (4.625,-1.125) -- (4.75,-1);
\draw[color=gray] (4.75,-1.125) -- (4.875,-1);
\draw[color=gray] (4.875,-1.125) -- (5,-1);

\draw[color=gray] (5,-1.125) -- (5.125,-1);
\draw[color=gray] (5.125,-1.125) -- (5.25,-1);
\draw[color=gray] (5.250,-1.125) -- (5.375,-1);
\draw[color=gray] (5.375,-1.125) -- (5.5,-1);
\draw[color=gray] (5.5,-1.125) -- (5.625,-1);
\draw[color=gray] (5.625,-1.125) -- (5.75,-1);
\draw[color=gray] (5.75,-1.125) -- (5.875,-1);
\draw[color=gray] (5.875,-1.125) -- (6,-1);

\draw[color=gray] (6,-1.125) -- (6.125,-1);
\draw[color=gray] (6.125,-1.125) -- (6.25,-1);
\draw[color=gray] (6.250,-1.125) -- (6.375,-1);
\draw[color=gray] (6.375,-1.125) -- (6.5,-1);
\draw[color=gray] (6.5,-1.125) -- (6.625,-1);
\draw[color=gray] (6.625,-1.125) -- (6.75,-1);
\draw[color=gray] (6.75,-1.125) -- (6.875,-1);
\draw[color=gray] (6.875,-1.125) -- (7,-1);

\draw[thick] (-1,-1) -- (7,-1) ;
\draw[thick] (-1,0.5) -- (0.125,0.5) ;
\draw[thick] (0.25,0.5) -- (1,0.5) ; 
\draw[thick] (-0.25,0.625) -- (0.25,0.625) ;
\draw[thick] (-0.25,0.375) -- (0.25,0.375) ;
\draw[thick] (0.25,0.375) -- (0.25,0.625) ;
\draw[thick] (0.125,0.4) -- (0.125,0.6) ;

\draw[thick] (-1,-0.25) -- (-0.25,-0.25) ;
\draw[thick] (0.25,-0.25) -- (1,-0.25) ;
\filldraw (0,-0.25) circle (5pt);
\draw[thick] (0,-0.25) node[] {$\white{1}$};

\draw[thick] (3,0.5) -- (3.75,0.5) ;
\draw[thick] (4.25,0.5) -- (5,0.5) ;
\filldraw (4,0.5) circle (5pt);
\draw[thick] (4,0.5) node[] {$\white{3}$};

\draw[thick] (3,-.25) -- (3.5,-.25) ;
\draw[snake=coil,segment length=4pt, thick] (3.5,-.25) -- (4.5,-.25) ; 
\draw[thick] (4.5,-.25) -- (5,-.25) ;

\draw[thick] (1,-0.75) -- (1,0.75);
\draw[thick] (3,-0.75) -- (3,0.75);
\draw[thick] (1,-0.75) -- (1,0.75);
\draw[thick] (3,-0.75) -- (3,0.75);
\draw[thick] (1,-0.75) -- (3,-0.75);
\draw[thick] (1,0.75) -- (3,0.75);
\draw[thick] (2,0) node[] {$m_1$};

\draw[thick] (2,1.75) node[] {forcing $\epsilon^{-1} \eta^\epsilon$};
\draw[thick,->]  (1.5,1.53) -- (2.5,1.53);

\draw (1.5,-0.875) circle (0.125cm);
\draw (2.5,-0.875) circle (0.125cm);
\draw[thick] (5,-0.75) -- (5,0.75);
\draw[thick] (7,-0.75) -- (7,0.75);
\draw[thick] (5,-0.75) -- (5,0.75);
\draw[thick] (7,-0.75) -- (7,0.75);
\draw[thick] (5,-0.75) -- (7,-0.75);
\draw[thick] (5,0.75) -- (7,0.75);
\draw[thick] (6,0) node[] {$m_3$};

\draw (5.5,-0.875) circle (0.125cm);
\draw (6.5,-0.875) circle (0.125cm);

\draw[thick] (2,0.75) -- (2,1.25);
\draw[->]  (2,1.25) -- (2.5,1.25);
\draw[thick] (2.5,1.25) node[right] {$X_1$};

\draw[thick] (6,0.75) -- (6,1.25);
\draw[->]  (6,1.25) -- (6.5,1.25);
\draw[thick] (6.5,1.25) node[right] {$X_3$};
\draw[thick] (0.15,0.85) node[] {$c_1$};
\draw[thick] (4,-0.6) node[] {$k_3$};
\end{tikzpicture}
\caption{A rheological model of a two-degree of freedom system. Two masses $m_1$ and $m_3$ are associated in series with elements which are themselves an association of dampers and springs. Elements \textcircled{1} and \textcircled{3} represent a spring and a damper respectively, both possibly nonlinear or hysteretic. Here $c_1$ is a damping coefficient associated to the linear damper connecting the mass $m_1$ to the foundation and $k_3$ is a stiffness coefficient of the linear spring linking the masses $m_1$ and $m_3$. A random forcing $\epsilon^{-1} \eta^\epsilon$ is applied to the mass $m_1$ (e.g. wind forces on a two-storey building).}
\label{fig:tdof}
\end{figure}
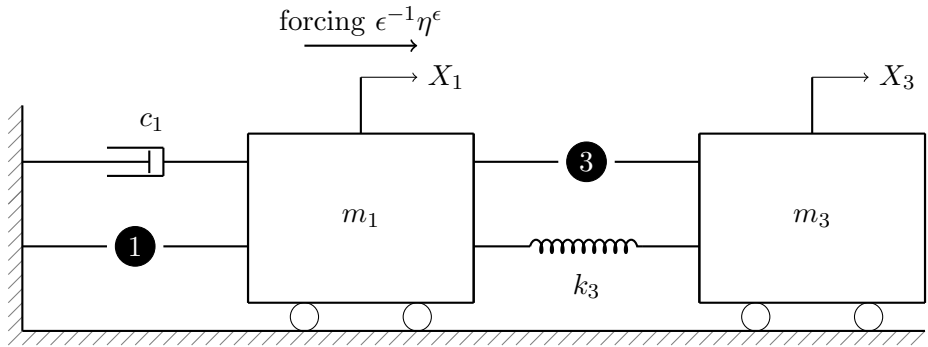

To demonstrate the efficiency of our method on a practical problem, we consider a two-degree of freedom (TDOF) system as shown in Figure \ref{fig:tdof}. It can describe a broad class of TDOF structures, 
including a two-storey building as presented in \cite[Figure 4.6(a)]{spanosbook}.

When the external force $\eta^\eps$ is a colored noise 
such as an Ornstein-Uhlenbeck process, the equation of motion can be written in the form of Equation \eqref{eq:sde1} with $n=4$,
where $(X^\eps_1,X_2^\eps)$, resp. $(X^\eps_3,X^\eps_4)$, represents the position and the velocity of the mass $m_1$, resp. $m_3$, shown in Figure~\ref{fig:tdof}. Many nonlinear behaviours enter into this framework, we have in mind a nonlinear spring of the linear-plus-quadratic cubic type and a nonlinear damper of the linear-plus-quadratic type (see Example \ref{example3c} and Figure~\ref{fig:i2}).
Similarly Equations \eqref{eq:mode1} and \eqref{eq:mode2} arise in the description of nonlinear behaviours with hysteresis such as elasto-plasticity and friction, see  \cite[Chapter 8]{spanosbook} and Section \ref{sec:nummulti} (see Example~\ref{example81c} and Figure~\ref{fig:i3}). 
\textcolor{black}{In Figures~\ref{fig:i2}-\ref{fig:i3} we compare the behaviors of the brute force Monte Carlo estimator $\hat{I}_N^\eps$  with the ones of the control variate estimators $\hat{J}^\eps_N$ and $\hat{K}^\eps_N$. 
We also plot the empirical estimators of the asymptotic variances of the estimators $\hat{I}_N^\eps$, $\hat{J}_N^\eps$, and $\hat{K}^\eps_N$
as described in Subsection \ref{subsec:num3multi}.
In addition, for each value of $\varepsilon \in \{0.1,0.5,0.9\}$, error bars (95\% confidence interval) are shown for each of the estimators (in $\hat{I}_N^\eps$, $\hat{J}_N^\eps$, $\hat{K}^\eps_N$ order from the left to the right).
Here $I^0$ is obtained by a massive Monte Carlo estimation of the limit process, which is possible with a coarse grid step as there is no large term involved.
We can observe that the control variate estimator $\hat{K}^\eps_N$ has always the minimal variance. When $\eps$ is small and the original system and the limit system are poorly correlated $\rho^\eps\simeq 0$, it behaves as the standard Monte Carlo estimator $\hat{I}^\eps_N$.
When $\eps$ is small and the original system and the limit system are strongly correlated $\rho^\eps\simeq 1$,
 it behaves as the control variate estimator estimator $\hat{J}^\eps_N$.
We can also observe that the variance reduction is by a factor of order $\eps^2$ when the quantity to be estimated is the expectation of a smooth function, 
while it is of order $\eps$  when the quantity to be estimated is the expectation of an indicator function.}



\begin{figure}[h!]
\centering
%
\begin{tikzpicture}[scale=0.56]
\begin{axis}[legend style={at={(0,0)},anchor=south west}, compat=1.3,
  xmin=0.01, xmax=1,ymin=0.0325,ymax=0.1025,
  xlabel= {$\varepsilon$},
  ylabel= {}]
\addplot[densely dashed,color=black,mark=none,mark size=1pt] table [x index=0, y index=5]{practical_nl-ou_dt1.0E-05_N1.0E+04_data.txt};
\addlegendentry{$\hat{I}_N^\varepsilon$}
\addplot[densely dotted,thick,color=blue,mark=none,mark size=1pt] table [x index=0, y index=6]{practical_nl-ou_dt1.0E-05_N1.0E+04_data.txt};
\addlegendentry{$\hat{J}_N^\varepsilon$}
\addplot[solid,color=red,mark=none,mark size=1pt] table [x index=0, y index=7]{practical_nl-ou_dt1.0E-05_N1.0E+04_data.txt};
\addlegendentry{$\hat{K}_N^\varepsilon$}
\addplot[thick,dotted,color=black,mark=none,mark size=4pt] table [x index=0, y index=8]{practical_nl-ou_dt1.0E-05_N1.0E+04_data.txt};
\addlegendentry{$I^0$}
\addplot[only marks,
  black, mark options={black,scale=0.25}, 
  error bars/.cd, 
    y fixed,
    y dir=both, 
    y explicit
] table [x=x, y=y,y error=error, col sep=comma] {
    x,  y,       error
0.0875,  0.0979195,  0.00270771   
0.4875,  0.0687944,  0.00185577 
0.8875,  0.0357051,  0.000973789
};
\addplot[only marks,
  blue, mark options={black,scale=0.25}, 
  error bars/.cd, 
    y fixed,
    y dir=both, 
    y explicit
] table [x=x, y=y,y error=error, col sep=comma] {
    x,  y,       error
0.1,  0.0950673,  0.000384563  
0.5,  0.0690914,  0.00179333  
0.9,  0.0381916,  0.00235276
};
\addplot[only marks,
  red, mark options={black,scale=0.25}, 
  error bars/.cd, 
    y fixed,
    y dir=both, 
    y explicit
] table [x=x, y=y,y error=error, col sep=comma] {
    x,  y,       error
0.1125,  0.0951224,  0.000380928  
0.5125,  0.0689481,  0.00130202 
0.9125,  0.0360534,  0.000907303 
};
\end{axis}
\node at (5cm,3.7cm) {\textcolor{black}{(a)}};
\end{tikzpicture}
\begin{tikzpicture}[scale=0.56]
\begin{axis}[legend style={at={(0.65,1)},anchor=north east}, compat=1.3,
  xmin=0.01, xmax=1,ymin=0,ymax=0.02,
  xlabel= {$\varepsilon$},
  ylabel= {}]
\addplot[densely  dashed, thick, color=black,mark=none,mark size=1pt] table [x index=0, y index=9]{practical_nl-ou_dt1.0E-05_N1.0E+04_data.txt};
\addlegendentry{$\hat \sigma_{I^\varepsilon,N}^2$}
\addplot[densely dotted, thick,color=blue,mark=none,mark size=1pt] table [x index=0, y index=10]{practical_nl-ou_dt1.0E-05_N1.0E+04_data.txt};
\addlegendentry{$\hat \sigma_{J^\varepsilon,N}^2$}
\addplot[solid,thick,color=red,mark=none,mark size=1pt] table [x index=0, y index=11]{practical_nl-ou_dt1.0E-05_N1.0E+04_data.txt};
\addlegendentry{$\hat \sigma_{K^\varepsilon,N}^2$}
\end{axis}
\node at (5cm,2.5cm) {\textcolor{black}{(b)}};
\end{tikzpicture}
\begin{tikzpicture}[scale=0.56]
\begin{axis}[legend style={at={(1,1)},anchor=north east}, compat=1.3,
  xmin=0.01, xmax=1,ymin=0,ymax=0.0425,
  xlabel= {$\varepsilon$},
  ylabel= {}]
\addplot[solid, thick, densely dotted,color=black,mark=none,mark size=1pt] table [x index=0, y index=3]{practical_nl-ou_dt1.0E-05_N1.0E+04_data.txt};
\addlegendentry{$\hat \sigma_{K^\varepsilon,N}^2/\varepsilon$}
\addplot[solid, thick, color=black,mark=none,mark size=1pt] table [x index=0, y index=4]{practical_nl-ou_dt1.0E-05_N1.0E+04_data.txt};
\addlegendentry{$\hat \sigma_{K^\varepsilon,N}^2/\varepsilon^2$}
\end{axis}
\node at (5cm,3.7cm) {\textcolor{black}{(c)}};
\end{tikzpicture}


\begin{tikzpicture}[scale=0.56]
\begin{axis}[legend style={at={(0,1)},anchor=north west}, compat=1.3,
  xmin=0.01, xmax=1,ymin=0.25,ymax=0.45,
  xlabel= {$\varepsilon$},
  ylabel= {}]
\addplot[densely dashed,color=black,mark=none,mark size=1pt] table [x index=0, y index=5]{proba-practical_nl-ou_dt1.0E-05_N1.0E+04_data.txt};
\addlegendentry{$\hat{I}_N^\varepsilon$}
\addplot[densely dotted,thick,color=blue,mark=none,mark size=1pt] table [x index=0, y index=6]{proba-practical_nl-ou_dt1.0E-05_N1.0E+04_data.txt};
\addlegendentry{$\hat{J}_N^\varepsilon$}
\addplot[solid,color=red,mark=none,mark size=1pt] table [x index=0, y index=7]{proba-practical_nl-ou_dt1.0E-05_N1.0E+04_data.txt};
\addlegendentry{$\hat{K}_N^\varepsilon$}
\addplot[thick,dotted,color=black,mark=none,mark size=4pt] table [x index=0, y index=8]{proba-practical_nl-ou_dt1.0E-05_N1.0E+04_data.txt};
\addlegendentry{$I^0$}
\addplot[only marks,
  black, mark options={black,scale=0.25}, 
  error bars/.cd, 
    y fixed,
    y dir=both, 
    y explicit
] table [x=x, y=y,y error=error, col sep=comma] {
    x,  y,       error
0.0875,  0.276,  0.00876153   
0.4875,  0.3165,  0.00911617
0.8875,  0.4267,  0.00969412  
};
\addplot[only marks,
  blue, mark options={black,scale=0.25}, 
  error bars/.cd, 
    y fixed,
    y dir=both, 
    y explicit
] table [x=x, y=y,y error=error, col sep=comma] {
    x,  y,       error
0.1,  0.27376,  0.00487548 
0.5,  0.31816,  0.0105489
0.9,  0.42606,  0.0122667
};
\addplot[only marks,
  red, mark options={black,scale=0.25}, 
  error bars/.cd, 
    y fixed,
    y dir=both, 
    y explicit
] table [x=x, y=y,y error=error, col sep=comma] {
    x,  y,       error
0.1125,  0.274098,  0.00469474  
0.5125,  0.317019,  0.00870337 
0.9125,  0.426619,  0.00963149 
};
\end{axis}
\node at (5cm,2.4cm) {\textcolor{black}{(a)}};
\end{tikzpicture}
\begin{tikzpicture}[scale=0.56]
\begin{axis}[legend style={at={(1,0)},anchor=south east}, compat=1.3,
  xmin=0.01, xmax=1,ymin=0,ymax=0.4125,
  xlabel= {$\varepsilon$},
  ylabel= {}]
\addplot[densely  dashed, thick, color=black,mark=none,mark size=1pt] table [x index=0, y index=9]{proba-practical_nl-ou_dt1.0E-05_N1.0E+04_data.txt};
\addlegendentry{$\hat \sigma_{I^\varepsilon,N}^2$}
\addplot[densely dotted, thick,color=blue,mark=none,mark size=1pt] table [x index=0, y index=10]{proba-practical_nl-ou_dt1.0E-05_N1.0E+04_data.txt};
\addlegendentry{$\hat \sigma_{J^\varepsilon,N}^2$}
\addplot[solid,thick,color=red,mark=none,mark size=1pt] table [x index=0, y index=11]{proba-practical_nl-ou_dt1.0E-05_N1.0E+04_data.txt};
\addlegendentry{$\hat \sigma_{K^\varepsilon,N}^2$}
\end{axis}
\node at (5cm,3.9cm) {\textcolor{black}{(b)}};
\end{tikzpicture}
\begin{tikzpicture}[scale=0.56]
\begin{axis}[legend style={at={(1,1)},anchor=north east}, compat=1.3,
  xmin=0.01, xmax=1,ymin=0,ymax=2.0,
  xlabel= {$\varepsilon$},
  ylabel= {}]
\addplot[solid, thick, densely dotted,color=black,mark=none,mark size=1pt] table [x index=0, y index=3]{proba-practical_nl-ou_dt1.0E-05_N1.0E+04_data.txt};
\addlegendentry{$\hat \sigma_{K^\varepsilon,N}^2/\varepsilon$}
\addplot[solid, thick, color=black,mark=none,mark size=1pt] table [x index=0, y index=4]{proba-practical_nl-ou_dt1.0E-05_N1.0E+04_data.txt};
\addlegendentry{$\hat \sigma_{K^\varepsilon,N}^2/\varepsilon^2$}
\end{axis}
\node at (5cm,3.7cm) {\textcolor{black}{(c)}};
\end{tikzpicture}

\caption{Example \ref{example3c} of a TDOF modeling a two-storey building with nonlinear spring of the linear-plus-quadratic cubic type and nonlinear damper of the linear-plus-quadratic type, driven by an Ornstein-Uhlenbeck noise.
In the top row the target is to estimate $I^\eps=\EE \big[ (X_{1,T}^\eps)^2+ (X_{3,T}^\eps)^2 \big]$ for $T=1$ and 
the expectation of the control variate $I^0=\EE [ (X_{1,T}^0)^2+ (X_{3,T}^0)^2  ]$ is obtained by an intensive Monte Carlo computation with coarse time step.
In the bottom row the target is to estimate $I^\eps=\PP( | X^\eps_{1,T} | \leq a, | X^\eps_{3,T} | \leq b)$ for $T=1,a=0.1,b=0.1$ and 
the expectation of the control variate $I^0=\PP(| X_{1,T}^0 | \leq a,  | X_{3,T}^0| \leq b)$ is also obtained by an intensive Monte Carlo computation with coarse time step.
\textcolor{black}{The numerical procedure is standard and is described in Section 4 (Euler-Maruyama time discretization with time step $\delta t = 10^{-5}$).}
Here the number of Monte Carlo samples is $N = 10^4$ and $m_1=m_3=c_1 = k_3 = 1$.
The details of the elements \textcircled{1} and \textcircled{3} with nonlinear behaviours can be found in Example \ref{example3c}.
}
\label{fig:i2}
\end{figure}
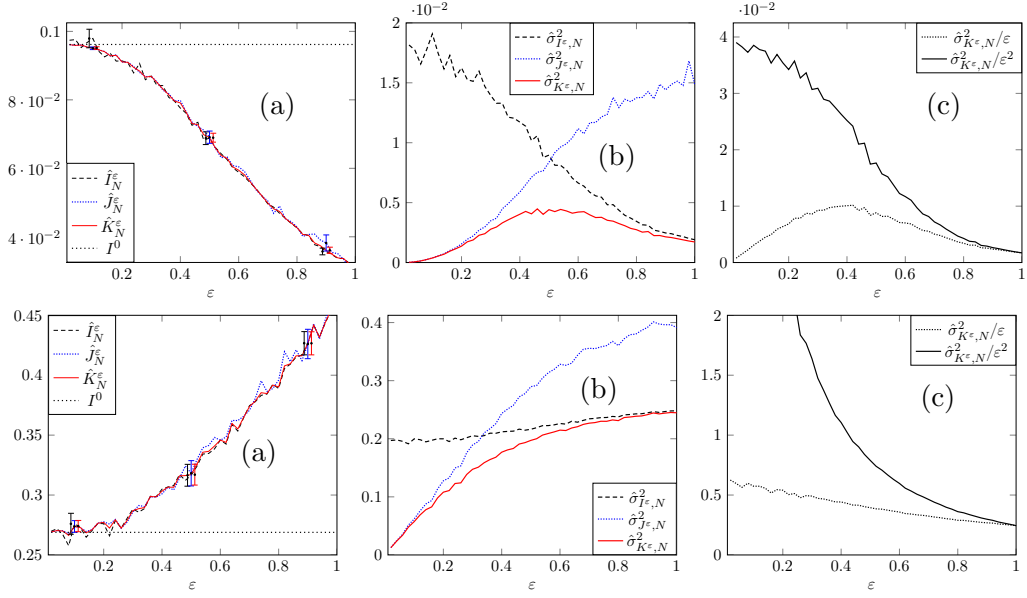


\begin{figure}[h!]
\centering
%
\begin{tikzpicture}[scale=0.56]
\begin{axis}[legend style={at={(0,0)},anchor=south west}, compat=1.3,
  xmin=0.01, xmax=1,ymin=0.025,ymax=0.06,
  xlabel= {$\varepsilon$},
  ylabel= {}]
\addplot[densely dashed,color=black,mark=none,mark size=1pt] table [x index=0, y index=5]{practical_ns-ou_dt1.0E-05_N1.0E+04_data.txt};
\addlegendentry{$\hat{I}_N^\varepsilon$}
\addplot[densely dotted,thick,color=blue,mark=none,mark size=1pt] table [x index=0, y index=6]{practical_ns-ou_dt1.0E-05_N1.0E+04_data.txt};
\addlegendentry{$\hat{J}_N^\varepsilon$}
\addplot[solid,color=red,mark=none,mark size=1pt] table [x index=0, y index=7]{practical_ns-ou_dt1.0E-05_N1.0E+04_data.txt};
\addlegendentry{$\hat{K}_N^\varepsilon$}
\addplot[thick,dotted,color=black,mark=none,mark size=4pt] table [x index=0, y index=8]{practical_ns-ou_dt1.0E-05_N1.0E+04_data.txt};
\addlegendentry{$I^0$}
error bars I
\addplot[only marks,
  black, mark options={black,scale=0.25}, 
  error bars/.cd, 
    y fixed,
    y dir=both, 
    y explicit
] table [x=x, y=y,y error=error, col sep=comma] {
    x,  y,       error
0.0875,  0.0529269,  0.000922953   
0.4875,  0.0432881,  0.000745692 
0.8875,  0.0284504,  0.000571672
};
\addplot[only marks,
  blue, mark options={black,scale=0.25}, 
  error bars/.cd, 
    y fixed,
    y dir=both, 
    y explicit
] table [x=x, y=y,y error=error, col sep=comma] {
    x,  y,       error
0.1,  0.0537329,  0.000167384  
0.5,  0.0423627,  0.000755432
0.9,  0.0285731,  0.00091615
};
\addplot[only marks,
  red, mark options={black,scale=0.25}, 
  error bars/.cd, 
    y fixed,
    y dir=both, 
    y explicit
] table [x=x, y=y,y error=error, col sep=comma] {
    x,  y,       error
0.1125,  0.0537116,  0.000165567 
0.5125,  0.0428331,  0.000586188
0.9125,  0.0284754,  0.000539526
};
\end{axis}
\node at (5cm,3.7cm) {\textcolor{black}{(a)}};
\end{tikzpicture}
\begin{tikzpicture}[scale=0.56]
\begin{axis}[legend style={at={(0.635,1)},anchor=north east}, compat=1.3,
  xmin=0.01, xmax=1,ymin=0,ymax=0.0029,
  xlabel= {$\varepsilon$},
  ylabel= {}]
\addplot[densely  dashed, thick, color=black,mark=none,mark size=1pt] table [x index=0, y index=9]{practical_ns-ou_dt1.0E-05_N1.0E+04_data.txt};
\addlegendentry{$\hat \sigma_{I^\varepsilon,N}^2$}
\addplot[densely dotted, thick,color=blue,mark=none,mark size=1pt] table [x index=0, y index=10]{practical_ns-ou_dt1.0E-05_N1.0E+04_data.txt};
\addlegendentry{$\hat \sigma_{J^\varepsilon,N}^2$}
\addplot[solid,thick,color=red,mark=none,mark size=1pt] table [x index=0, y index=11]{practical_ns-ou_dt1.0E-05_N1.0E+04_data.txt};
\addlegendentry{$\hat \sigma_{K^\varepsilon,N}^2$}
\end{axis}
\node at (5cm,2.8cm) {\textcolor{black}{(b)}};
\end{tikzpicture}
\begin{tikzpicture}[scale=0.56]
\begin{axis}[legend style={at={(1,1)},anchor=north east}, compat=1.3,
  xmin=0.01, xmax=1,ymin=0,ymax=0.009,
  xlabel= {$\varepsilon$},
  ylabel= {}]
\addplot[solid, thick, densely dotted,color=black,mark=none,mark size=1pt] table [x index=0, y index=3]{practical_ns-ou_dt1.0E-05_N1.0E+04_data.txt};
\addlegendentry{$\hat \sigma_{K^\varepsilon,N}^2/\varepsilon$}
\addplot[solid, thick, color=black,mark=none,mark size=1pt] table [x index=0, y index=4]{practical_ns-ou_dt1.0E-05_N1.0E+04_data.txt};
\addlegendentry{$\hat \sigma_{K^\varepsilon,N}^2/\varepsilon^2$}
\end{axis}
\node at (5cm,3.7cm) {\textcolor{black}{(c)}};
\end{tikzpicture}


\begin{tikzpicture}[scale=0.56]
\begin{axis}[legend style={at={(0,1)},anchor=north west}, compat=1.3,
  xmin=0.01, xmax=1,ymin=0.25,ymax=0.45,
  xlabel= {$\varepsilon$},
  ylabel= {}]
\addplot[densely dashed,color=black,mark=none,mark size=1pt] table [x index=0, y index=5]{proba-practical_ns-ou_dt1.0E-05_N1.0E+04_data.txt};
\addlegendentry{$\hat{I}_N^\varepsilon$}
\addplot[densely dotted,thick,color=blue,mark=none,mark size=1pt] table [x index=0, y index=6]{proba-practical_ns-ou_dt1.0E-05_N1.0E+04_data.txt};
\addlegendentry{$\hat{J}_N^\varepsilon$}
\addplot[solid,color=red,mark=none,mark size=1pt] table [x index=0, y index=7]{proba-practical_ns-ou_dt1.0E-05_N1.0E+04_data.txt};
\addlegendentry{$\hat{K}_N^\varepsilon$}
\addplot[thick,dotted,color=black,mark=none,mark size=4pt] table [x index=0, y index=8]{proba-practical_ns-ou_dt1.0E-05_N1.0E+04_data.txt};
\addlegendentry{$I^0$}
\addplot[only marks,
  black, mark options={black,scale=0.25}, 
  error bars/.cd, 
    y fixed,
    y dir=both, 
    y explicit
] table [x=x, y=y,y error=error, col sep=comma] {
    x,  y,       error
0.0875,  0.2696,  0.00869754  
0.4875,  0.3235,  0.00916911
0.8875,  0.4245,  0.00968763  
};
\addplot[only marks,
  blue, mark options={black,scale=0.25}, 
  error bars/.cd, 
    y fixed,
    y dir=both, 
    y explicit
] table [x=x, y=y,y error=error, col sep=comma] {
    x,  y,       error
0.1,  0.26327,  0.00474802 
0.5,  0.32617,  0.0104662
0.9,  0.41627,  0.0122829
};
\addplot[only marks,
  red, mark options={black,scale=0.25}, 
  error bars/.cd, 
    y fixed,
    y dir=both, 
    y explicit
] table [x=x, y=y,y error=error, col sep=comma] {
    x,  y,       error
0.1125,  0.264234,  0.00455764
0.5125,  0.32438,  0.00871465 
0.9125,  0.423444,  0.00962228
};
\end{axis}
\node at (5cm,2.5cm) {\textcolor{black}{(a)}};
\end{tikzpicture}
\begin{tikzpicture}[scale=0.56]
\begin{axis}[legend style={at={(1,0)},anchor=south east}, compat=1.3,
  xmin=0.01, xmax=1,ymin=0,ymax=0.4125,
  xlabel= {$\varepsilon$},
  ylabel= {}]
\addplot[densely  dashed, thick, color=black,mark=none,mark size=1pt] table [x index=0, y index=9]{proba-practical_ns-ou_dt1.0E-05_N1.0E+04_data.txt};
\addlegendentry{$\hat \sigma_{I^\varepsilon,N}^2$}
\addplot[densely dotted, thick,color=blue,mark=none,mark size=1pt] table [x index=0, y index=10]{proba-practical_ns-ou_dt1.0E-05_N1.0E+04_data.txt};
\addlegendentry{$\hat \sigma_{J^\varepsilon,N}^2$}
\addplot[solid,thick,color=red,mark=none,mark size=1pt] table [x index=0, y index=11]{proba-practical_ns-ou_dt1.0E-05_N1.0E+04_data.txt};
\addlegendentry{$\hat \sigma_{K^\varepsilon,N}^2$}
\end{axis}
\node at (5cm,3.9cm) {\textcolor{black}{(b)}};
\end{tikzpicture}
\begin{tikzpicture}[scale=0.56]
\begin{axis}[legend style={at={(1,1)},anchor=north east}, compat=1.3,
  xmin=0.01, xmax=1,ymin=0,ymax=2.0,
  xlabel= {$\varepsilon$},
  ylabel= {}]
\addplot[solid, thick, densely dotted,color=black,mark=none,mark size=1pt] table [x index=0, y index=3]{proba-practical_ns-ou_dt1.0E-05_N1.0E+04_data.txt};
\addlegendentry{$\hat \sigma_{K^\varepsilon,N}^2/\varepsilon$}
\addplot[solid, thick, color=black,mark=none,mark size=1pt] table [x index=0, y index=4]{proba-practical_ns-ou_dt1.0E-05_N1.0E+04_data.txt};
\addlegendentry{$\hat \sigma_{K^\varepsilon,N}^2/\varepsilon^2$}
\end{axis}
\node at (5cm,3.7cm) {\textcolor{black}{(c)}};
\end{tikzpicture}
\caption{Example \ref{example81c}  of a TDOF modeling a two-storey building with an hysteretic spring of elasto-plastic type, driven by an Ornstein-Uhlenbeck noise. The quantities presented here are similar to those presented in Figure \ref{fig:i2}. The details of the elements \textcircled{1} and \textcircled{3} with nonlinear and hysteretic behaviours  can be found in Example \ref{example81c}.
}
\label{fig:i3}
\end{figure}
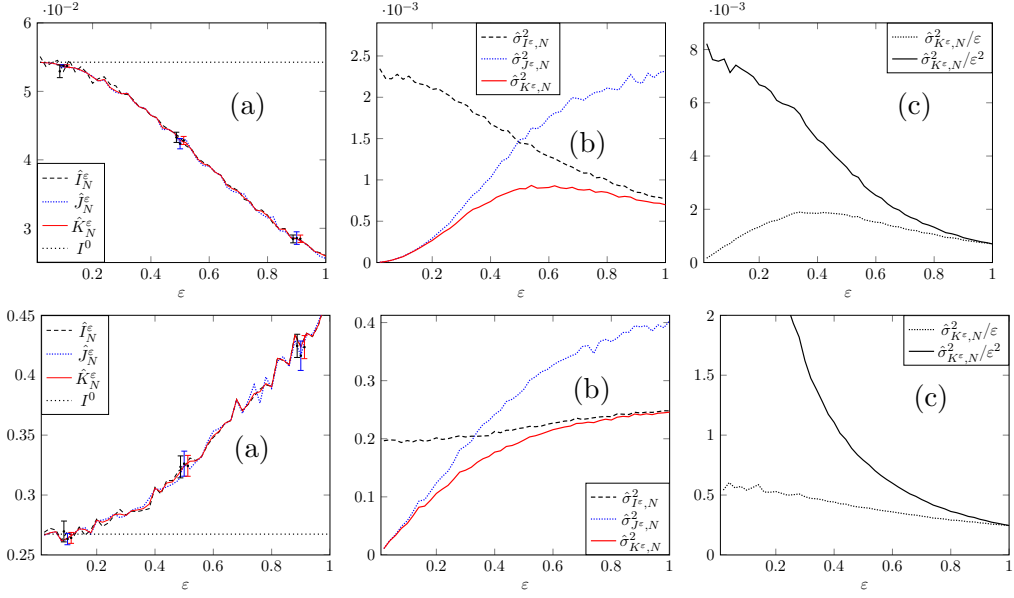
 
 \section{Diffusion approximation for a driving multivariate Ornstein-Uhlenbeck process}
\label{sec:adif1}
We consider the $\RR^n$-valued process $\bX^\eps$ solution of the ODE
(\ref{eq:sde1})
when $\boeta^\eps$ is the multivariate $d$-dimensional Ornstein-Uhlenbeck process (\ref{eq:sde2}).
We give several explicit examples.

\begin{example}
$\eta^\eps$ is a one-dimensional Ornstein-Uhlenbeck process, $d=d'=1$, $A,K>0$,
\begin{equation}
\label{eq:ou}
 {\rm d}\eta^\eps = -\frac{A}{\eps^2} \eta^\eps {\rm d}t + \frac{K}{\eps} {\rm d}W_t.
\end{equation}
\end{example}

\begin{example}
$\boeta^\eps$ is a  Langevin process
\begin{align}
\label{eq:lang1a}
{\rm d} \eta_1^\eps &=\frac{1}{\eps^2}  \eta_2^\eps  {\rm d}t ,\\
{\rm d} \eta_2^\eps &= - \frac{1}{\eps^2}
\big[ \mu \eta_1^\eps +\gamma \eta_2^\eps  \big] {\rm d}t +
\frac{K}{\eps} {\rm d}W_t,
\label{eq:lang1b}
\end{align}
which corresponds to $d=2$, $d'=1$, ${\bf A} = \begin{pmatrix} 0 & -1\\ \mu &\gamma \end{pmatrix}$, and ${\bf K} = \begin{pmatrix} 0  \\ K \end{pmatrix}$.
The process $\eta^\eps_{1}$ is a white-noise driven linear oscillator with stiffness $\mu>0$ and damping $\gamma>0$.
It can be encountered in earthquake engineering because it is considered to be a realistic type of random forcing to represent seismic
excitation (it is the so-called Kanai-Tajimi model  \cite{lin}).
\end{example}

\begin{example}
If $\tilde{\eta}^\eps$ is a real-valued zero-mean stationary Gaussian process with power spectral density 
${\rm PSD}^\eps (\omega) = \eps^2 {\rm PSD}(\eps^2 \omega)$, 
${\rm PSD}(\omega) =  \sum_{k=1}^q \frac{\sigma_k^2}{1+\omega^2/\Delta \Omega_k^2} $, then it has the same distribution as the process $\sum_{k=1}^q \sigma_k \eta_{k}^\eps$ where $\boeta^\eps$ is solution of (\ref{eq:sde2}) with  $d=d'=q$ and 
$$
{\bf A} = {\bf K} = {\rm diag} \big(   \Delta \Omega_k  , k=1,\ldots,q\big)  .
$$
This shows that any zero-mean stationary Gaussian process with power spectral density 
that can be decomposed as a sum of centered Lorentzians belongs to the model (\ref{eq:sde2}).
\end{example}

\begin{example}
If $\tilde{\eta}^\eps$ is a real-valued zero-mean stationary Gaussian process with power spectral density 
${\rm PSD}^\eps (\omega) = \eps^2 {\rm PSD}(\eps^2 \omega)$, 
${\rm PSD}(\omega) = \frac{1}{2} \sum_{k=1}^q \frac{\sigma_k^2}{1+(\omega-\omega_k)^2/\Delta \Omega_k^2} + \frac{\sigma_k^2}{1+(\omega+\omega_k)^2/\Delta \Omega_k^2} $, then it has the same distribution as the process $\sum_{k=1}^q \sigma_k \eta_{2k-1}^\eps$ where $\boeta^\eps$ is solution of (\ref{eq:sde2}) with  $d=d'=2q$ and 
$$
{\bf A} = \oplus_{k=1}^q \begin{pmatrix} \Delta \Omega_k & -\omega_k \\  \omega_k &\Delta \Omega_k  \end{pmatrix} ,
\quad 
{\bf K} =  \oplus_{k=1}^q \begin{pmatrix}  \Delta \Omega_k & 0 \\  0 &  \Delta \Omega_k  \end{pmatrix}  .
$$
This shows that any zero-mean stationary Gaussian process with power spectral density 
that can be decomposed as a sum of non-centered Lorentzian functions belongs to the model (\ref{eq:sde2}).
\end{example}

\noindent We also consider the limiting $\RR^n$-valued process $\bU$ solution of the SDE (\ref{eq:sde3}).
The continuous process $(\bX^\eps-\bU)$ converges in probability to zero
as $\eps\to 0$ as stated in the following proposition.

\begin{prop}
\label{prop:0}
If $\bX^\eps_0=\bU_0$, then 
the continuous process $(\bX^\eps-\bU)$ converges in probability to zero as $\eps \to 0$.
The convergence holds in the space of continuous functions equipped with the topology associated to the 
uniform norm over compact intervals.
\end{prop}

The proof of Proposition \ref{prop:0} is based on the perturbed test function method as described first in \cite[Chapter 7]{kushner90} or in \cite[Chapter 6]{book}. It is given in Appendix \ref{app:prooflemmaptf}.


\begin{example}
We consider the process $\bX^\eps$ solution of the ODE \eqref{eq:sde1}
where $\boeta^\eps$ is the rapidly varying mean-reverting process \eqref{eq:ou}.
We also consider the limiting process
$$
{\rm d}\bU = \bb(\bU) {\rm d}t + \frac{K}{A} \bsigma (\bU) {\rm d}W_t +\frac{K^2}{2 A^2}
 ( \bsigma (\bU) \cdot \nabla_{\bx^0})\bsigma(\bU) {\rm d}t,
$$
driven by the same Brownian motion. The continuous process $(\bX^\eps-\bU)$ converges in probability to zero
as $\eps\to 0$.
\end{example}

\begin{example}
We consider the process $\bX^\eps$ solution of the ODE \eqref{eq:sde1}
where $\boeta^\eps$ is the rapidly varying mean-reverting process (\ref{eq:lang1a}-\ref{eq:lang1b}).
We also consider the limiting process
$$
{\rm d}\bU = \bb(\bU) {\rm d}t + \frac{K}{\mu} \bsigma(\bU) {\rm d}W_t +\frac{K^2}{2 \mu^2}( \bsigma (\bU) \cdot \nabla_{\bx^0})\bsigma(\bU) {\rm d}t,
$$
driven by the same Brownian motion.  
The continuous process $(\bX^\eps-\bU)$ converges in probability to zero
as $\eps\to 0$.
\end{example}

\textcolor{black}{
 The proof that the optimal control variate estimator $\hat{K}^\eps_N$ and the control variate estimator $\hat{J}^\eps_N$ have asymptotic variances of the order of $\eps^2$ as stated in Proposition \ref{prop:21}
 follows from the following lemma.
\begin{lem}
\label{prop:3}
Let ${f},{g}$ be smooth functions from $\RR^n$ to $\RR$ with bounded derivatives.
Let $T>0$.
There exists $C>0$ such that, for any $t \in [\eps,T]$,
\begin{equation}
\label{eq:estimeffeps2}
\big| \EE \big[{g}(\bU_t) \big({f}(\bX^\eps_t)-{f}(\bU_t)\big) \big] \big| \leq C \eps^2 ,\quad \quad
\EE \big[ \big({f}(\bX^\eps_t)-{f}(\bU_t)\big)^2\big] \leq C \eps^2 .
\end{equation}
\end{lem}
The important hypothesis is that $f$ should be smooth. We could certainly relax the hypothesis 
on the bounded derivatives by using uniform estimates of high-order moments of the process $\bX^\eps$.
Lemma \ref{prop:3} is proved in Appendix \ref{app:proofprop3}.
}

\section{Numerical simulations}
\label{sec:num}
In this section, we illustrate our control variate method and report the numerical results on different types of dynamical systems driven by colored noises.
The two  examples are smooth oscillators that can be described by Equation \eqref{eq:sde1} (one being linear with time-dependent coefficients and the other being of Van der Pol type). 
Other examples with non-smooth dynamical systems will be addressed in Section \ref{sec:nummulti}.

We use the Euler-Maruyama approximation method to compute the approximate numerical solution of a SDE \cite{kloeden}.
In Subsection \ref{subsec:num1}, we recall the two types of colored noise that we consider and provide their time discretization.
Then, in Subsection \ref{subsec:num2}, some details and discretization of the dynamical systems under consideration are given.
Finally, in Subsection \ref{subsec:num3}, numerical experiments on the control variate estimators are provided and discussed in each case.

\subsection{Colored noise models and their discretization}
\label{subsec:num1}
\noindent The two models of noise are shown in Eq.~\eqref{eq:ou} (OU) and in the system of equations (\ref{eq:lang1a}-\ref{eq:lang1b}) (Langevin).
The OU noise has two parameters $A,K_{\textup{ou}} > 0$ whereas the Langevin has three parameters $\mu,\gamma,K_{\textup{lan}} >0$.
Their discretization works as follows. 
Let $T>0$ and $N_T \in \mathbb{N}$ be the number of time steps such that $T=N_T \delta t$. 
Let $N$ be the number of Monte Carlo samples. 
Consider a sequence of independent and identically distributed  standard Gaussian variables 
$$
\{ \Delta W_n^{k} \sim  \mathcal{N} (0,1) , \: 0 \leq n \leq N_T-1, \: 1 \leq {k} \leq N \}.
$$
Let $\eps>0$. 
For each $1 \leq {k} \leq  N$, 
we overload the notation by denoting the discretized noise in both cases by $\{ \hat \eta_n^{\eps,{k}}, 0 \leq n \leq N_T \}$.
\begin{itemize}
\item
\textbf{Ornstein-Uhlenbeck noise:}
$\hat \eta_0^{\eps,{k}} \sim \mathcal{N}\left (0, \dfrac{K_{\textup{ou}}^2}{2 A} \right )$
and for $0 \leq n \leq N_T-1$,
$$
\hat \eta_{n+1}^{\eps, {k}} = \hat \eta_n^{\eps, {k}} \left (1 - \delta t \dfrac{A}{\eps^2} \right ) + \sqrt{\delta t} \dfrac{K_{\textup{ou}}}{\eps} \Delta W_n^{k}.
$$ 
\item
\textbf{Langevin noise:} 
$\hat \eta_0^{\eps,{k}}$ and $\hat \eta_{2,0}^{\eps,{k}}$ are independent variables with
$$
\hat \eta_0^{\eps,{k}} \sim \mathcal{N}\left (0, \dfrac{K_{\textup{lan}}^2}{2 \gamma} \right ), 
\quad \hat \eta_{2,0}^{\eps,{k}} \sim \mathcal{N}\left (0, \dfrac{K_{\textup{lan}}^2}{2 \gamma \mu} \right )
$$
and for $0 \leq n \leq N_T-1$,
$$
\hat \eta_{n+1}^{\eps,{k}} = \hat \eta_{n}^{\eps,{k}} + \dfrac{\delta t}{\eps^2} \hat \eta_{2,n}^{\eps,{k}},
\quad \hat \eta_{2,n+1}^{\eps,{k}} = \hat \eta_{2,n}^{\eps,{k}} 
- \dfrac{\delta t}{\eps^2} \big[ \mu \hat \eta_n^{\eps,{k}} +\gamma \hat \eta_{2,n}^{\eps,{k}}  \big] 
+ \sqrt{\delta t} \dfrac{K_{\textup{lan}}}{\eps} \Delta W_n^{k}.
$$
\end{itemize}

\subsection{Details and discretization of the illustrative dynamical systems}
\label{subsec:num2}
We consider systems of the form of \eqref{eq:sde1}.
We first consider the case of smooth systems that can have time-dependent coefficients, 
\begin{equation}
\label{eq:smoothX}
\dfrac{{\rm d} X_1^\eps}{{\rm d}t} = X_2^\eps,\\
\quad \dfrac{{\rm d} X_2^\eps}{ {\rm d}t} = - h (X_1^\eps ,X_2^\eps,t) + \frac{1}{\eps} \eta^\eps  .
\end{equation}
Here we are interested in $\mathbb{E} [ \|\bX^\eps_T\|^2 ]$ and in $\PP( | X^\eps_{1,T} | \leq 1)$ for $T=1$. 
Note that the second case corresponds to an expectation $\PP ( | X^\eps_{1,T} | \leq 1   ) = \EE[ f(\bX^\eps_T)]$ 
with a non-smooth function $f(\bx) = {\bf 1}_{ |x_1| \leq 1 }$.
As $\eps \to 0$, $\bX^\eps = (X_1^\eps ,X_2^\eps) $ converges to $ \bU = (X^0_1,X^0_2)$ where 
\begin{equation}
\label{eq:smoothU}
{\rm d} X^0_1 = X^0_2 {\rm d} t ,
\quad {\rm d} X^0_2 = - h (X^0_1,X^0_2,t) {\rm d} t  + C {\rm d} W,
\end{equation}
$C = K_{\textup{ou}} A^{-1}$ for an OU noise 
and $C = K_{\textup{lan}} \mu^{-1}$ for a Langevin noise.
For the stochastic simulation of \eqref{eq:smoothX} and \eqref{eq:smoothU}, we proceed as follows:
\begin{itemize}
\item
$\hat X_{1,0}^{\eps,{k}} = x_{1,0}, \: \hat X_{2,0}^{\eps,{k}} = x_{2,0}$
and for $0 \leq n \leq N_T-1$,
$$
\begin{cases}
\hat X_{1,n+1}^{\eps, {k}} = \hat X_{1,n}^{\eps, {k}} + \delta t  \hat X_{2,n}^{\eps, {k}} ,\\
\hat X_{2,n+1}^{\eps, {k}} = 
\hat X_{2,n}^{\eps, {k}} 
- \delta t  h ( \hat X_{1,n}^{\eps, {k}}, \hat X_{2,n}^{\eps, {k}}, n \delta t )
+ \frac{\delta t}{\eps} \hat \eta_n^{\eps, {k}}.
\end{cases}
$$
\item
$\hat X_{1,0}^{0, {k}} = x_{1,0}, \: \hat X_{2,0}^{0, {k}} = x_{2,0}$ and for $0 \leq n \leq N_T-1$,
$$
\begin{cases}
\hat X_{1,n+1}^{0, {k}} = \hat X_{1,n}^{0, {k}} + \delta t  \hat X_{2,n}^{0,  {k}} ,\\
\hat X_{2,n+1}^{0, {k}} = 
\hat X_{2,n}^{0, {k}} 
- \delta t  h ( \hat X_{1,n}^{0,{k}}, \hat X_{2,n}^{0, {k}}, n \delta t )
+ C \sqrt{\delta t} \Delta W_n^{k}.
\end{cases}
$$ 
$\hat \bX_n^{\eps,{k}}$ and $\hat \bX_n^{0,{k}}$ are independent (in ${k}$) copies that are meant to approximate $\bX^\eps_{n \delta t}$ and $\bU_{n \delta t}$.
\end{itemize}

\begin{example}[\textbf{linear oscillator with time-dependent coefficients}]
\label{example3a}
We take $h(x_1,x_2,t) \triangleq p(t) x_1 + q(t) x_2$ where 
$p(t) \triangleq 1+ \cos(t)$ and $q(t) \triangleq 1+ \sin(t)$ (the choice is purely arbitrary).
Here, in both OU and Langevin cases, the limiting process $\bU=(X^0_1,X^0_2)$ is a Gaussian process provided that the initial condition is deterministic or Gaussian. 
This is useful to derive the expectation of the control variate.
The distribution of $\bU_t=(X^0_{1,t},X^0_{2,t})$ is characterized by its first-order moment ${\itbf m}(t) \triangleq \mathbb{E} [\bU_{t}] \in \RR^2$ and second-order moment ${\bf M}(t) \triangleq (\EE  [ X^0_{i,t} X^0_{j,t}])_{i,j=1}^2 \in {\cal M}_{2,2}( \RR)$ which satisfy the following systems of differential equations:
\begin{itemize}
\item
\textbf{first-order moment}
\begin{itemize}
\item
$(m_1(0),m_2(0)) = (x_0,\dot x_0)$,
\item
$\dot m_1(t) = m_2(t)$, 
\item
$\dot m_2(t) = -p(t)m_1(t) - q(t)m_2(t)$.
\end{itemize}
\item
\textbf{second-order moment}
\begin{equation}
\label{eq:M}
\begin{cases}
& (M_{11}(0),M_{22}(0),M_{12}(0)) = (x_0^2,\dot x_0^2,x_0 \dot x_0),\\
& \dot M_{11}(t) = 2 M_{12}(t),\\
& \dot M_{22}(t) = -2 p(t) M_{12}(t) - 2 q(t) M_{22}(t) + C^2,\\
& \dot M_{12}(t) =  M_{22}(t) - p(t) M_{11}(t) - q(t) M_{12}(t).
\end{cases}
\end{equation}
\end{itemize}
The expectation of the control variate $\EE [ \|\bU_T\|^2 ]$ with $T=1$ is estimated by solving numerically, with an Euler method, the differential equations for the first-  and second-order moments. 
\end{example}

\begin{example}[\textbf{Van der Pol oscillator}]
\label{example3b}
We take $h(x_1,x_2) =  x_1 - \nu (1-x_1^2) x_2$ where $\nu > 0$.
The expectation of the control variate can be represented by $\EE[\| \bU_T\|^2 ] = \mathfrak{c}(\bx_{0},0)$ with $T=1$, 
where $\mathfrak{c}$ satisfies the following backward in time PDE 
\begin{equation}
\label{eq:pde_vdp}
\begin{cases}
& \partial_t \mathfrak{c} +\frac{C^2}{2} \partial_{x_2}^2 \mathfrak{c}  - h (x_1,x_2) \partial_{x_2} \mathfrak{c}  + x_2 \partial_{x_1} \mathfrak{c} =0 ,  \: \textup{in} \: \mathbb{R}^2 \times [0,1)\\
& \mathfrak{c}(\bx,1) =  \| \bx\|^2\: \textup{in} \: \mathbb{R}^2  .
\end{cases}
\end{equation}
The expectation $\EE[\| \bU_T\|^2 ] $ is estimated by solving the PDE (\ref{eq:pde_vdp}) with a finite difference method.
\end{example}

\begin{example}[\textbf{nonlinear TDOF}]
\label{example3c}
A TDOF modeling a nonlinear spring of the linear-plus-quadratic cubic type and a nonlinear damper of the linear-plus-quadratic type, can be seen as a coupling between two systems of the form \eqref{eq:smoothX}
\begin{equation}
\label{eq:smoothXcouple}
\begin{cases}
\dfrac{{\rm d} X_1^\eps}{{\rm d}t} = X_2^\eps, 
\quad \dfrac{{\rm d} X_2^\eps}{ {\rm d}t} = - g_2 (X_1^\eps ,X_2^\eps,X_3^\eps ,X_4^\eps) + \dfrac{1}{\eps} \eta^\eps , \\[2mm]
\dfrac{{\rm d} X_3^\eps}{{\rm d}t} = X_4^\eps,
\quad \dfrac{{\rm d} X_4^\eps}{ {\rm d}t} = - g_4(X_1^\eps ,X_2^\eps,X_3^\eps ,X_4^\eps) .
\end{cases}
\end{equation}
In addition to $\EE \big[ (X_{1,T}^\eps )^2+ (X_{3,T}^\eps )^2 \big]$ for $T=1$, we are interested in 
$\PP( | X^\eps_{1,T} | \leq a, | X^\eps_{3,T} | \leq b)$.
As $\eps \to 0$, $\bX^\eps \to \bU$ where
\begin{equation}
\label{eq:smoothUcouple}
\begin{cases}
{\rm d} X^0_1 = X^0_2  {\rm d} t,
\quad {\rm d} X^0_2 = - g_2 (X^0_1,X^0_2,X^0_3,X^0_4) {\rm d} t + C {\rm d}W, \\[2mm]
{\rm d} X^0_3 = X^0_4  {\rm d} t,
\quad {\rm d} X^0_4 = - g_4 (X^0_1,X^0_2,X^0_3,X^0_4)  {\rm d} t.
\end{cases}
\end{equation}
Here 
$$
g_2(x_1,x_2,x_3,x_4) \triangleq 
k_1 x_1 (1+\min(x_1^2,\textcolor{black}{\tilde L}) )+ c_1 x_2
- k_3 (x_3-x_1)
- c_3(x_4-x_2) (1 + \min(|x_4-x_2|,\textcolor{black}{\tilde L}))
$$
and 
$$
g_4(x_1,x_2,x_3,x_4) \triangleq 
k_3 (x_3-x_1) + c_3(x_4-x_2) (1 + \min(|x_4-x_2|, \textcolor{black}{\tilde L})).
$$
In the original model of Spanos \textcolor{black}{$\tilde L$} $= \infty$, \textcolor{black}{see pages 189-190 in \cite{spanosbook}}. For any positive finite value of $\textcolor{black}{\tilde L}$, the system above enters into the scope of our results. The simulation of \eqref{eq:smoothXcouple} and \eqref{eq:smoothUcouple} is similar to what is done for \eqref{eq:smoothX} and \eqref{eq:smoothU}.
We take $\textcolor{black}{\tilde L}=1000$, $c_1=c_3=k_1=k_3=1$. 
\end{example}

\subsection{Numerical experiments}
\label{subsec:num3}%
We report our numerical results for the two systems mentioned above.
In each of the two figures below, there are four subfigures (a)-(b)-(c)-(d). For subfigures  (a) and (b), the driving force is an Ornstein-Uhlenbeck noise \eqref{eq:ou} with $A=K=1$. 
In subfigure (a), 
\textcolor{black}{the dashed black, dotted blue, and solid red lines represent the  standard MC estimator $\hat{I}_N^\eps $ and the control variate estimators $\hat{J}_N^\eps$ and $\hat{K}_N^\eps$, respectively.
For each value of $\varepsilon \in \{0.1,0.5,0.9\}$, error bars (95\% confidence interval) are shown for each of the estimators (in $\hat{I}_N^\eps$, $\hat{J}_N^\eps$, $\hat{K}^\eps_N$ order from the left to the right).}
The black dotted line represents the expectation of the control variate. The objective of the subfigure (b) is to illustrate the bound \eqref{eq:boundepsilon2}
and to show that the $\eps^2$-behavior is actually sharp.
The same description applies to (c) and (d), except they correspond to the case of a Langevin noise (\ref{eq:lang1a}-\ref{eq:lang1b})
with $\mu=\gamma=K=1$.
\textcolor{black}{
In the figures 
 the asymptotic variance of the standard MC estimator $\hat{I}^\eps_N$ is estimated by
\begin{equation}
\label{eq:estimvarJneps}
\widehat{\sigma}_{I^\eps,N}^2  =  \frac{1}{N} 
\sum_{k=1}^N \big(F (\bX^\eps(\bW^k) ) \big)^2  - ({\hat{I}^\eps_N})^2 ,
\end{equation}
the asymptotic variance of the control variate estimator $\hat{J}^\eps_N$ is estimated by
\begin{equation}
\label{eq:estimvarIneps}
\widehat{\sigma}_{J^\eps,N}^2  =  \frac{1}{N} 
\sum_{k=1}^N \big( F (\bX^\eps(\bW^k) ) - F ( \bU  (\bW^k) )  +I^0 \big)^2 - (\hat{J}^\eps_N)^2 ,
\end{equation}
and
the asymptotic variance of the optimal control variate estimator $\hat{K}^\eps_N$ is estimated by
\begin{equation}
\label{eq:estimvarKneps}
\widehat{\sigma}_{K^\eps,N}^2  =  \frac{1}{N} 
\sum_{k=1}^N \big( F(\bX^\eps(\bW^k) ) -  \hat{\rho}_N^\eps F ( \bU(\bW^k) )  +\hat{\rho}_N^\eps I^0\big)^2 - (\hat{K}^\eps_N)^2 ,
\end{equation}
with $\hat{\rho}_N^\eps$ defined by (\ref{def:hatrhoNeps}).
$\widehat{\sigma}_{I^\eps,N}^2$, $\widehat{\sigma}_{J^\eps,N}^2 $, and $\widehat{\sigma}_{K^\eps,N}^2 $ are consistent estimators of 
${\sigma}_{I^\eps}^2$, ${\sigma}_{J^\eps}^2$, and ${\sigma}_{K^\eps}^2$, 
respectively.
}

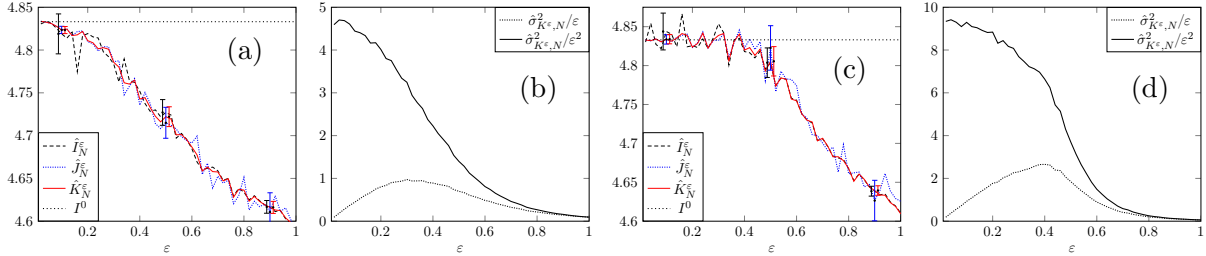
\begin{figure}[h!]
\centering
\begin{tikzpicture}[scale=0.5]
\begin{axis}[legend style={at={(0,0)},anchor=south west}, compat=1.3,
  xmin=0.01, xmax=1,ymin=4.6,ymax=4.85,
  xlabel= {$\varepsilon$},
  ylabel= {}]
\addplot[densely dashed,color=black,mark=none,mark size=1pt] table [x index=0, y index=5]{time_dep_oscillator_noise-ou_dt1.0E-05_N1.0E+04_data.txt};
\addlegendentry{$\hat{I}_N^\varepsilon$}
\addplot[densely dotted,thick,color=blue,mark=none,mark size=1pt] table [x index=0, y index=6]{time_dep_oscillator_noise-ou_dt1.0E-05_N1.0E+04_data.txt};
\addlegendentry{$\hat{J}_N^\varepsilon$}
\addplot[solid,color=red,mark=none,mark size=1pt] table [x index=0, y index=7]{time_dep_oscillator_noise-ou_dt1.0E-05_N1.0E+04_data.txt};
\addlegendentry{$\hat{K}_N^\varepsilon$}
\addplot[thick,dotted,color=black,mark=none,mark size=4pt] table [x index=0, y index=8]{time_dep_oscillator_noise-ou_dt1.0E-05_N1.0E+04_data.txt};
\addlegendentry{$I^0$}
\addplot[only marks,
  black, mark options={black,scale=0.25}, 
  error bars/.cd, 
    y fixed,
    y dir=both, 
    y explicit
] table [x=x, y=y,y error=error, col sep=comma] {
    x,  y,       error
0.0875, 4.81887,  0.0233356    
0.4875,  4.72708,  0.0150348 
0.8875,  4.61676,  0.00743379 
};
\addplot[only marks,
  blue, mark options={black,scale=0.25}, 
  error bars/.cd, 
    y fixed,
    y dir=both, 
    y explicit
] table [x=x, y=y,y error=error, col sep=comma] {
    x,  y,       error
0.1,  4.82352,  0.00421625  
0.5,  4.71498,  0.0180721 
0.9,  4.61068,  0.0225433  
};
\addplot[only marks,
  red, mark options={black,scale=0.25}, 
  error bars/.cd, 
    y fixed,
    y dir=both, 
    y explicit
] table [x=x, y=y,y error=error, col sep=comma] {
    x,  y,       error
0.1125, 4.82335,  0.00413095  
0.5125, 4.72213,  0.0115796  
0.9125,  4.61615,  0.00704399 
};
\end{axis}
\node at (5.5cm,4.5cm) {\textcolor{black}{(a)}};
\end{tikzpicture}
\begin{tikzpicture}[scale=0.5]
\begin{axis}[legend style={at={(1,1)},anchor=north east}, compat=1.3,
  xmin=0.01, xmax=1,ymin=0,ymax=5,
  xlabel= {$\varepsilon$},
  ylabel= {}]
\addplot[solid, thick, densely dotted,color=black,mark=none,mark size=1pt] table [x index=0, y index=3]{time_dep_oscillator_noise-ou_dt1.0E-05_N1.0E+04_data.txt};
\addlegendentry{$\hat \sigma_{K^\varepsilon,N}^2/\varepsilon$}
\addplot[solid, thick, color=black,mark=none,mark size=1pt] table [x index=0, y index=4]{time_dep_oscillator_noise-ou_dt1.0E-05_N1.0E+04_data.txt};
\addlegendentry{$\hat \sigma_{K^\varepsilon,N}^2/\varepsilon^2$}
\end{axis}
\node at (5.5cm,3.5cm) {\textcolor{black}{(b)}};
\end{tikzpicture}
\centering
\begin{tikzpicture}[scale=0.5]
\begin{axis}[legend style={at={(0,0)},anchor=south west}, compat=1.3,
  xmin=0.01, xmax=1,ymin=4.6,ymax=4.875,
  xlabel= {$\varepsilon$},
  ylabel= {}]
\addplot[densely dashed,color=black,mark=none,mark size=1pt] table [x index=0, y index=5]{time_dep_oscillator_noise-lan_dt1.0E-05_N1.0E+04_data.txt};
\addlegendentry{$\hat{I}_N^\varepsilon$}
\addplot[densely dotted,thick,color=blue,mark=none,mark size=1pt] table [x index=0, y index=6]{time_dep_oscillator_noise-lan_dt1.0E-05_N1.0E+04_data.txt};
\addlegendentry{$\hat{J}_N^\varepsilon$}
\addplot[solid,color=red,mark=none,mark size=1pt] table [x index=0, y index=7]{time_dep_oscillator_noise-lan_dt1.0E-05_N1.0E+04_data.txt};
\addlegendentry{$\hat{K}_N^\varepsilon$}
\addplot[thick,dotted,color=black,mark=none,mark size=4pt] table [x index=0, y index=8]{time_dep_oscillator_noise-lan_dt1.0E-05_N1.0E+04_data.txt};
\addlegendentry{$I^0$}
\addplot[only marks,
  black, mark options={black,scale=0.25}, 
  error bars/.cd, 
    y fixed,
    y dir=both, 
    y explicit
] table [x=x, y=y,y error=error, col sep=comma] {
    x,  y,       error
0.0875, 4.84379,  0.0236637   
0.4875,  4.80372,  0.0189879 
0.8875,  4.63857,  0.00635131
};
\addplot[only marks,
  blue, mark options={black,scale=0.25}, 
  error bars/.cd, 
    y fixed,
    y dir=both, 
    y explicit
] table [x=x, y=y,y error=error, col sep=comma] {
    x,  y,       error
0.1,  4.83363,  0.00602099   
0.5,  4.8227,  0.0285681
0.9,  4.62651,  0.026038  
};
\addplot[only marks,
  red, mark options={black,scale=0.25}, 
  error bars/.cd, 
    y fixed,
    y dir=both, 
    y explicit
] table [x=x, y=y,y error=error, col sep=comma] {
    x,  y,       error
0.1125,  4.83395,  0.0059747
0.5125,  4.80552,  0.0188529 
0.9125,  4.6394,  0.00613961 
};
\end{axis}
\node at (5.5cm,4.0cm) {\textcolor{black}{(c)}};
\end{tikzpicture}
\begin{tikzpicture}[scale=0.5]
\begin{axis}[legend style={at={(1,1)},anchor=north east}, compat=1.3,
  xmin=0.01, xmax=1,ymin=0,ymax=10,
  xlabel= {$\varepsilon$},
  ylabel= {}]
\addplot[solid, thick, densely dotted,color=black,mark=none,mark size=1pt] table [x index=0, y index=3]{time_dep_oscillator_noise-lan_dt1.0E-05_N1.0E+04_data.txt};
\addlegendentry{$\hat \sigma_{K^\varepsilon,N}^2/\varepsilon$}
\addplot[solid, thick, color=black,mark=none,mark size=1pt] table [x index=0, y index=4]{time_dep_oscillator_noise-lan_dt1.0E-05_N1.0E+04_data.txt};
\addlegendentry{$\hat \sigma_{K^\varepsilon,N}^2/\varepsilon^2$}
\end{axis}
\node at (5.5cm,3.5cm) {\textcolor{black}{(d)}};
\end{tikzpicture}
\caption{Example \ref{example3a} (linear oscillator with time-dependent coefficients) with 
$h(x_1,x_2,t) \triangleq p(t) x_1 + q(t) x_2$, $p(t) \triangleq 1+ \cos(t)$ and $q(t) \triangleq 1+ \sin(t)$.
The target is to estimate $I^\eps=\EE[\|\bX^\eps_T\|^2]$ for $T=1$ and 
the expectation of the control variate ${I^0=}\EE[\| \bU_T\|^2 ]$ is obtained by solving the set of differential equations \eqref{eq:M}.}
\label{fig:1}
\end{figure}

\begin{figure}[h!]
\centering
\begin{tikzpicture}[scale=0.5]
\begin{axis}[legend style={at={(0,0)},anchor=south west}, compat=1.3,
  xmin=0.01, xmax=1,ymin=0.8,ymax=1.95,
  xlabel= {$\varepsilon$},
  ylabel= {}]
\addplot[densely dashed,color=black,mark=none,mark size=1pt] table [x index=0, y index=5]{vdp_noise-ou_dt1.0E-05_N1.0E+04_data.txt};
\addlegendentry{$\hat{I}_N^\varepsilon$}
\addplot[densely dotted,thick,color=blue,mark=none,mark size=1pt] table [x index=0, y index=6]{vdp_noise-ou_dt1.0E-05_N1.0E+04_data.txt};
\addlegendentry{$\hat{J}_N^\varepsilon$}
\addplot[solid,color=red,mark=none,mark size=1pt] table [x index=0, y index=7]{vdp_noise-ou_dt1.0E-05_N1.0E+04_data.txt};
\addlegendentry{$\hat{K}_N^\varepsilon$}
\addplot[thick,dotted,color=black,mark=none,mark size=4pt] table [x index=0, y index=8]{vdp_noise-ou_dt1.0E-05_N1.0E+04_data.txt};
\addlegendentry{$I^0$}
\addplot[only marks,
  black, mark options={black,scale=0.25}, 
  error bars/.cd, 
    y fixed,
    y dir=both, 
    y explicit
] table [x=x, y=y,y error=error, col sep=comma] {
    x,  y,       error
0.0875, 1.81387,  0.0374355
0.4875,  1.53632,  0.0337291
0.8875,  0.939618,  0.0227258
};
\addplot[only marks,
  blue, mark options={black,scale=0.25}, 
  error bars/.cd, 
    y fixed,
    y dir=both, 
    y explicit
] table [x=x, y=y,y error=error, col sep=comma] {
    x,  y,       error
0.1,  1.83173,  0.00561489 
0.5,  1.50726,  0.0277936  
0.9,  0.899007,  0.0371341  
};
\addplot[only marks,
  red, mark options={black,scale=0.25}, 
  error bars/.cd, 
    y fixed,
    y dir=both, 
    y explicit
] table [x=x, y=y,y error=error, col sep=comma] {
    x,  y,       error
0.1125,  1.83145,  0.0055845
0.5125,  1.51817,  0.0237851 
0.9125,  0.930788,  0.0210787 
};
\end{axis}
\node at (5.5cm,4.5cm) {\textcolor{black}{(a)}};
\end{tikzpicture}
\begin{tikzpicture}[scale=0.5]
\begin{axis}[legend style={at={(1,1)},anchor=north east}, compat=1.3,
  xmin=0.01, xmax=1,ymin=0,ymax=9,
  xlabel= {$\varepsilon$},
  ylabel= {}]
\addplot[solid, thick, densely dotted,color=black,mark=none,mark size=1pt] table [x index=0, y index=3]{vdp_noise-ou_dt1.0E-05_N1.0E+04_data.txt};
\addlegendentry{$\hat \sigma_{K^\varepsilon,N}^2/\varepsilon$}
\addplot[solid, thick, color=black,mark=none,mark size=1pt] table [x index=0, y index=4]{vdp_noise-ou_dt1.0E-05_N1.0E+04_data.txt};
\addlegendentry{$\hat \sigma_{K^\varepsilon,N}^2/\varepsilon^2$}
\end{axis}
\node at (5.5cm,3.75cm) {\textcolor{black}{(b)}};
\end{tikzpicture}
\centering
\centering
%
%
\begin{tikzpicture}[scale=0.5]
\begin{axis}[legend style={at={(0,0)},anchor=south west}, compat=1.3,
  xmin=0.01, xmax=1,ymin=0.8,ymax=1.95,
  xlabel= {$\varepsilon$},
  ylabel= {}]
\addplot[densely dashed,color=black,mark=none,mark size=1pt] table [x index=0, y index=5]{vdp_noise-lan_dt1.0E-05_N1.0E+04_data.txt};
\addlegendentry{$\hat{I}_N^\varepsilon$}
\addplot[densely dotted,thick,color=blue,mark=none,mark size=1pt] table [x index=0, y index=6]{vdp_noise-lan_dt1.0E-05_N1.0E+04_data.txt};
\addlegendentry{$\hat{J}_N^\varepsilon$}
\addplot[solid,color=red,mark=none,mark size=1pt] table [x index=0, y index=7]{vdp_noise-lan_dt1.0E-05_N1.0E+04_data.txt};
\addlegendentry{$\hat{K}_N^\varepsilon$}
\addplot[thick,dotted,color=black,mark=none,mark size=4pt] table [x index=0, y index=8]{vdp_noise-lan_dt1.0E-05_N1.0E+04_data.txt};
\addlegendentry{$I^0$}
\addplot[only marks,
  black, mark options={black,scale=0.25}, 
  error bars/.cd, 
    y fixed,
    y dir=both, 
    y explicit
] table [x=x, y=y,y error=error, col sep=comma] {
    x,  y,       error
0.0875, 1.83581,  0.0373308
0.4875,  1.88702,  0.0391248
0.8875,  1.17196,  0.0276693
};
\addplot[only marks,
  blue, mark options={black,scale=0.25}, 
  error bars/.cd, 
    y fixed,
    y dir=both, 
    y explicit
] table [x=x, y=y,y error=error, col sep=comma] {
    x,  y,       error
0.1,  1.85017,  0.00799738 
0.5,  1.89044,  0.0390406 
0.9,  1.16736,  0.0464236  
};
\addplot[only marks,
  red, mark options={black,scale=0.25}, 
  error bars/.cd, 
    y fixed,
    y dir=both, 
    y explicit
] table [x=x, y=y,y error=error, col sep=comma] {
    x,  y,       error
0.1125,  1.84985,  0.00795346
0.5125,  1.88874,  0.0342966
0.9125,  1.1718,  0.027636 
};
\end{axis}
\node at (5.5cm,4.5cm) {\textcolor{black}{(c)}};
\end{tikzpicture}
%
%
\begin{tikzpicture}[scale=0.5]
\begin{axis}[legend style={at={(1,1)},anchor=north east}, compat=1.3,
  xmin=0.01, xmax=1,ymin=0,ymax=19,
  xlabel= {$\varepsilon$},
  ylabel= {}]
\addplot[solid, thick, densely dotted,color=black,mark=none,mark size=1pt] table [x index=0, y index=3]{vdp_noise-lan_dt1.0E-05_N1.0E+04_data.txt};
\addlegendentry{$\hat \sigma_{K^\varepsilon,N}^2/\varepsilon$}
\addplot[solid, thick, color=black,mark=none,mark size=1pt] table [x index=0, y index=4]{vdp_noise-lan_dt1.0E-05_N1.0E+04_data.txt};
\addlegendentry{$\hat \sigma_{K^\varepsilon,N}^2/\varepsilon^2$}
\end{axis}
\node at (5.5cm,3.75cm) {\textcolor{black}{(d)}};
\end{tikzpicture}
\begin{tikzpicture}[scale=0.5]
\begin{axis}[legend style={at={(0,0)},anchor=south west}, compat=1.3,
  xmin=0.01, xmax=1,ymin=0.65,ymax=1.,
  xlabel= {$\varepsilon$},
  ylabel= {}]
\addplot[densely dashed,color=black,mark=none,mark size=1pt] table [x index=0, y index=5]{proba-vdp_noise-ou_dt1.0E-05_N1.0E+04_data.txt};
\addlegendentry{$\hat{I}_N^\varepsilon$}
\addplot[densely dotted,thick,color=blue,mark=none,mark size=1pt] table [x index=0, y index=6]{proba-vdp_noise-ou_dt1.0E-05_N1.0E+04_data.txt};
\addlegendentry{$\hat{J}_N^\varepsilon$}
\addplot[solid,color=red,mark=none,mark size=1pt] table [x index=0, y index=7]{proba-vdp_noise-ou_dt1.0E-05_N1.0E+04_data.txt};
\addlegendentry{$\hat{K}_N^\varepsilon$}
\addplot[thick,dotted,color=black,mark=none,mark size=4pt] table [x index=0, y index=8]{proba-vdp_noise-ou_dt1.0E-05_N1.0E+04_data.txt};
\addlegendentry{$I^0$}
\addplot[only marks,
  black, mark options={black,scale=0.25}, 
  error bars/.cd, 
    y fixed,
    y dir=both, 
    y explicit
] table [x=x, y=y,y error=error, col sep=comma] {
    x,  y,       error
0.0875, 0.827,  0.00741365
0.4875,  0.9052,  0.00574159
0.8875,  0.9828,  0.00254831  
};
\addplot[only marks,
  blue, mark options={black,scale=0.25}, 
  error bars/.cd, 
    y fixed,
    y dir=both, 
    y explicit
] table [x=x, y=y,y error=error, col sep=comma] {
    x,  y,       error
0.1,  0.829356,  0.00353298  
0.5,  0.904356,  0.00711109  
0.9,  0.974756,  0.00739902  
};
\addplot[only marks,
  red, mark options={black,scale=0.25}, 
  error bars/.cd, 
    y fixed,
    y dir=both, 
    y explicit
] table [x=x, y=y,y error=error, col sep=comma] {
    x,  y,       error
0.1125,  0.829076,  0.00342036 
0.5125,  0.904913,  0.0051608
0.9125,  0.982436,  0.00252688 
};
\end{axis}
\node at (5.5cm,4.5cm) {\textcolor{black}{(a)}};
\end{tikzpicture}
\begin{tikzpicture}[scale=0.5]
\begin{axis}[legend style={at={(1,1)},anchor=north east}, compat=1.3,
  xmin=0.01, xmax=1,ymin=0,ymax=3,
  xlabel= {$\varepsilon$},
  ylabel= {}]
\addplot[solid, thick, densely dotted,color=black,mark=none,mark size=1pt] table [x index=0, y index=3]{proba-vdp_noise-ou_dt1.0E-05_N1.0E+04_data.txt};
\addlegendentry{$\hat \sigma_{K^\varepsilon,N}^2/\varepsilon$}
\addplot[solid, thick, color=black,mark=none,mark size=1pt] table [x index=0, y index=4]{proba-vdp_noise-ou_dt1.0E-05_N1.0E+04_data.txt};
\addlegendentry{$\hat \sigma_{K^\varepsilon,N}^2/\varepsilon^2$}
\end{axis}
\node at (5.5cm,3.75cm) {\textcolor{black}{(b)}};
\end{tikzpicture}
\centering
\centering
%
%
\begin{tikzpicture}[scale=0.5]
\begin{axis}[legend style={at={(0,0)},anchor=south west}, compat=1.3,
  xmin=0.01, xmax=1,ymin=0.65,ymax=1,
  xlabel= {$\varepsilon$},
  ylabel= {}]
\addplot[densely dashed,color=black,mark=none,mark size=1pt] table [x index=0, y index=5]{proba-vdp_noise-lan_dt1.0E-05_N1.0E+04_data.txt};
\addlegendentry{$\hat{I}_N^\varepsilon$}
\addplot[densely dotted,thick,color=blue,mark=none,mark size=1pt] table [x index=0, y index=6]{proba-vdp_noise-lan_dt1.0E-05_N1.0E+04_data.txt};
\addlegendentry{$\hat{J}_N^\varepsilon$}
\addplot[solid,color=red,mark=none,mark size=1pt] table [x index=0, y index=7]{proba-vdp_noise-lan_dt1.0E-05_N1.0E+04_data.txt};
\addlegendentry{$\hat{K}_N^\varepsilon$}
\addplot[thick,dotted,color=black,mark=none,mark size=4pt] table [x index=0, y index=8]{proba-vdp_noise-lan_dt1.0E-05_N1.0E+04_data.txt};
\addlegendentry{$I^0$}
\addplot[only marks,
  black, mark options={black,scale=0.25}, 
  error bars/.cd, 
    y fixed,
    y dir=both, 
    y explicit
] table [x=x, y=y,y error=error, col sep=comma] {
    x,  y,       error
0.0875, 0.8274,  0.00740686
0.4875,  0.8489,  0.00701967
0.8875,  0.9683,  0.00343392
};
\addplot[only marks,
  blue, mark options={black,scale=0.25}, 
  error bars/.cd, 
    y fixed,
    y dir=both, 
    y explicit
] table [x=x, y=y,y error=error, col sep=comma] {
    x,  y,       error
0.1,  0.828856,  0.00410171  
0.5,  0.846056,  0.00893932  
0.9,  0.963856,  0.00804129
};
\addplot[only marks,
  red, mark options={black,scale=0.25}, 
  error bars/.cd, 
    y fixed,
    y dir=both, 
    y explicit
] table [x=x, y=y,y error=error, col sep=comma] {
    x,  y,       error
0.1125,  0.828627,  0.00393132
0.5125,  0.848279,  0.00683253  
0.9125,  0.968254,  0.0034331 
};
\end{axis}
\node at (5.5cm,3.5cm) {\textcolor{black}{(c)}};
\end{tikzpicture}
%
\begin{tikzpicture}[scale=0.5]
\begin{axis}[legend style={at={(1,1)},anchor=north east}, compat=1.3,
  xmin=0.01, xmax=1,ymin=0,ymax=3,
  xlabel= {$\varepsilon$},
  ylabel= {}]
\addplot[solid, thick, densely dotted,color=black,mark=none,mark size=1pt] table [x index=0, y index=3]{proba-vdp_noise-lan_dt1.0E-05_N1.0E+04_data.txt};
\addlegendentry{$\hat \sigma_{K^\varepsilon,N}^2/\varepsilon$}
\addplot[solid, thick, color=black,mark=none,mark size=1pt] table [x index=0, y index=4]{proba-vdp_noise-lan_dt1.0E-05_N1.0E+04_data.txt};
\addlegendentry{$\hat \sigma_{K^\varepsilon,N}^2/\varepsilon^2$}
\end{axis}
\node at (5.5cm,3.75cm) {\textcolor{black}{(d)}};
\end{tikzpicture}

\caption{Example \ref{example3b} (Van der Pol oscillator) with $h(x_1,x_2) \triangleq  x_1 - (1-x_1^2) x_2$.
In the top row the target is to estimate $I^\eps=\EE[\|\bX^\eps_T\|^2]$ for $T=1$ and 
the expectation of the control variate ${I^0=}\mathbb{E} [ \|\bU_T\|^2 ]$ is obtained by solving the PDE \eqref{eq:pde_vdp}.
In the bottom row the target is to estimate $I^\eps=\PP( | X^\eps_{1,T} | \leq 1)$ for $T=1$ and 
the expectation of the control variate ${I^0=}\PP( | X^0_{1,T} | \leq 1)$ is obtained by solving the PDE \eqref{eq:pde_vdp} with the suitable final condition.}
\label{fig:2}
\end{figure}
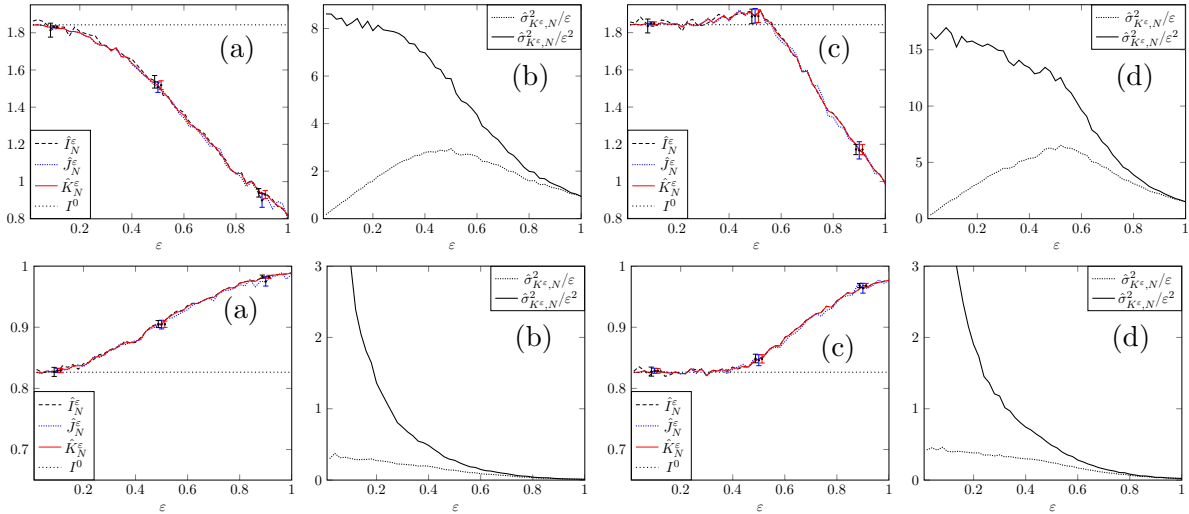

We use $N=10^4$ samples with a time step of 
\textcolor{black}{$\delta t = 10^{-5}$ (note that $\delta t /\eps^2 =0.1$ for the smallest $\eps=10^{-2}$ used  in the numerical results). }
We report the numerical results for the linear oscillator with time-dependent coefficients in Figure \ref{fig:1} 
and for the Van der Pol oscillator in Figure \ref{fig:2}.
The numerical results concern the estimation of $I^\eps = \mathbb{E} [ \|\bX^\eps_T\|^2 ]$ or $I^\eps=\PP(|X^\eps_{1,T}| \leq 1)$ 
with $T=1$ where $\bX^\eps$ satisfies \eqref{eq:smoothX} and thus the expectation of the control variate is $I^0=\EE[\| \bU_T\|^2 ]$ or
$\PP(|X^0_{1,T}| \leq 1)$  where $\bU$ satisfies \eqref{eq:smoothU}.


The theoretical predictions provided by Eq.~(\ref{eq:boundepsilon2}) and Proposition \ref{prop:21} are based on the condition that $f$ has bounded derivatives.
As we have discussed above, the assumption that $f$ is smooth is important but the hypothesis on the boundedness of the derivatives can certainly be relaxed.
The numerical results shown in Figures \ref{fig:1}-\ref{fig:2} are actually in good agreement with the theoretical predictions:
\textcolor{black}{the asymptotic
variances $\sigma^2_{J^\eps}$ and $\sigma^2_{K^\eps}$ behave as $O(\eps^2)$.}
The only cases where the behavior is $O(\eps)$, and not $O(\eps^2)$, are when the quantity of interest is of the form $\EE [ f(\bX^\eps_T)]$ with a function $f$ that is not smooth, which is not surprising.

\section{Diffusion approximation for the multivalued case}
\label{sec:adif2}
In this section we consider multivalued ODEs of the form \eqref{eq:mode1} or \eqref{eq:mode2}.

\subsection{Basic properties on differential inclusions}
\label{sec:inc}
We recall that the  subdifferential of a convex function $F : \RR^r \to (-\infty,\infty]$ such that Dom($F$)$\triangleq \{ \bx \in \RR^r , \: F(\bx) < \infty \}$ is not empty, is the map from $\RR^r$ to $\mathcal{P}(\RR^r)$ (the set of subsets of $\RR^r$) defined by $\partial F (\bx) \triangleq \{ \boldsymbol{\xi} \in \RR^r, \: \forall \bz \in \RR^r , \:  \boldsymbol{\xi}\cdot ( \bz-\bx) + F(\bx) \leq F(\bz) \}$ for $\bx \in$ Dom($F$) and $\partial F (\bx) =\emptyset$
for $\bx \not\in$ Dom($F$).
To grasp quickly the idea when $r=1$, $\partial F (x)$ can be seen as the set of sub-slopes of $F$ at the point $x$ and when $F$ is differentiable at the point $x$, $\partial F (x) = \{ F'(x)\}$. See \cite{MR0348562} for more details.

One way to construct a solution to a multivalued ODE of the form \eqref{eq:mode1} or \eqref{eq:mode2} is to proceed by penalization.
The inclusion is replaced by an equality involving the Moreau-Yosida regularisation of $F:\RR^r\to (-\infty,+\infty]$ (with $F=\varphi, r = n$ or $F=\psi, r = m$), that is
\begin{equation}
\label{def:yosida}
\forall p \geq 1, \: \forall \bx \in \RR^r , \: F_p(\bx) \triangleq \inf \limits_{\bz \in \RR^r} \left \{ F(\bz) + \frac{p}{2} \|\bx-\bz\|^2 \right \}.
\end{equation}
We recall from Annex B in \cite{MR3308895} some properties of $F_p$:
\begin{enumerate}
\item
$F_p: \RR^r \mapsto \mathbb{R}$ is a convex differentiable function,
\item
$\forall \bx \in \RR^r ,\: \partial F_p(\bx) = \{ \nabla F_p(\bx) \}$
and $\nabla F_p(\bx) \in \partial F (J_p \bx)$ where $J_p \bx \triangleq \bx - \frac{1}{p} \nabla F_p (\bx)$,
\item
$\exists C > 0, \: \forall \bx \in \RR^r, \: \forall p, \:  \| J_p \bx \| \leq \| \bx \| + C$,
\item
$\forall \bx , \by \in \RR^r, \: \|\nabla F_p(\bx) - \nabla F_p(\by) \| \leq p \|\bx-\by\|$,
\item
$\forall \bx, \by \in \RR^r , \big( \nabla F_p(\bx) - \nabla F_p(\by)\big) \cdot (\bx-\by ) \geq 0$,
\item
$\forall \bx \in \RR^r$, 
\begin{equation}
\label{eq:prop2phi}
 \bx \cdot \nabla F_p(\bx)  \geq 0,
 \end{equation}
\item
$\forall \bx ,\by \in \RR^r$, 
\begin{equation}
\label{eq:prop1phi}
\big( \nabla F_p(\bx)-\nabla F_{p'}(\by) \big) \cdot( \bx-\by ) \geq - \left ( \frac{1}{p} + \frac{1}{p'} \right ) 
 \nabla F_p(\bx) \cdot \nabla F_{p'}(\by) ,
\end{equation}
\item
as a consequence of properties 2 and 3 above, we also have 
\begin{equation}
\label{eq:condplastic}
\sup \limits_{p \geq 1} \sup \limits_{\bx \in \RR^r} \frac{ \| \nabla F_p(\bx) \|}{p (1+ \| \bx \| )} < \infty .
\end{equation}
\end{enumerate}
Thus, the penalized versions of \eqref{eq:mode1} and \eqref{eq:mode2} are 
\begin{equation}
\label{eq:pmode1}
\frac{{\rm d}\bX^{\eps,p}}{{\rm d}t} + \nabla \varphi_p (\bX^{\eps,p}) = \bb(\bX^{\eps,p}) +\frac{1}{\eps}  \bsigma \boeta^\eps  , \quad \bX_0^{\eps,p}=\bx_0
,
\end{equation}
and
\begin{equation}
\label{eq:pmode2}
\begin{dcases}
\frac{{\rm d}\bX^{\eps,p}}{{\rm d}t} + \nabla \varphi_p (\bX^{\eps,p}) = \bb^X(\bX^{\eps,p},\bZ^{\eps,p}) +\frac{1}{\eps}  \bsigma \boeta^\eps , \quad \bX_0^{\eps,p}=\bx_0,\\
\frac{{\rm d}\bZ^{\eps,p}}{{\rm d}t} + \nabla \psi_p (\bZ^{\eps,p}) = \bb^Z(\bX^{\eps,p},\bZ^{\eps,p}) , \quad \bZ_0^{\eps,p}=\bz_0.
\end{dcases}
\end{equation}
It can be shown \cite{MR0348562} that,
if $\varphi$ satisfies the condition:
\begin{equation}
\label{eq:condfriction}
 \sup \limits_{p \geq 1} \sup \limits_{\bx \in \RR^n} \| \nabla \varphi_p(\bx) \| < \infty ,
\end{equation}
where $\varphi_p$ is the Yosida approximation (\ref{def:yosida}) of $\varphi$,  
then the sequence of solutions of \eqref{eq:pmode1} $\{ \bX^{\eps,p}, p \geq 1 \}$ 
is a Cauchy sequence in $\mathcal{C}([0,T];\RR^n)$, 
the limit $\bX^{\eps}$ satisfies the differential inclusion \eqref{eq:mode1} and its solution is unique. 

A similar statement using the sequence of solutions of \eqref{eq:pmode2} $\{ (\bX^{\eps,p},\bZ^{\eps,p}), p \geq 1 \}$ in $\mathcal{C}([0,T];\RR^n \times \RR^m)$, holds for the existence and uniqueness of a solution for \eqref{eq:mode2} when $\varphi$ (but not necessarily $\psi$) 
satisfies the condition (\ref{eq:condfriction}), while $\psi$ satisfies the assumption:
\begin{equation}
\label{eq:condphipz0}
\sup \limits_{p \geq 1} \psi_p(\bz_0) < \infty.  
\end{equation}
For the convenience of the reader we give the proofs of these results in Appendix \ref{app:A}.

\subsection{Diffusion approximation for Equation \eqref{eq:mode1}}

We consider the $\RR^n$-valued process $\bX^\eps $ solution of the multivalued ODE \eqref{eq:mode1}
when $\boeta^\eps$ is given by \eqref{eq:sde2}. 
We assume that ${\itbf b}$ is Lipschitz and that $\varphi$ satisfies  the condition (\ref{eq:condfriction}).
We also consider the limiting $\RR^n$-valued process $\bX^0 $ solution of the multivalued SDE
\begin{equation}
{\rm d}\bX^0 + \partial \varphi (\bX^0) {\rm d}t \ni  \bb(\bX^0) {\rm d} t + \bGamma  {\rm d} \bW_t   ,
\label{eq:msde_limit1}
\end{equation}
driven by the same Brownian motion, with
$\bGamma= \bsigma  {\bf A}^{-1} {\bf K}$.
Existence and uniqueness of the solution of \eqref{eq:msde_limit1} is the same one as in Proposition \ref{thmA1} and is discussed in Appendix \ref{app:D}.
The following proposition gives the convergence of the process $(\bX^\eps-\bX^0)$
to zero. It is proved in Appendix \ref{app:proofprop1}.

\begin{prop}
\label{prop:1}
{\it 1.}
We have for all $p \geq 1$:
\begin{equation}
\label{eq:lem1}
 \sup \limits_{\eps} \mathbb{E} \left[ \sup \limits_{ t \leq T} \| \bX_t^{\eps,p} - \bX_t^\eps \|^2\right] \leq \frac{C_T}{p} ,
\end{equation}
where $\bX^{\eps,p}$ is the approximation (\ref{eq:pmode1}) of $\bX^\eps$.\\
{\it 2.}
The continuous process $(\bX^\eps-\bX^0)$  converges in probability to zero
as $\eps\to 0$.
\end{prop}

\subsection{Diffusion approximation for Equation \eqref{eq:mode2}}

We consider the $\RR^n \times \RR^m$-valued process $(\bX^\eps,\bZ^\eps)$ solution of the multivalued ODE \eqref{eq:mode2}
when $\boeta^\eps$ is given by \eqref{eq:sde2}. 
We assume that ${\itbf b}^X$ and ${\itbf b}^Z$ are Lipschitz, that $\varphi$ satisfies \eqref{eq:condfriction}, and that 
$\psi$ satisfies (\ref{eq:condphipz0}).
We also consider the limiting $\RR^n \times \RR^m$-valued process $(\bX^0,\bZ^0)$ solution of the multivalued SDE
\begin{equation}
{\rm d}\bX^0 + \partial \varphi (\bX^0) {\rm d}t \ni  \bb^X(\bX^0,\bZ^0) {\rm d} t + \bGamma  {\rm d} \bW_t,
\quad {\rm d} \bZ^0 + \partial \psi (\bZ^0) {\rm d}t \ni \bb^Z(\bX^0,\bZ^0) {\rm d} t, 
\label{eq:msde_limit}
\end{equation}
driven by the same Brownian motion, with
$\bGamma= \bsigma  {\bf A}^{-1} {\bf K} $.
Existence and uniqueness of the solution of \eqref{eq:msde_limit} is the same one as in Proposition \ref{thmA2}.
The following proposition gives the convergence of the  continuous process $(\bX^\eps-\bX^0,\bZ^\eps-\bZ^0)$ to zero.
It is proved in Appendix \ref{app:proofprop2}.

\begin{prop}
\label{prop:2}
{\it 1.}
For all $p\geq 1$, we have 
\begin{equation}
 \sup \limits_{\eps} \mathbb{E} \left [ \sup \limits_{ t \leq T} \left \{ \| \bX_t^{\eps,p} - \bX_t^\eps \|^2 + \| \bZ_t^{\eps,p} - \bZ_t^\eps \|^2\right  \} \right]
\leq \frac{C_T}{p} ,
\label{eq:lem2}
\end{equation}
where $(\bX^{\eps,p},\bZ^{\eps,p})$ is the approximation (\ref{eq:pmode2}) of $(\bX^\eps,\bZ^\eps)$.\\
{\it 2.}
The continuous process $(\bX^\eps-\bX^0,\bZ^\eps-\bZ^0)$ converges in probability to zero
as $\eps\to 0$.
\end{prop}

\section{Control variate method in the multivalued case}
\label{sec:cvmulti}
We here consider the multivalued case. 
Let $\bX^{\eps}$ satisfy (\ref{eq:mode1}) or $(\bX^{\eps},\bZ^{\eps})$ satisfy (\ref{eq:mode2}).
\textcolor{black}{
We want to estimate $I^\eps$ defined by (\ref{def:Ieps}) when $\bX^{\eps}$ satisfies (\ref{eq:mode1}) or 
$I^\eps=\EE[ F(\bX^{\eps},\bZ^{\eps})]$
when $(\bX^{\eps},\bZ^{\eps})$ satisfies (\ref{eq:mode2}).
The control variate method can be applied in this framework as in the ODE case addressed in Section \ref{sec:ode}.
The control variate estimator $\hat{J}^\eps_N$ and the optimal control variate estimator $\hat{K}^\eps_N$ are defined by
(\ref{def:IepsN}) and  (\ref{def:KepsN}), respectively, 
for $\bX^{\eps}$ satisfying (\ref{eq:mode1}), they are asymptotically normal and their asymptotic variances are (\ref{def:sigmaIeps}) and (\ref{def:sigmaKeps}), respectively.
For $(\bX^{\eps},\bZ^{\eps})$ satisfying (\ref{eq:mode2}) the control variate estimator 
\begin{equation}
\hat{J}_N^\eps  \triangleq \frac{1}{N}\sum_{k=1}^N F (\bX^\eps(\bW^k), \bZ^\eps(\bW^k) ) - F (\bX^0 (\bW^k),\bZ^0(\bW^k) )  +I^0,
\end{equation}
with $I^0=\EE[F(\bX^0,\bZ^0)]$,
is asymptotically normal with an asymptotic variance given by
\begin{equation}
\label{eq:prop5:0}
\sigma^2_{J^\eps}= {\rm Var}\big( F (\bX^\eps , \bZ^\eps  ) - F(\bX^0  ,\bZ^0  )\big)  .
\end{equation}
The pratical optimal control variate estimator is
\begin{equation}
\label{def:KepsNXZ}
\hat{K}^{\eps}_N \triangleq  \hat{\rho}^\eps_N I^0+ \frac{1}{N} \sum_{k=1}^N F( \bX^\eps (\bW^{k}), \bZ^\eps(\bW^k) ) - \hat\rho^\eps_N F( \bU(\bW^{k}), \bZ^0(\bW^k) ), 
\end{equation}
where $\hat\rho^\eps_N$ is the empirical correlation 
\begin{equation}
\hat\rho^\eps_N = \frac{
 \sum_{k=1}^N ( F( \bX^\eps(\bW^{k}), \bZ^\eps(\bW^k) ) - \hat{I}^\eps_N) (F( \bU(\bW^{k}), \bZ^0(\bW^k)) - \hat{I}^0_N)
 }
 {
  \sum_{k=1}^N ( F(\bU(\bW^{k}), \bZ^0(\bW^k)) - \hat{I}^0_N )^2
},
\label{def:hatrhoNepsXZ}
\end{equation}
$\hat{I}^\eps_N$ is the standard Monte Carlo estimator
\begin{equation}
\hat{I}^\eps_N  \triangleq  \frac{1}{N}  \sum_{k=1}^N F( \bX^\eps(\bW^{k}), \bZ^\eps(\bW^k))  ,
\end{equation}
and
$ \hat{I}^0_N= \frac{1}{N}  \sum_{k=1}^N F( \bX^0(\bW^{k}), \bZ^0(\bW^k))$.
The estimator $\hat{K}^\eps_N$ is asymptotically normal with an asymptotic variance given by
\begin{equation}
\sigma^2_{K^\eps}= {\rm Var}\big( F(\bX^\eps , \bZ^\eps) - \rho^\eps F(\bX^0  ,\bZ^0  )\big)  ,
\end{equation}
with $\rho^\eps = {\rm Cov}(F (\bX^\eps , \bZ^\eps  ) , F (\bX^0  ,\bZ^0  ) )/{\rm Var}(  F (\bX^0  ,\bZ^0  ))$.
The asymptotic variance of the estimator $\hat{I}_N^\eps$ is
\begin{equation}
\sigma^2_{I^\eps}= {\rm Var}\big( F (\bX^\eps , \bZ^\eps  ) \big)  .
\end{equation}
For small $\eps$
the asymptotic variance  $\sigma^2_{I^\eps}$ is approximately equal to $ {\rm Var}\big( F (\bX^0, \bZ^0  ) \big)$ and the asymptotic variances $\sigma^2_{J^\eps}$ and $\sigma^2_{K^\eps}$  are  small, 
by Propositions \ref{prop:1} and \ref{prop:2}.
More quantitiatvely, if  $F(\bX,\bZ)=f(\bX_T,\bZ_T)$ and $f$ is a smooth function with bounded derivatives, then  the asymptotic variances $\sigma^2_{J^\eps}$ and $\sigma^2_{K^\eps}$  are
of order $\eps^2$ when  $\bX^{\eps}$ satisfies (\ref{eq:mode1}):
\begin{equation}
\label{eq:prop5:1}
\sigma^2_{K^\eps} \leq C \eps^2,\quad \quad
0\leq \sigma^2_{J^\eps} - \sigma^2_{K^\eps} \leq C \eps^4 ,
\end{equation}
or of order $\eps$ when $(\bX^{\eps},\bZ^{\eps})$ satisfies (\ref{eq:mode2}):
\begin{equation}
\label{eq:prop5:2}
\sigma^2_{K^\eps} \leq C \eps,\quad \quad
0\leq \sigma^2_{J^\eps} - \sigma^2_{K^\eps} \leq C \eps^2.
\end{equation}
Eqs.~(\ref{eq:prop5:1}-\ref{eq:prop5:2}) are consequences of the following lemma proved in Appendix \ref{app:proofprop5}.
}

\begin{lem}
\label{prop:5}
\begin{enumerate}
\item
Let $\bX^{\eps}$ satisfy (\ref{eq:mode1}). If $f,g:\RR^n\to\RR$  are smooth functions with bounded derivatives and $T>0$,
then there exists $C>0$ such that, for any $t \in [\eps,T]$,
\begin{equation}
\big| \EE \big[  g (\bZ_t^0  ) \big( f (\bX^\eps_t   ) - f (\bU_t   )\big) \big] \big| \leq C \eps^2  ,
\quad \quad 
\EE \big[   \big( f (\bX^\eps_t   ) - f (\bU_t   )\big)^2 \big] \leq C \eps^2  .
\end{equation}

\item
Let $(\bX^{\eps},\bZ^{\eps})$ satisfy (\ref{eq:mode2}). If $f,g:\RR^{n+m}\to\RR$ are smooth functions with bounded derivatives and $T>0$,
then there exists $C>0$ such that, for any $t \in [\eps,T]$,
\begin{equation}
\big| \EE \big[  g (\bX_t^0  ,\bZ_t^0  ) \big( f (\bX_t^\eps , \bZ_t^\eps  ) - f (\bX_t^0  ,\bZ_t^0  )\big)\big]\big| \leq C \eps ,\quad \quad
\EE \big[ \big( f (\bX_t^\eps , \bZ_t^\eps  ) - f (\bX_t^0  ,\bZ_t^0  )\big)^2\big] \leq C \eps .
\end{equation}
\end{enumerate}
\end{lem}

\section{Numerical simulations in the multi-valued case}
\label{sec:nummulti}
We present examples which are non-smooth dynamical systems that are prevalent in engineering mechanics.
Examples \ref{example81a}-\ref{example81b}
are oscillators involving friction or/and elasto-plastic behaviours,
Exemple \ref{example81c} is a nonlinear and nonsmooth two-degree of freedom (TDOF) system,
they can be described by Eqs.~\eqref{eq:mode1} and \eqref{eq:mode2}.
Examples \ref{example82a}-\ref{example82b} 
which do not fall within the scope of any aforementioned case correspond to an obstacle problem and to the reflection of the integral of a colored noise, respectively.

\subsection{Non-smooth systems in the form of Equations \eqref{eq:mode1} and \eqref{eq:mode2}}
\label{subsubsec:mod2}


\begin{example}[\textbf{friction behaviour}]
\label{example81a}%
With Equation \eqref{eq:mode1} in mind, we take $\forall x \in \RR$,  $\varphi(x) \triangleq c_{\rm f} |x|$ where $c_{\rm f}>0$ is a friction coefficient.
The $\RR$-valued process $X^\eps$ represents the velocity of a material point (stick-slip motion) subjected to friction and colored noise.
See for instance \cite{MITVIDEOS} for an explanation of the physics behind and \cite{bernardin} for the use of SDEs with multivalued drift for modeling.
As $\eps \to 0$, $\bX^\eps \to \bU$ where $\bU$ satisfies Equation \eqref{eq:msde_limit1}.
For the stochastic simulation, we proceed as follows:
\begin{itemize}
\item
\textcolor{black}{
$\hat X_0^{\eps,{k}} = x_0$ and for $0 \leq n \leq N_T-1$,
$$
\hat X_{n+1}^{\eps, {k}} 
=  \hat X_n^{\eps, {k}} + \frac{\delta t}{\eps} \hat{\eta}^{\eps,{k}}_n -  \delta t \textup{proj}_{[-c_{\rm f}, c_{\rm f}]} \left ( \frac{\hat X_n^{\eps, {k}}}{\delta t} + \frac{1}{\eps}  \hat{\eta}^{\eps,{k}}_n \right )   ,
$$
with $ \hat{\eta}^{\eps,{k}}_n$ described in Subsection \ref{subsec:num1}.
}
\item
$\hat X^{0,{k}}_0 = x_0$ and for $0 \leq n \leq N_T-1$,
$$
\hat X^{0,{k}}_{n+1} 
=  \hat X^{0,{k}}_n + (\delta t)^{1/2} C \Delta W_n^{k} -  \delta t \textup{proj}_{[-c_{\rm f}, c_{\rm f}]} \left ( \frac{\hat X^{0,{k}}_n}{\delta t} + 
\frac{C}{\sqrt{\delta t}} \Delta W_n^{k}  \right )   .
$$ 
\end{itemize}
We are interested in $\mathbb{E} [ \left ( X_T^\eps \right )^2 ]$ for $T=1$. The expectation of the control variate is $\mathbb{E} [(X^0_T)^2]$. The latter can be represented as $\mathfrak{c}(x_0,0)$ where $\mathfrak{c}$ satisfies the following backward in time partial differential inclusion
\begin{equation}
\label{eq:pdi_friction}
\begin{cases}
& \partial_t \mathfrak{c}(x,t) + \frac{C^2}{2}\partial_x^2  \mathfrak{c}(x,t)  \in \partial \varphi(x)  \partial_x \mathfrak{c}(x,t) = 0, \: \: \textup{for} \: \: (x,t) \in  \mathbb{R} \times [0,1),\\
& \mathfrak{c}(x,1) = x^2,\: \: \textup{for} \: \: x \in \RR.
\end{cases}
\end{equation}
It can be estimated by solving this  partial differential inclusion with a finite difference method.
We proceed as follows. For every $t>0$, the function $x \mapsto \mathfrak{c}(x,t)$ is smooth and even, provided that the initial condition is smooth and even.
Indeed, this comes from the probabilistic representation and the fact that, for any starting point $x \in \mathbb{R}, \: \{ X_t^x, t \geq 0 \}$ and $\{ X_t^{-x}, t \geq 0 \}$ have the same distribution because $\varphi$ is even.
Therefore we must have $\forall t>0, \: \partial_x \mathfrak{c}(0,t) = 0$. The solution of \eqref{eq:pdi_friction} is thus estimated by applying a finite difference method to 
\begin{equation}
\label{eq:pde_friction}
\begin{cases}
& \partial_t \mathfrak{c}(x,t) + \frac{C^2}{2}\partial_x^2  \mathfrak{c}(x,t) - c_{\rm f}  \partial_x \mathfrak{c}(x,t) = 0 , \: \: \textup{for} \: \: (x,t) \in  (0,\infty) \times [0,1),\\
& \partial_x \mathfrak{c}(0,t) = 0,\: \: \textup{for} \: \: t \in [0,1),\\
& \mathfrak{c}(x,1) = x^2,\: \: \textup{for} \: \: x \in [0,+\infty).
\end{cases}
\end{equation}
The whole function $x \mapsto \mathfrak{c}(x,t)$ can be recovered by using the symmetry property. 
\end{example}

\begin{example}[\textbf{elasto-plastic behaviour}]
\label{example81b}
With Eq.~\eqref{eq:mode2}, we consider $\varphi \triangleq 0$, $\psi \triangleq \chi_D$ the indicator function of $D \triangleq [-c_{\rm ep},c_{\rm ep}]$ in the sense of convex analysis, that is $\chi_D (x)= 0$ if $x \in D$ and $+\infty$ otherwise. Here $c_{\rm ep}>0$ is an elasto-plastic coefficient. 
The real-valued process $\bX^\eps$ represents the velocity of a material point subjected to an elasto-plastic restoring force and colored noise.
The process $\bZ^\eps$ taking values in $[-c_{\rm ep},c_{\rm ep}]$ represents the restoring force.
See for instance \cite{MYVIDEO} for an explanation of the physics and the use of SDEs with multivalued drift for modeling.
Here we are interested in $\EE[ (\bX^\eps_T)^2 + (\bZ^\eps_T)^2]$ and in $\PP \left ( | \bZ^\eps_T | = c_{\rm ep} \right )$ for $T=1$. 
As $\eps \to 0$, $(\bX^\eps,\bZ^\eps) \to (\bX^0,\bZ^0)$ where $(\bX^0,\bZ^0)$ satisfies  \eqref{eq:msde_limit}.
For the stochastic simulation, we proceed as follows:
\begin{itemize}
\item
$\hat X_0^{\eps,{k}} = x_0$ and $\hat Z_0^{\eps,{k}} = z_0$ and for $0 \leq n \leq N_T-1$,
$$
\begin{cases}
\hat Z_{n+1}^{\eps, {k}} = \textup{proj}_{[-c_{\rm ep},c_{\rm ep}]} \left ( \hat Z_n^{\eps, {k}} + \delta t  \hat X_n^{\eps, {k}} \right) ,\\
\hat X_{n+1}^{\eps, {k}} = \hat X_n^{\eps, {k}} - \delta t \hat Z_n^{\eps, {k}} + \frac{\delta t}{\eps} \hat \eta_n^{\eps, {k}}.
\end{cases}
$$
\item
$\hat X_0^{0,{k}} = x_0$ and $\hat Z_0^{0,{k}} = z_0$ and for $0 \leq n \leq N_T-1$,
$$
\begin{cases}
\hat Z_{n+1}^{0, {k}} = \textup{proj}_{[-c_{\rm ep},c_{\rm ep}]} \left ( \hat Z_n^{0, {k}} + \delta t  \hat X_n^{0, {k}} \right) , \\
\hat X_{n+1}^{0, {k}} = \hat X_n^{0, {k}} - \delta t \hat Z_n^{0, {k}} + \sqrt{\delta t} C \Delta W_n^{k}.
\end{cases}
$$
\end{itemize}
The expectations of the control variates are $\EE [ (\bX^0_T)^2 +(\bZ^0_T)^2]$ and  $\mathbb{P} ( | \bZ^0_T | = c_{\rm ep} )$ for $T=1$. They are estimated using the PDE method of \cite{msw19}.
\end{example}

\begin{example}[\textbf{nonlinear and nonsmooth TDOF}]
\label{example81c}
A TDOF with an elasto-plastic element can be represented as a system of the form \eqref{eq:mode1} which becomes
\begin{equation}
\label{eq:nonsmoothXcouple}
\begin{cases}
\dfrac{{\rm d} X_1^\eps}{{\rm d}t} + \partial \chi_D(X_1^\eps) = X_2^\eps, 
\quad \dfrac{{\rm d} X_2^\eps}{ {\rm d}t} = - \tilde{g}_2 (X_1^\eps ,X_2^\eps,X_3^\eps ,X_4^\eps) + \dfrac{1}{\eps} \eta^\eps, \\[2mm]
\dfrac{{\rm d} X_3^\eps}{{\rm d}t} = X_4^\eps,
\quad \dfrac{{\rm d} X_4^\eps}{ {\rm d}t} = - \tilde{g}_4(X_1^\eps ,X_2^\eps,X_3^\eps ,X_4^\eps).
\end{cases}
\end{equation}
Here $D = [-c_p,c_p], c_p = 0.25$.
In addition to $\EE \big[(X_{1,T}^\eps )^2+ (X_{3,T}^\eps)^2 \big]$ for $T=1$, we are interested in 
$\PP( | X^\eps_{1,T} | \leq a, | X^\eps_{3,T} | \leq b)$.
As $\eps \to 0$, $\bX^\eps \to \bU$ where
\begin{equation}
\label{eq:nonsmoothUcouple}
\begin{cases}
{\rm d} X^0_1 + \partial \chi_D(X^0_1) {\rm d} t = X^0_2  {\rm d} t,
\quad {\rm d} X^0_2 = - \tilde{g}_2 (X^0_1,X^0_2,X^0_3,X^0_4) {\rm d} t + C {\rm d}W, \\[2mm]
{\rm d} X^0_3 = X^0_4  {\rm d} t,
\quad {\rm d} X^0_4 = - \tilde{g}_4 (X^0_1,X^0_2,X^0_3,X^0_4)  {\rm d} t.
\end{cases}
\end{equation}
Here, for simplicity the other elements are linear,  
$$
\tilde{g}_2(x_1,x_2,x_3,x_4) \triangleq 
k_1 x_1 + c_1 x_2
- k_3 (x_3-x_1)
- c_3(x_4-x_2) 
$$
and 
$$
\tilde{g}_4(x_1,x_2,x_3,x_4) \triangleq 
k_3 (x_3-x_1) + c_3(x_4-x_2).
$$
The simulation of \eqref{eq:nonsmoothXcouple} and \eqref{eq:nonsmoothUcouple} is similar to what is explained above.
We take $c_1=c_3=k_1=k_3=1$.
\end{example}

\subsection{Non-smooth systems : beyond Equations \eqref{eq:mode1} and \eqref{eq:mode2}}
\label{subsubsec:mod3}
The two models presented in this subsection do not fall in the scope of our theoretical results, though they are not too far off. The presentation of the impact problem remains formal.
The behaviour of the control variate estimator is investigated via numerical experiments.

\begin{example}[\textbf{impact problem}] 
\label{example82a}
The pair displacement-velocity $\bX^\eps=(X_1^\eps,X_2^\eps)$ (taking values in $\RR^2$) of a colored noise driven oscillator constrained by an obstacle can be formulated 
in terms of an equation of the form \eqref{eq:smoothX} when $|X_{1,t}^\eps| < P_{\rm O}$ with the condition (that expresses the switch of the velocity at collision): for all $t$, 
$|X_{1,t}^\eps| = P_{\rm O} \implies  X_{2,t+}^\eps = - \mathfrak{e} X_{2,t-}^\eps$ 
where $P_{\rm O}$ is the location of the obstacle and $\mathfrak{e} \in [0,1]$ is the coefficient of restitution of energy.
The notations $X_{2,t\pm}^\eps$ stand for the velocity immediately before and after the collision.
Here we are interested in $\mathbb{E} [ (X_{2,T}^\eps)^2 ]$ for $T=1$.
Formally, as $\eps \to 0$, the $\RR^2$-valued limit process $\bU=(X^0_1,X^0_2)$ is a white noise driven oscillator constrained by an obstacle that can be formulated similarly to the former case, 
except that  we replace \eqref{eq:smoothX} by \eqref{eq:smoothU}. When $\mathfrak{e} = 1$ (resp. $0 \leq \mathfrak{e} < 1$), we say that the collisions are elastic (resp. inelastic). 
It is important to stress that obstacle problems with inelastic collisions deserve more attention for practical purposes since in real world phenomena kinetic energy
is dissipated through heat or plastic deformation. With elastic collisions, there is no loss of kinetic energy.
For the stochastic simulation, we use the same numerical procedure as for \eqref{eq:smoothX} and \eqref{eq:smoothU}, 
except that if we find out that the $(n+1)^{st}$ point does not satisfy the obstacle condition, 
i.e. $|\hat X_{1,n+1}^{\eps,{k}}| > P_{\rm O}$,
we adjust the time step length to $\theta_{n+1} \delta t$
with $\theta_{n+1} \triangleq \frac{\pm P_{\rm O} - \hat X_{1,n}^{\eps,{k}}}{\hat X_{1,n+1}^{\eps,{k}}-\hat X_{1,n}^{\eps,{k}}}$
and set $t_{n+1} \triangleq t_n + \theta_{n+1} \delta t, \hat X_{1,n+1}^{\eps,{k}} \triangleq P_{\rm O}$, 
$\hat X_{2,n+1}^{\eps,{k}} \triangleq - \mathfrak{e} \left (  
\hat X_{2,n}^{\eps, {k}} 
- \theta_{n+1} \delta t  f ( \hat X_{1,n}^{\eps, {k}}, \hat X_{2,n}^{\eps, {k}} )
+ \theta_{n+1} \frac{\delta t}{\eps} \hat \eta_n^{\eps, {k}} \right )$
and 
$\hat \eta_{2,n+1}^{\eps,{k}} \triangleq \hat \eta_n^{\eps, {k}} \left (1 - \theta_{n+1} \delta t \dfrac{A}{\eps^2} \right ) + \sqrt{\theta_{n+1} \delta t} \dfrac{K_{\textup{ou}}}{\eps} \Delta W_n^{k}$.
A similar adjustment is done in the other cases with Langevin and white noises. 
The expectation of the control variate is $\mathbb{E} [ (X^0_{2,T})^2 ]$ for $T=1$ which is estimated using the PDE method of \cite{msw19}.
\end{example}

\begin{example}[\textbf{reflection of an integrated colored noise}]
\label{example82b}
Define $E \triangleq [0,\infty)$ and consider the indicator function of $E$, that is $\chi_E( x)=0$ if $x \in E$ and $+\infty$ otherwise.
The reflection of an integrated colored noise corresponds to the case where $\bX^\eps$ satisfies
\begin{equation}
\label{eq:ref_icn}
\frac{{\rm d} X^\eps}{{\rm d}t} + \partial \chi_E (X^\eps) \ni \frac{1}{\eps} \eta^\eps ,
\end{equation}
and $\bU$, the limit process as $\eps \to 0$, is a reflected Brownian motion  
\begin{equation}
\label{eq:ref_bm}
{\rm d} X^0 + \partial \chi_E (X^0) {\rm d}t  \ni C {\rm d} W. 
\end{equation}
We are interested in $\mathbb{E} [X^\eps_T ]$ for $T=1$.
For the stochastic simulation of \eqref{eq:ref_icn} and \eqref{eq:ref_bm}, we use the following scheme: $\hat X_0^{\eps,{k}} = x_{0}$, 
$\hat X^{0,{k}} = x_{0}$ 
and for $0 \leq n \leq N_T-1$, 
\begin{itemize}
\item
$\hat X_{n+1}^{\eps,{k}} = \textup{proj}_E \left ( \hat X_n^{\eps,{k}} + \frac{\delta t}{\eps} \hat \eta_n^{\eps,{k}} \right )$,
\item
$\hat X_{n+1}^{0,{k}} = \textup{proj}_E \left ( \hat X_n^{0,{k}} + \sqrt{\delta t} C \Delta W_n^{k} \right )$.
\end{itemize}
The expectation of the control variate is given by an explicit formula $\EE [X^0_1] = \sqrt{{2}/{\pi}}$. 
Indeed, the backward Kolmogorov equation for the reflected Brownian motion in \eqref{eq:ref_bm} is
$$
\partial_t w =  C \partial_{x}^2 w , \: x > 0, \: t > 0,
\quad w(x,t = 0) = x, \: x >0,
\quad \partial_x w(0,t) = 0, \: t >0.
$$
It has an explicit solution
$$
w(x,t) = \frac{1}{\sqrt{4 C \pi t}} \int_0^\infty y \left ( \exp\big(  -\dfrac{(x-y)^2}{4 C t} \big) + \exp\big(  -\dfrac{(x+y)^2}{4 C t}\big) \right ) {\rm d} y,
$$
which gives $ \mathbb{E} [X^0_T]=w(0,T) = \sqrt{{2}/{\pi}}$ for $T=1$.
In this case, we can provide an ad hoc proof to get an estimate similar to \eqref{eq:estimeffeps2} (see Appendix \ref{app:C}):
\begin{equation}
\label{eq:estimadhoc}
\EE \big[ (X^\eps_T - X^0_T)^2 \big] \leq C \eps^2 |\log \eps| .
\end{equation}
The $\log \eps$ correction comes from a maximal inequality for the OU process \cite{MR1664394} and a standard result on the maxima of Gaussian processes \cite{massart}.
\end{example}

\subsection{Numerical experiments}
\label{subsec:num3multi}%
We report our numerical results for the four systems mentioned above.
The convention is as in Subsection \ref{subsec:num3}.
In each of the four figures below, there are four subfigures (a)-(b)-(c)-(d). For subfigures  (a) and (b), the driving force is an Ornstein-Uhlenbeck noise \eqref{eq:ou} with $A=K=1$. 
\textcolor{black}{
In subfigure (a), the dashed black lines, the dotted blue lines, and the solid red lines represent 
 the standard MC estimator $\hat{I}_N^\eps $ and the control variate estimators $\hat{J}_N^\eps$  and $\hat{K}_N^\eps$, respectively.}
The dotted black line represents the expectation of the control variate $I^0$. The objective of the subfigure (b) is to illustrate the bounds 
(\ref{eq:prop5:1}) and (\ref{eq:prop5:2}).
The same description applies to (c) and (d), except they correspond to the case of a Langevin noise (\ref{eq:lang1a}-\ref{eq:lang1b})
with $\mu=\gamma=K=1$.
\textcolor{black}{
For the examples in which $\bX^{\eps}$ satisfies (\ref{eq:mode1}), 
  the asymptotic variance the standard MC estimator  $\hat{I}^\eps_N$ is estimated by (\ref{eq:estimvarJneps}),
 the asymptotic variance of the control variate estimator $\hat{J}^\eps_N$ is estimated by (\ref{eq:estimvarIneps}),
  and
 the asymptotic variance of the optimal control variate estimator $\hat{K}^\eps_N$ is estimated by (\ref{eq:estimvarKneps}).
For the examples in which $(\bX^{\eps},\bZ^{\eps})$ satisfies (\ref{eq:mode2}),
the asymptotic variance of the standard MC estimator $\hat{I}^\eps_N$ is estimated by
\begin{equation}
\widehat{\sigma}_{I^\eps,N}^2  =\frac{1}{N} 
\sum_{k=1}^N \big( F(\bX^\eps(\bW^k), \bZ^\eps(\bW^k) )  ) \big)^2  - ({\hat{I}^\eps_N})^2 ,
\end{equation}
the asymptotic variance of the control variate estimator $\hat{J}^\eps_N$ is estimated by
\begin{equation}
\widehat{\sigma}_{J^\eps,N}^2  =  \frac{1}{N} 
\sum_{k=1}^N \big( F(\bX^\eps(\bW^k), \bZ^\eps(\bW^k) ) - F (\bX^0 (\bW^k),\bZ^0(\bW^k) )  +I^0\big)^2 - (\hat{J}^\eps_N)^2 ,
\end{equation}
and
the asymptotic variance of the optimal control variate estimator $\hat{K}^\eps_N$ is estimated by
\begin{equation}
\widehat{\sigma}_{K^\eps,N}^2  =  \frac{1}{N} 
\sum_{k=1}^N \big( F (\bX^\eps(\bW^k), \bZ^\eps(\bW^k) ) -  \hat{\rho}_N^\eps F ( \bU (\bW^k) , \bZ^0(\bW^k))  +\hat{\rho}_N^\eps I^0 \big)^2 - (\hat{K}^\eps_N)^2 ,
\end{equation}
with $\hat{\rho}_N^\eps$ defined by (\ref{def:hatrhoNepsXZ}).
$\widehat{\sigma}_{I^\eps,N}^2$, $\widehat{\sigma}_{J^\eps,N}^2 $, and $\widehat{\sigma}_{K^\eps,N}^2 $ are consistent estimators of 
${\sigma}_{I^\eps}^2$, ${\sigma}_{J^\eps}^2$, and ${\sigma}_{K^\eps}^2$, 
respectively.}

\begin{figure}[h!]
\centering
\begin{tikzpicture}[scale=0.5]
\begin{axis}[legend style={at={(0,0)},anchor=south west}, compat=1.3,
  xmin=0.01, xmax=1,ymin=0.0,ymax=0.35,
  xlabel= {$\varepsilon$},
  ylabel= {Ornstein-Uhlenbeck}]
\addplot[densely dashed,color=black,mark=none,mark size=1pt] table [x index=0, y index=5]{friction_noise-ou_dt1.0E-05_N1.0E+04_data.txt};
\addlegendentry{$\hat{I}_N^\varepsilon$}
\addplot[densely dotted,thick,color=blue,mark=none,mark size=1pt] table [x index=0, y index=6]{friction_noise-ou_dt1.0E-05_N1.0E+04_data.txt};
\addlegendentry{$\hat{J}_N^\varepsilon$}
\addplot[solid,color=red,mark=none,mark size=1pt] table [x index=0, y index=7]{friction_noise-ou_dt1.0E-05_N1.0E+04_data.txt};
\addlegendentry{$\hat{K}_N^\varepsilon$}
\addplot[thick,dotted,color=black,mark=none,mark size=4pt] table [x index=0, y index=8]{friction_noise-ou_dt1.0E-05_N1.0E+04_data.txt};
\addlegendentry{$I^0$}
\addplot[only marks,
  black, mark options={black,scale=0.25}, 
  error bars/.cd, 
    y fixed,
    y dir=both, 
    y explicit
] table [x=x, y=y,y error=error, col sep=comma] {
    x,  y,       error
0.0875,  0.324519,  0.0116878   
0.4875,  0.124403,  0.0065813 
0.8875,  0.0215634,  0.0019716
};
\addplot[only marks,
  blue, mark options={black,scale=0.25}, 
  error bars/.cd, 
    y fixed,
    y dir=both, 
    y explicit
] table [x=x, y=y,y error=error, col sep=comma] {
    x,  y,       error
0.1,  0.327004,  0.00171263
0.5,  0.119132,  0.00882551  
0.9,  0.0314939,  0.0112926
};
\addplot[only marks,
  red, mark options={black,scale=0.25}, 
  error bars/.cd, 
    y fixed,
    y dir=both, 
    y explicit
] table [x=x, y=y,y error=error, col sep=comma] {
    x,  y,       error
0.1125,  0.326926,  0.00167046  
0.5125,  0.122364,  0.00452465 
0.9125,  0.0220548,  0.00188444 
};
\end{axis}
\node at (5.5cm,4.0cm) {\textcolor{black}{(a)}};
\end{tikzpicture}
\begin{tikzpicture}[scale=0.5]
\begin{axis}[legend style={at={(1,1)},anchor=north east}, compat=1.3,
  xmin=0.01, xmax=1,ymin=0,ymax=0.875,
  xlabel= {$\varepsilon$},
  ylabel= {}]
\addplot[solid, thick, densely dotted,color=black,mark=none,mark size=1pt] table [x index=0, y index=3]{friction_noise-ou_dt1.0E-05_N1.0E+04_data.txt};
\addlegendentry{$\hat \sigma_{K^\varepsilon,N}^2/\varepsilon$}
\addplot[solid, thick, color=black,mark=none,mark size=1pt] table [x index=0, y index=4]{friction_noise-ou_dt1.0E-05_N1.0E+04_data.txt};
\addlegendentry{$\hat \sigma_{K^\varepsilon,N}^2/\varepsilon^2$}
\end{axis}
\node at (5.5cm,3.5cm) {\textcolor{black}{(b)}};
\end{tikzpicture}
\centering
\begin{tikzpicture}[scale=0.49]
\begin{axis}[legend style={at={(0,0)},anchor=south west}, compat=1.3,
  xmin=0.01, xmax=1,ymin=0.0,ymax=0.35,
  xlabel= {$\varepsilon$},
  ylabel= {Langevin}]
\addplot[densely dashed,color=black,mark=none,mark size=1pt] table [x index=0, y index=5]{friction_noise-lan_dt1.0E-05_N1.0E+04_data.txt};
\addlegendentry{$\hat{I}_N^\varepsilon$}
\addplot[densely dotted,thick,color=blue,mark=none,mark size=1pt] table [x index=0, y index=6]{friction_noise-lan_dt1.0E-05_N1.0E+04_data.txt};
\addlegendentry{$\hat{J}_N^\varepsilon$}
\addplot[solid,color=red,mark=none,mark size=1pt] table [x index=0, y index=7]{friction_noise-lan_dt1.0E-05_N1.0E+04_data.txt};
\addlegendentry{$\hat{K}_N^\varepsilon$}
\addplot[thick,dotted,color=black,mark=none,mark size=4pt] table [x index=0, y index=8]{friction_noise-lan_dt1.0E-05_N1.0E+04_data.txt};
\addlegendentry{$I^0$}
\addplot[only marks,
  black, mark options={black,scale=0.25}, 
  error bars/.cd, 
    y fixed,
    y dir=both, 
    y explicit
] table [x=x, y=y,y error=error, col sep=comma] {
    x,  y,       error
0.0875,  0.333822,  0.0123062   
0.4875,  0.21652,  0.0103896
0.8875,   0.0366184,  0.00314677
};
\addplot[only marks,
  blue, mark options={black,scale=0.25}, 
  error bars/.cd, 
    y fixed,
    y dir=both, 
    y explicit
] table [x=x, y=y,y error=error, col sep=comma] {
    x,  y,       error
0.1,  0.336495,  0.00243837
0.5,  0.220044,  0.0118184  
0.9,  0.044319,  0.0120353
};
\addplot[only marks,
  red, mark options={black,scale=0.25}, 
  error bars/.cd, 
    y fixed,
    y dir=both, 
    y explicit
] table [x=x, y=y,y error=error, col sep=comma] {
    x,  y,       error
0.1125,  0.336438,  0.00242403 
0.5125,  0.217917,  0.00915955 
0.9125,  0.0367544,  0.00313984 
};
\end{axis}
\node at (5.5cm,4.0cm) {\textcolor{black}{(c)}};
\end{tikzpicture}
\begin{tikzpicture}[scale=0.49]
\begin{axis}[legend style={at={(1,1)},anchor=north east}, compat=1.3,
  xmin=0.01, xmax=1,ymin=0,ymax=1.75,
  xlabel= {$\varepsilon$},
  ylabel= {}]
\addplot[solid, thick, densely dotted,color=black,mark=none,mark size=1pt] table [x index=0, y index=3]
{friction_noise-lan_dt1.0E-05_N1.0E+04_data.txt};
\addlegendentry{$\hat \sigma_{K^\varepsilon,N}^2/\varepsilon$}
\addplot[solid, thick, color=black,mark=none,mark size=1pt] table [x index=0, y index=4]
{friction_noise-lan_dt1.0E-05_N1.0E+04_data.txt};
\addlegendentry{$\hat \sigma_{K^\varepsilon,N}^2/\varepsilon^2$}
\end{axis}
\node at (5.5cm,3.5cm) {\textcolor{black}{(d)}};
\end{tikzpicture}
\caption{Example \ref{example81a} (friction problem).  The target is to estimate $I^\eps = \EE [ ( X_T^\eps )^2 ]$ for $T=1$ where $X^\eps$ satisfies \eqref{eq:mode1} with  $\varphi(x) \triangleq c_{\rm f} |x|$ with $c_{\rm f}>0$. The expectation of the control variate ${I^0=}\mathbb{E} [(X^0_T)^2]$, where $X^0$ satisfies \eqref{eq:msde_limit1}, is obtained by solving the  partial differential inclusion \eqref{eq:pdi_friction}.}
\label{fig:3}
\end{figure}
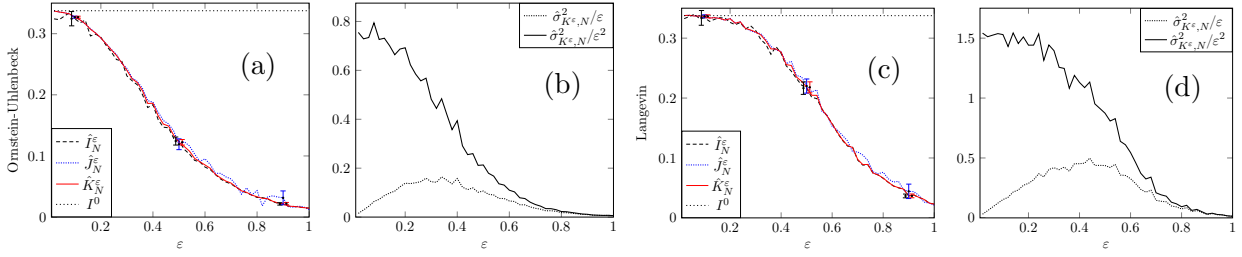

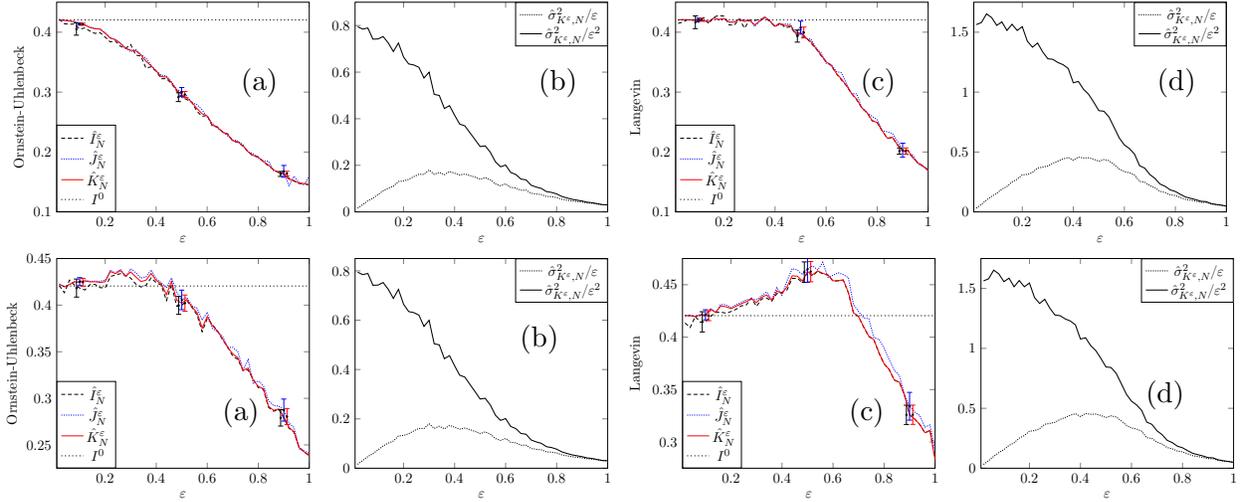
\begin{figure}[h!]
\centering
\begin{tikzpicture}[scale=0.49]
\begin{axis}[legend style={at={(0,0)},anchor=south west}, compat=1.3,
  xmin=0.01, xmax=1,ymin=0.1,ymax=0.45,
  xlabel= {$\varepsilon$},
  ylabel= {Ornstein-Uhlenbeck}]
\addplot[densely dashed,color=black,mark=none,mark size=1pt] table [x index=0, y index=5]{epp_noise-ou_dt1.0E-05_N1.0E+04_data.txt};
\addlegendentry{$\hat{I}_N^\varepsilon$}
\addplot[densely dotted,thick,color=blue,mark=none,mark size=1pt] table [x index=0, y index=6]{epp_noise-ou_dt1.0E-05_N1.0E+04_data.txt};
\addlegendentry{$\hat{J}_N^\varepsilon$}
\addplot[solid,color=red,mark=none,mark size=1pt] table [x index=0, y index=7]{epp_noise-ou_dt1.0E-05_N1.0E+04_data.txt};
\addlegendentry{$\hat{K}_N^\varepsilon$}
\addplot[thick,dotted,color=black,mark=none,mark size=4pt] table [x index=0, y index=8]{epp_noise-ou_dt1.0E-05_N1.0E+04_data.txt};
\addlegendentry{$I^0$}
\addplot[only marks,
  black, mark options={black,scale=0.25}, 
  error bars/.cd, 
    y fixed,
    y dir=both, 
    y explicit
] table [x=x, y=y,y error=error, col sep=comma] {
    x,  y,       error
0.0875,  0.405295,  0.0102782   
0.4875,  0.291524,  0.00733287  
0.8875, 0.163363,  0.00410669
};
\addplot[only marks,
  blue, mark options={black,scale=0.25}, 
  error bars/.cd, 
    y fixed,
    y dir=both, 
    y explicit
] table [x=x, y=y,y error=error, col sep=comma] {
    x,  y,       error
0.1,  0.413555,  0.00171893  
0.5,  0.30032,  0.00735843  
0.9,  0.167737,  0.00986134  
};
\addplot[only marks,
  red, mark options={black,scale=0.25}, 
  error bars/.cd, 
    y fixed,
    y dir=both, 
    y explicit
] table [x=x, y=y,y error=error, col sep=comma] {
    x,  y,       error
0.1125,  0.413302,  0.00168887  
0.5125,  0.295907,  0.00518969 
0.9125,  0.163992,  0.00381172 
};
\end{axis}
\node at (5.5cm,3.5cm) {\textcolor{black}{(a)}};
\end{tikzpicture}
\begin{tikzpicture}[scale=0.49]
\begin{axis}[legend style={at={(1,1)},anchor=north east}, compat=1.3,
  xmin=0.01, xmax=1,ymin=0,ymax=0.9,
  xlabel= {$\varepsilon$},
  ylabel= {}]
\addplot[solid, thick, densely dotted,color=black,mark=none,mark size=1pt] table [x index=0, y index=3]{epp_noise-ou_dt1.0E-05_N1.0E+04_data.txt};
\addlegendentry{$\hat \sigma_{K^\varepsilon,N}^2/\varepsilon$}
\addplot[solid, thick, color=black,mark=none,mark size=1pt] table [x index=0, y index=4]{epp_noise-ou_dt1.0E-05_N1.0E+04_data.txt};
\addlegendentry{$\hat \sigma_{K^\varepsilon,N}^2/\varepsilon^2$}
\end{axis}
\node at (5.5cm,3.5cm) {\textcolor{black}{(b)}};
\end{tikzpicture}
\centering
\begin{tikzpicture}[scale=0.49]
\begin{axis}[legend style={at={(0,0)},anchor=south west}, compat=1.3,
  xmin=0.01, xmax=1,ymin=0.1,ymax=0.45,
  xlabel= {$\varepsilon$},
  ylabel= {Langevin}]
\addplot[densely dashed,color=black,mark=none,mark size=1pt] table [x index=0, y index=5]{epp_noise-lan_dt1.0E-05_N1.0E+04_data.txt};
\addlegendentry{$\hat{I}_N^\varepsilon$}
\addplot[densely dotted,thick,color=blue,mark=none,mark size=1pt] table [x index=0, y index=6]{epp_noise-lan_dt1.0E-05_N1.0E+04_data.txt};
\addlegendentry{$\hat{J}_N^\varepsilon$}
\addplot[solid,color=red,mark=none,mark size=1pt] table [x index=0, y index=7]{epp_noise-lan_dt1.0E-05_N1.0E+04_data.txt};
\addlegendentry{$\hat{K}_N^\varepsilon$}
\addplot[thick,dotted,color=black,mark=none,mark size=4pt] table [x index=0, y index=8]{epp_noise-lan_dt1.0E-05_N1.0E+04_data.txt};
\addlegendentry{$I^0$}
\addplot[only marks,
  black, mark options={black,scale=0.25}, 
  error bars/.cd, 
    y fixed,
    y dir=both, 
    y explicit
] table [x=x, y=y,y error=error, col sep=comma] {
    x,  y,       error
0.0875,  0.416569,  0.0106306   
0.4875,  0.393438,  0.0100313  
0.8875, 0.201527,  0.00526664
};
\addplot[only marks,
  blue, mark options={black,scale=0.25}, 
  error bars/.cd, 
    y fixed,
    y dir=both, 
    y explicit
] table [x=x, y=y,y error=error, col sep=comma] {
    x,  y,       error
0.1,  0.420973,  0.00246508  
0.5,  0.407591,  0.0107992  
0.9,  0.202944,  0.0117355  
};
\addplot[only marks,
  red, mark options={black,scale=0.25}, 
  error bars/.cd, 
    y fixed,
    y dir=both, 
    y explicit
] table [x=x, y=y,y error=error, col sep=comma] {
    x,  y,       error
0.1125,  0.420845,  0.00244568  
0.5125,  0.39946,  0.00900979 
0.9125,  0.201566,  0.00525804 
};
\end{axis}
\node at (5.5cm,3.5cm) {\textcolor{black}{(c)}};
\end{tikzpicture}
\begin{tikzpicture}[scale=0.49]
\begin{axis}[legend style={at={(1,1)},anchor=north east}, compat=1.3,
  xmin=0.01, xmax=1,ymin=0,ymax=1.75,
  xlabel= {$\varepsilon$},
  ylabel= {}]
\addplot[solid, thick, densely dotted,color=black,mark=none,mark size=1pt] table [x index=0, y index=3]{epp_noise-lan_dt1.0E-05_N1.0E+04_data.txt};
\addlegendentry{$\hat \sigma_{K^\varepsilon,N}^2/\varepsilon$}
\addplot[solid, thick, color=black,mark=none,mark size=1pt] table [x index=0, y index=4]{epp_noise-lan_dt1.0E-05_N1.0E+04_data.txt};
\addlegendentry{$\hat \sigma_{K^\varepsilon,N}^2/\varepsilon^2$}
\end{axis}
\node at (5.5cm,3.5cm) {\textcolor{black}{(d)}};
\end{tikzpicture}
\begin{tikzpicture}[scale=0.49]
\begin{axis}[legend style={at={(0,0)},anchor=south west}, compat=1.3,
  xmin=0.01, xmax=1,ymin=0.225,ymax=0.45,
  xlabel= {$\varepsilon$},
  ylabel= {Ornstein-Uhlenbeck}]
\addplot[densely dashed,color=black,mark=none,mark size=1pt] table [x index=0, y index=5]{proba_epp_noise-ou_dt1.0E-05_N1.0E+04_data.txt};
\addlegendentry{$\hat{I}_N^\varepsilon$}
\addplot[densely dotted,thick,color=blue,mark=none,mark size=1pt] table [x index=0, y index=6]{proba_epp_noise-ou_dt1.0E-05_N1.0E+04_data.txt};
\addlegendentry{$\hat{J}_N^\varepsilon$}
\addplot[solid,color=red,mark=none,mark size=1pt] table [x index=0, y index=7]
{proba_epp_noise-ou_dt1.0E-05_N1.0E+04_data.txt};
\addlegendentry{$\hat{K}_N^\varepsilon$}
\addplot[thick,dotted,color=black,mark=none,mark size=4pt] table [x index=0, y index=8]
{proba_epp_noise-ou_dt1.0E-05_N1.0E+04_data.txt};
\addlegendentry{$I^0$}
\addplot[only marks,
  black, mark options={black,scale=0.25}, 
  error bars/.cd, 
    y fixed,
    y dir=both, 
    y explicit
] table [x=x, y=y,y error=error, col sep=comma] {
    x,  y,       error
0.0875,  0.4183,  0.00966829   
0.4875,  0.3997,  0.0096008  
0.8875, 0.279,  0.00879075
};
\addplot[only marks,
  blue, mark options={black,scale=0.25}, 
  error bars/.cd, 
    y fixed,
    y dir=both, 
    y explicit
] table [x=x, y=y,y error=error, col sep=comma] {
    x,  y,       error
0.1,  0.425194,  0.00431547  
0.5,  0.405594,  0.0104796  
0.9,  0.287394,  0.0118217  
};
\addplot[only marks,
  red, mark options={black,scale=0.25}, 
  error bars/.cd, 
    y fixed,
    y dir=both, 
    y explicit
] table [x=x, y=y,y error=error, col sep=comma] {
    x,  y,       error
0.1125,  0.424516,  0.00420984  
0.5125,  0.402089,  0.00876667
0.9125,  0.28038,  0.00864653 
};
\end{axis}
\node at (5cm,1.5cm) {\textcolor{black}{(a)}};
\end{tikzpicture}
\begin{tikzpicture}[scale=0.49]
\begin{axis}[legend style={at={(1,1)},anchor=north east}, compat=1.3,
  xmin=0.01, xmax=1,ymin=0,ymax=0.85,
  xlabel= {$\varepsilon$},
  ylabel= {}]
\addplot[solid, thick, densely dotted,color=black,mark=none,mark size=1pt] table [x index=0, y index=3]{epp_noise-ou_dt1.0E-05_N1.0E+04_data.txt};
\addlegendentry{$\hat \sigma_{K^\varepsilon,N}^2/\varepsilon$}
\addplot[solid, thick, color=black,mark=none,mark size=1pt] table [x index=0, y index=4]{epp_noise-ou_dt1.0E-05_N1.0E+04_data.txt};
\addlegendentry{$\hat \sigma_{K^\varepsilon,N}^2/\varepsilon^2$}
\end{axis}
\node at (5cm,3.5cm) {\textcolor{black}{(b)}};
\end{tikzpicture}
\begin{tikzpicture}[scale=0.49]
\begin{axis}[legend style={at={(0,0)},anchor=south west}, compat=1.3,
  xmin=0.01, xmax=1,ymin=0.275,ymax=0.475,
  xlabel= {$\varepsilon$},
  ylabel= {Langevin}]
\addplot[densely dashed,color=black,mark=none,mark size=1pt] table [x index=0, y index=5]{proba_epp_noise-lan_dt1.0E-05_N1.0E+04_data.txt};
\addlegendentry{$\hat{I}_N^\varepsilon$}
\addplot[densely dotted,thick,color=blue,mark=none,mark size=1pt] table [x index=0, y index=6]{proba_epp_noise-lan_dt1.0E-05_N1.0E+04_data.txt};
\addlegendentry{$\hat{J}_N^\varepsilon$}
\addplot[solid,color=red,mark=none,mark size=1pt] table [x index=0, y index=7]
{proba_epp_noise-lan_dt1.0E-05_N1.0E+04_data.txt};
\addlegendentry{$\hat{K}_N^\varepsilon$}
\addplot[thick,dotted,color=black,mark=none,mark size=4pt] table [x index=0, y index=8]
{proba_epp_noise-lan_dt1.0E-05_N1.0E+04_data.txt};
\addlegendentry{$I^0$}
\addplot[only marks,
  black, mark options={black,scale=0.25}, 
  error bars/.cd, 
    y fixed,
    y dir=both, 
    y explicit
] table [x=x, y=y,y error=error, col sep=comma] {
    x,  y,       error
0.0875,  0.4146,  0.009656   
0.4875,  0.462,  0.00977166   
0.8875,  0.3259,  0.00918672
};
\addplot[only marks,
  blue, mark options={black,scale=0.25}, 
  error bars/.cd, 
    y fixed,
    y dir=both, 
    y explicit
] table [x=x, y=y,y error=error, col sep=comma] {
    x,  y,       error
0.1,  0.421294,  0.00486457  
0.5,  0.464894,  0.0125659  
0.9,  0.334094,  0.0132258  
};
\addplot[only marks,
  red, mark options={black,scale=0.25}, 
  error bars/.cd, 
    y fixed,
    y dir=both, 
    y explicit
] table [x=x, y=y,y error=error, col sep=comma] {
    x,  y,       error
0.1125,  0.420446,  0.00470833 
0.5125,  0.46248,  0.00963905 
0.9125,  0.326013,  0.00918576 
};
\end{axis}
\node at (5cm,1.5cm) {\textcolor{black}{(c)}};
\end{tikzpicture}
\begin{tikzpicture}[scale=0.49]
\begin{axis}[legend style={at={(1,1)},anchor=north east}, compat=1.3,
  xmin=0.01, xmax=1,ymin=0,ymax=1.75,
  xlabel= {$\varepsilon$},
  ylabel= {}]
\addplot[solid, thick, densely dotted,color=black,mark=none,mark size=1pt] table [x index=0, y index=3]{epp_noise-lan_dt1.0E-05_N1.0E+04_data.txt};
\addlegendentry{$\hat \sigma_{K^\varepsilon,N}^2/\varepsilon$}
\addplot[solid, thick, color=black,mark=none,mark size=1pt] table [x index=0, y index=4]{epp_noise-lan_dt1.0E-05_N1.0E+04_data.txt};
\addlegendentry{$\hat \sigma_{K^\varepsilon,N}^2/\varepsilon^2$}
\end{axis}
\node at (5cm,2.0cm) {\textcolor{black}{(d)}};
\end{tikzpicture}
\caption{Example \ref{example81b} (elasto-plastic problem). 
In the top row the target is to estimate $I^\eps = \EE[  ( X_T^\eps )^2 + ( Z_T^\eps   )^2 ]$ with $T=1$ where $(X^\eps,Z^\eps)$ satisfies \eqref{eq:mode2} with  $\varphi(x) \triangleq 0$ and $\psi(x) \triangleq 0$ if $|x| \leq c_{\rm ep}$ and $\infty$ otherwise. Here $c_{\rm ep}= 0.25$. The expectation of the control variate ${I^0=}\mathbb{E} [  (X^0_T )^2 + ( Z^0_T  )^2]$, where $(X^0,Z^0)$ satisfies \eqref{eq:msde_limit}, is obtained by using the PDE method of \cite{msw19}.
In the bottom row the target is to estimate $I^\eps = \PP(|Z^\eps_T|=c_{\rm ep})$.}
\label{fig:4}
\end{figure}

\textcolor{black}{Similarly to what was presented in Section 4,} we use $N=10^4$ samples with a time step of $\delta t = 10^{-5}$. 
In Figures \ref{fig:3} and \ref{fig:4}, we report the numerical results for the friction and elasto-plastic problems,
which are of the form (\ref{eq:mode1}) and (\ref{eq:mode2}), respectively.
In Figures \ref{fig:5} and \ref{fig:6}, we report the numerical results for the obstacle problem and for the reflection of the integral of a colored noise.
\textcolor{black}{The numerical results include errors bars on the estimators for each value of $\varepsilon \in \{0.1,0.5,0.9 \}$ in the $\hat{I}^\eps_N,\hat{J}^\eps_N,\hat{K}^\eps_N$ order.}\\ 

The theoretical predictions provided by (\ref{eq:prop5:1}) and (\ref{eq:prop5:2}) are based on the condition that $f$ has bounded derivatives.
The assumption that $f$ is smooth is important but the hypothesis on the boundedness of the derivatives can certainly be relaxed.
The numerical results shown in Figures \ref{fig:3} and \ref{fig:4} are in good agreement with the theoretical predictions: the 
asymptotic variances $\sigma^2_{J^\eps}$  and $\sigma^2_{K^\eps}$ behave as $O(\eps^2)$. The only cases where the behavior is $O(\eps)$, and not $O(\eps^2)$, are when the quantity of interest is of the form $\EE [ f(\bX^\eps_T)]$ or $\EE [ f(\bX^\eps_T,\bZ^\eps_T)]$ with a function $f$ that is not smooth, which is not surprising.
In Figure \ref{fig:4}, we also observe that $\sigma^2_{J^\eps}$  and $\sigma^2_{K^\eps}$ behave as $O(\eps^2)$ which is better than the behaviour $O(\eps)$ expected  from (\ref{eq:prop5:2}) (which is an upper bound).
In Figures \ref{fig:5} and \ref{fig:6}, the numerical results concern two problems which do not fall within the scope of our theoretical predictions. 
The first one (Figure \ref{fig:5}) is the impact problem that cannot be formulated in the form a differential inclusion of the form \eqref{eq:mode1} or \eqref{eq:mode2}.
The function $f$ is smooth but the behavior of $\sigma^2_{J^\eps}$  and $\sigma^2_{K^\eps}$  is not of order $O(\eps^2)$, only of order $O(\eps)$. 
The second one (Figure \ref{fig:6}) is the reflection of an integrated colored noise that can be formulated with a differential inclusion which is similar to \eqref{eq:mode1} but the multivalued drift does not satisfy the condition \eqref{eq:condfriction}. However, $\sigma^2_{J^\eps}$  and $\sigma^2_{K^\eps}$  behave as $O(\eps^2)$.
To summarize, the numerical simulations indicate that the $O(\eps^2)$ behavior of the asymptotic variance is observed in the cases predicted by the theory and also slightly beyond. The smoothness of the function $f$ that appears in the quantity of interest is, however, an important condition
to ensure the $O(\eps^2)$-behavior, otherwise one only observes a  $O(\eps)$-behavior.


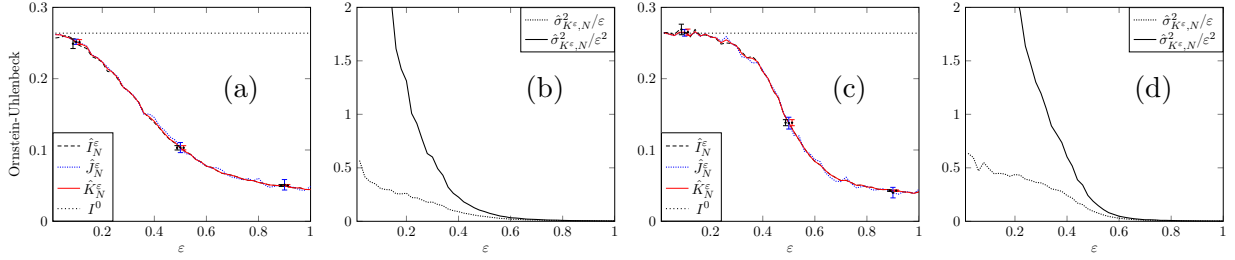
\begin{figure}[h!]
\centering
%
%
\begin{tikzpicture}[scale=0.5]
\begin{axis}[legend style={at={(0,0)},anchor=south west}, compat=1.3,
  xmin=0.01, xmax=1,ymin=0.0,ymax=0.3,
  xlabel= {$\varepsilon$},
  ylabel= {Ornstein-Uhlenbeck}]
\addplot[densely dashed,color=black,mark=none,mark size=1pt] table [x index=0, y index=5]{obstacle_noise-ou_dt1.0E-05_N1.0E+04_data.txt};
\addlegendentry{$\hat{I}_N^\varepsilon$}
\addplot[densely dotted,thick,color=blue,mark=none,mark size=1pt] table [x index=0, y index=6]{obstacle_noise-ou_dt1.0E-05_N1.0E+04_data.txt};
\addlegendentry{$\hat{J}_N^\varepsilon$}
\addplot[solid,color=red,mark=none,mark size=1pt] table [x index=0, y index=7]{obstacle_noise-ou_dt1.0E-05_N1.0E+04_data.txt};
\addlegendentry{$\hat{K}_N^\varepsilon$}
\addplot[thick,dotted,color=black,mark=none,mark size=4pt] table [x index=0, y index=8]{obstacle_noise-ou_dt1.0E-05_N1.0E+04_data.txt};
\addlegendentry{$I^0$}
\addplot[only marks,
  black, mark options={black,scale=0.25}, 
  error bars/.cd, 
    y fixed,
    y dir=both, 
    y explicit
] table [x=x, y=y,y error=error, col sep=comma] {
    x,  y,       error
0.0875,  0.249059,  0.00681488  
0.4875,  0.103066,  0.00306103
0.8875,  0.0501863,  0.00131598
};
\addplot[only marks,
  blue, mark options={black,scale=0.25}, 
  error bars/.cd, 
    y fixed,
    y dir=both, 
    y explicit
] table [x=x, y=y,y error=error, col sep=comma] {
    x,  y,       error
0.1,  0.251731,  0.00365398 
0.5,  0.103389,  0.00707759 
0.9,  0.0512204,  0.00734227
};
\addplot[only marks,
  red, mark options={black,scale=0.25}, 
  error bars/.cd, 
    y fixed,
    y dir=both, 
    y explicit
] table [x=x, y=y,y error=error, col sep=comma] {
    x,  y,       error
0.1125,  0.251268,  0.00344009  
0.5125, 0.103102,  0.00295411 
0.9125,  0.0501874,  0.00131596 
};
\end{axis}
\node at (5cm,3.5cm) {\textcolor{black}{(a)}};
\end{tikzpicture}
%
%
\begin{tikzpicture}[scale=0.5]
\begin{axis}[legend style={at={(1,1)},anchor=north east}, compat=1.3,
  xmin=0.01, xmax=1,ymin=0,ymax=2,
  xlabel= {$\varepsilon$},
  ylabel= {}]
\addplot[solid, thick, densely dotted,color=black,mark=none,mark size=1pt] table [x index=0, y index=3]{obstacle_noise-ou_dt1.0E-05_N1.0E+04_data.txt};
\addlegendentry{$\hat \sigma_{K^\varepsilon,N}^2/\varepsilon$}
\addplot[solid, thick, color=black,mark=none,mark size=1pt] table [x index=0, y index=4]{obstacle_noise-ou_dt1.0E-05_N1.0E+04_data.txt};
\addlegendentry{$\hat \sigma_{K^\varepsilon,N}^2/\varepsilon^2$}
\end{axis}
\node at (5cm,3.5cm) {\textcolor{black}{(b)}};
\end{tikzpicture}
\centering
%
%
\begin{tikzpicture}[scale=0.5]
\begin{axis}[legend style={at={(0,0)},anchor=south west}, compat=1.3,
  xmin=0.01, xmax=1,ymin=0.0,ymax=0.3,
  xlabel= {$\varepsilon$},
  ylabel= {}]
\addplot[densely dashed,color=black,mark=none,mark size=1pt] table [x index=0, y index=5]{obstacle_noise-lan_dt1.0E-05_N1.0E+04_data.txt};
\addlegendentry{$\hat{I}_N^\varepsilon$}
\addplot[densely dotted,thick,color=blue,mark=none,mark size=1pt] table [x index=0, y index=6]{obstacle_noise-lan_dt1.0E-05_N1.0E+04_data.txt};
\addlegendentry{$\hat{J}_N^\varepsilon$}
\addplot[solid,color=red,mark=none,mark size=1pt] table [x index=0, y index=7]{obstacle_noise-lan_dt1.0E-05_N1.0E+04_data.txt};
\addlegendentry{$\hat{K}_N^\varepsilon$}
\addplot[thick,dotted,color=black,mark=none,mark size=4pt] table [x index=0, y index=8]{obstacle_noise-lan_dt1.0E-05_N1.0E+04_data.txt};
\addlegendentry{$I^0$}
\addplot[only marks,
  black, mark options={black,scale=0.25}, 
  error bars/.cd, 
    y fixed,
    y dir=both, 
    y explicit
] table [x=x, y=y,y error=error, col sep=comma] {
    x,  y,       error
0.0875,  0.268996,  0.00739578    
0.4875,  0.138292,  0.00425503 
0.8875,  0.0429911,  0.00108723
};
\addplot[only marks,
  blue, mark options={black,scale=0.25}, 
  error bars/.cd, 
    y fixed,
    y dir=both, 
    y explicit
] table [x=x, y=y,y error=error, col sep=comma] {
    x,  y,       error
0.1,  0.264236,  0.00460209 
0.5,  0.137642,  0.00828022 
0.9,  0.0402229,  0.00739251
};
\addplot[only marks,
  red, mark options={black,scale=0.25}, 
  error bars/.cd, 
    y fixed,
    y dir=both, 
    y explicit
] table [x=x, y=y,y error=error, col sep=comma] {
    x,  y,       error
0.1125,  0.265197,  0.00434649 
0.5125,  0.138276,  0.00425099 
0.9125,  0.0429867,  0.00108716 
};
\end{axis}
\node at (5cm,3.5cm) {\textcolor{black}{(c)}};
\end{tikzpicture}
%
%
\begin{tikzpicture}[scale=0.5]
\begin{axis}[legend style={at={(1,1)},anchor=north east}, compat=1.3,
  xmin=0.01, xmax=1,ymin=0,ymax=2,
  xlabel= {$\varepsilon$},
  ylabel= {}]
\addplot[solid, thick, densely dotted,color=black,mark=none,mark size=1pt] table [x index=0, y index=3]{obstacle_noise-lan_dt1.0E-05_N1.0E+04_data.txt};
\addlegendentry{$\hat \sigma_{K^\varepsilon,N}^2/\varepsilon$}
\addplot[solid, thick, color=black,mark=none,mark size=1pt] table [x index=0, y index=4]{obstacle_noise-lan_dt1.0E-05_N1.0E+04_data.txt};
\addlegendentry{$\hat \sigma_{K^\varepsilon,N}^2/\varepsilon^2$}
\end{axis}
\node at (5cm,3.5cm) {\textcolor{black}{(d)}};
\end{tikzpicture}
\caption{Example \ref{example82a} (impact problem). 
The target is to estimate $I^\eps = \EE [ ( X_{2,T}^\eps )^2 ]$ with $T=1$ where $\bX^\eps$ satisfies the impact problem with a colored noise forcing. Here $P_O = 0.25$. The expectation of the control variate ${I^0=} \EE [ (X^0_{2,T} )^2 ]$, where $\bU$ satisfies the impact problem with a colored noise forcing, is obtained by using the PDE method of \cite{msw19}.}
\label{fig:5}
\end{figure}
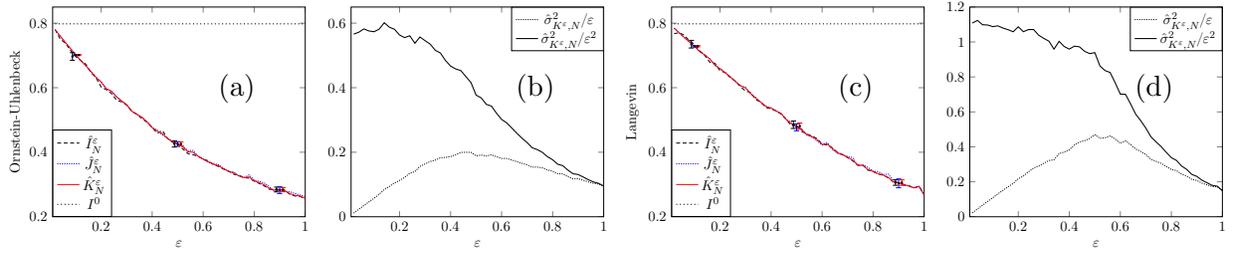
\begin{figure}[h!]
\centering
%
%
\begin{tikzpicture}[scale=0.49]
\begin{axis}[legend style={at={(0,0)},anchor=south west}, compat=1.3,
  xmin=0.01, xmax=1,ymin=0.2,ymax=0.85,
  xlabel= {$\varepsilon$},
  ylabel= {Ornstein-Uhlenbeck}]
\addplot[densely dashed,color=black,mark=none,mark size=1pt] table [x index=0, y index=5]{rbm_noise-ou_dt1.0E-05_N1.0E+04_data.txt};
\addlegendentry{$\hat{I}_N^\varepsilon$}
\addplot[densely dotted,thick,color=blue,mark=none,mark size=1pt] table [x index=0, y index=6]{rbm_noise-ou_dt1.0E-05_N1.0E+04_data.txt};
\addlegendentry{$\hat{J}_N^\varepsilon$}
\addplot[solid,color=red,mark=none,mark size=1pt] table [x index=0, y index=7]{rbm_noise-ou_dt1.0E-05_N1.0E+04_data.txt};
\addlegendentry{$\hat{K}_N^\varepsilon$}
\addplot[thick,dotted,color=black,mark=none,mark size=4pt] table [x index=0, y index=8]{rbm_noise-ou_dt1.0E-05_N1.0E+04_data.txt};
\addlegendentry{$I^0$}
\addplot[only marks,
  black, mark options={black,scale=0.25}, 
  error bars/.cd, 
    y fixed,
    y dir=both, 
    y explicit
] table [x=x, y=y,y error=error, col sep=comma] {
    x,  y,       error
0.0875,  0.696784,  0.0116888   
0.4875,  0.425856,  0.00982864 
0.8875,  0.284518,  0.00745026
};
\addplot[only marks,
  blue, mark options={black,scale=0.25}, 
  error bars/.cd, 
    y fixed,
    y dir=both, 
    y explicit
] table [x=x, y=y,y error=error, col sep=comma] {
    x,  y,       error
0.1,  0.701216,  0.00149796 
0.5,  0.42535,  0.00727561  
0.9,  0.282435,  0.0103088
};
\addplot[only marks,
  red, mark options={black,scale=0.25}, 
  error bars/.cd, 
    y fixed,
    y dir=both, 
    y explicit
] table [x=x, y=y,y error=error, col sep=comma] {
    x,  y,       error
0.1125,  0.701141,  0.0014843  
0.5125,  0.425524,  0.00601745 
0.9125,  0.283848, 0.00638318 
};
\end{axis}
\node at (5cm,3.5cm) {\textcolor{black}{(a)}};
\end{tikzpicture}
%
%
\begin{tikzpicture}[scale=0.49]
\begin{axis}[legend style={at={(1,1)},anchor=north east}, compat=1.3,
  xmin=0.01, xmax=1,ymin=0,ymax=0.65,
  xlabel= {$\varepsilon$},
  ylabel= {}]
\addplot[solid, thick, densely dotted,color=black,mark=none,mark size=1pt] table [x index=0, y index=3]{rbm_noise-ou_dt1.0E-05_N1.0E+04_data.txt};
\addlegendentry{$\hat \sigma_{K^\varepsilon,N}^2/\varepsilon$}
\addplot[solid, thick, color=black,mark=none,mark size=1pt] table [x index=0, y index=4]{rbm_noise-ou_dt1.0E-05_N1.0E+04_data.txt};
\addlegendentry{$\hat \sigma_{K^\varepsilon,N}^2/\varepsilon^2$}
\end{axis}
\node at (5cm,3.5cm) {\textcolor{black}{(b)}};
\end{tikzpicture}
\centering
%
%
\begin{tikzpicture}[scale=0.49]
\begin{axis}[legend style={at={(0,0)},anchor=south west}, compat=1.3,
  xmin=0.01, xmax=1,ymin=0.2,ymax=0.85,
  xlabel= {$\varepsilon$},
  ylabel= {Langevin}]
\addplot[densely dashed,color=black,mark=none,mark size=1pt] table [x index=0, y index=5]{rbm_noise-lan_dt1.0E-05_N1.0E+04_data.txt};
\addlegendentry{$\hat{I}_N^\varepsilon$}
\addplot[densely dotted,thick,color=blue,mark=none,mark size=1pt] table [x index=0, y index=6]{rbm_noise-lan_dt1.0E-05_N1.0E+04_data.txt};
\addlegendentry{$\hat{J}_N^\varepsilon$}
\addplot[solid,color=red,mark=none,mark size=1pt] table [x index=0, y index=7]{rbm_noise-lan_dt1.0E-05_N1.0E+04_data.txt};
\addlegendentry{$\hat{K}_N^\varepsilon$}
\addplot[thick,dotted,color=black,mark=none,mark size=4pt] table [x index=0, y index=8]{rbm_noise-lan_dt1.0E-05_N1.0E+04_data.txt};
\addlegendentry{$I^0$}
\addplot[only marks,
  black, mark options={black,scale=0.25}, 
  error bars/.cd, 
    y fixed,
    y dir=both, 
    y explicit
] table [x=x, y=y,y error=error, col sep=comma] {
    x,  y,       error
0.0875,  0.734865,  0.011861   
0.4875,  0.485254,  0.0116129 
0.8875,  0.306029,  0.00854819
};
\addplot[only marks,
  blue, mark options={black,scale=0.25}, 
  error bars/.cd, 
    y fixed,
    y dir=both, 
    y explicit
] table [x=x, y=y,y error=error, col sep=comma] {
    x,  y,       error
0.1, 0.728343,  0.00205244  
0.5,  0.47739,  0.0108463 
0.9,  0.303532,  0.013696
};
\addplot[only marks,
  red, mark options={black,scale=0.25}, 
  error bars/.cd, 
    y fixed,
    y dir=both, 
    y explicit
] table [x=x, y=y,y error=error, col sep=comma] {
    x,  y,       error
0.1125,  0.728443,  0.00204443 
0.5125, 0.480846,  0.00949818 
0.9125,  0.305794,  0.0084748
};
\end{axis}
\node at (5cm,3.5cm) {\textcolor{black}{(c)}};
\end{tikzpicture}
%
%
\begin{tikzpicture}[scale=0.49]
\begin{axis}[legend style={at={(1,1)},anchor=north east}, compat=1.3,
  xmin=0.01, xmax=1,ymin=0,ymax=1.2,
  xlabel= {$\varepsilon$},
  ylabel= {}]
\addplot[solid, thick, densely dotted,color=black,mark=none,mark size=1pt] table [x index=0, y index=3]{rbm_noise-lan_dt1.0E-05_N1.0E+04_data.txt};
\addlegendentry{$\hat \sigma_{K^\varepsilon,N}^2/\varepsilon$}
\addplot[solid, thick, color=black,mark=none,mark size=1pt] table [x index=0, y index=4]{rbm_noise-lan_dt1.0E-05_N1.0E+04_data.txt};
\addlegendentry{$\hat \sigma_{K^\varepsilon,N}^2/\varepsilon^2$}
\end{axis}
\node at (5cm,3.5cm) {\textcolor{black}{(d)}};
\end{tikzpicture}

\caption{Example \ref{example82b} (reflection of an integrated colored noise).
The target is to estimate $I^\eps = \mathbb{E} [X^\eps_T ]$ for $T=1$ where $X^\eps$ satisfies $\eqref{eq:ref_icn}$. 
The expectation of the control variate is $I^0=\EE[ X^0_T] = \sqrt{{2}/{\pi}}$.}
\label{fig:6}
\end{figure}

\section{Concluding remark}
\label{sec:conclu}
\textcolor{black}{When the expectation of the limit process $\EE[f(\bX^0_T)]$ cannot be computed by a PDE method but is estimated by a massive Monte Carlo method,
the control variate method with $\hat{J}^\eps_N$ (or $\hat{K}_N^\eps$) shares an important similarity with a two-level Monte Carlo method \cite{R3349310} in the sense that massive but cheap simulations are performed (samples of $\bX^0$ used to estimate $\EE[f(\bX^0_T)]$) together with a few expensive simulations (samples of $(\bX^\eps,\bX^0)$ used to estimate $\EE[f(\bX^\eps_T)-f(\bX^0_T)]$).
More generally, Multilevel Monte Carlo (MLMC) methods rely on random samples taken on different levels of accuracy, when several approximations with different costs and accuracies are available. The overall idea of MLMC methods is to reduce the computational cost of standard Monte Carlo methods by taking most samples with a low accuracy and corresponding low cost, and by taking only few samples with a high accuracy and corresponding high cost \cite{R3349310}. } 

\textcolor{black}{
In this two-level Monte Carlo framework, the total cost of computing the control variate estimator $\hat{J}^\eps_N$ (or $\hat{K}_N^\eps$) is
 $C_J^\eps= N_0 C_0 + N_1 C_1$ where $C_0$ is the cost of computing one realization of $f(\bX_T^0)$, $C_1$ is the cost of computing one realization of $(f(\bX_T^\eps) - f(\bX_T^0))$, $N_0$ is the number of samples of $f(\bX_T^0)$, and $N_1$ is the number of samples of $(f(\bX_T^\eps) - f(\bX_T^0))$. 
The variance of the  control variate estimator $\hat{J}^\eps_N$ is $V_J^\eps=N_0^{-1} V_0 + N_1^{-1} V_1$,  where $V_0$ is the variance of $f(\bX_T^0)$ and $V_1$ is the variance of $(f(\bX_T^\eps) - f(\bX_T^0))$.
For a fixed total budget $C_{\rm tot}$, the variance is minimized when  $(N_1/N_0)^2 = (V_1/V_0) (C_0/C_1)$  and it is then equal to 
$V_J^\eps=(\sqrt{V_0C_0}+\sqrt{V_1C_1})^2/C_{\rm tot}$.
This can be compared to the brute force Monte Carlo method: the cost is $C_I^\eps=N_I C_1$ where $N_I$ is the number of samples of $f(\bX_T^\eps)$ (we neglect the difference of cost between $(f(\bX_T^\eps) - f(\bX_T^0))$ and $f(\bX_T^\eps)$, which is very small because $\bX^\eps$ is more difficult to simulate than $\bX^0$) and the variance is $V_I^\eps = V_0 / N_I$ (by (\ref{eq:expandJepsN})), so that for the total budget $C_{\rm tot}$, we have $V_I^\eps = V_0C_1/C_{\rm tot}$.
If $V_1/V_0=O(\eps^2)$
 (by Proposition~\ref{prop:21}) and $C_0/C_1=O(\eps^2)$ (because 
the time step used to simulate $\bX^\eps$ should be $O(\eps^2)$ smaller than the one used to simulate $\bX^0$),
then we find that $N_1/N_0=O(\eps^2)$ (hence the ``massive Monte Carlo" strategy for the control variate) and the ratio of the 
variance of the control variate estimator over the one of the brute force Monte Carlo estimator is finally $V_J^\eps/V_I^\eps=O(\eps^2)$. The control variate method is very advantageous in this context.
}

\appendix

\section{Proof of Proposition \ref{prop:0}}
\label{app:prooflemmaptf}
Let us first study the driving noise.
The process $\boeta^1$ is a Gaussian, Markov process. It has the form
$$
\boeta^1_t = e^{-{\bf A} t} \boeta^1_0 +  \int_0^t e^{-{\bf A}(t-s)} {\bf K} {\rm d}\bW_s .
$$
Its infinitesimal generator is:
\begin{equation}
\label{def:Q}
Q = \frac{1}{2} \sum_{k,k'=1}^d \sum_{k''=1}^{d'} K_{kk''}K_{k'k''} {\partial_{\eta_k\eta_{k'}}^2} - \sum_{k,k'=1}^d A_{kk'} \eta_{k'}  {\partial_{\eta_k}} .
\end{equation}
The properties of the matrix ${\bf A}$ show that the process $\boeta^1$ is stationary and ergodic; 
its unique invariant probability measure is the Gaussian measure 
with mean zero and variance ${\bf C}$ given by~(\ref{def:C}).

\noindent 
The process $(\boeta^\eps_t)_{t\geq 0}$ has the same distribution as $(\boeta^1_{t/\eps^2})_{t\geq 0}$ because $(\eps^{-1}W_t)_{t\geq 0}$ has the same distribution as $(W_{t/\eps^2})_{t\geq 0}$. Therefore it is a Markov process 
with generator $\eps^{-2} Q$.

\noindent The process $(\bX^\eps,\bU,\boeta^\eps)$ is Markov with generator ${\cal L}^\eps$ given by:
\begin{align}
\nonumber
{\cal L}^\eps =& \frac{1}{\eps^2} Q
+\frac{1}{\eps} \Big[   \sum_{j=1}^n \sum_{i=1}^d \sigma_{ji}(\bx) \eta_i \partial_{x_j}
 +   \sum_{j=1}^n \sum_{i=1}^d\sum_{k=1}^{d'} \Gamma_{jk}(\bx^0) K_{ik} \partial^2_{\eta_i x_j^0} \Big] \\
 &
+ \Big[ \sum_{j=1}^n b_j (\bx)\partial_{x_j} +\frac{1}{2} 
\sum_{i,j=1}^n \sum_{k=1}^{d'} \Gamma_{ik}(\bx^0)\Gamma_{jk}(\bx^0)\partial^2_{x^0_i x^0_j}   +\sum_{j=1}^n \widetilde{b}_j(\bx^0)\partial_{x^0_j}  \Big] ,
\label{def:Leps}
\end{align}
where $Q$ is the generator (\ref{def:Q}).

\begin{lem}
\label{lem:ptf}
{\it For any smooth and bounded test function $\phi:\RR^n \times \RR^n \to \RR$, for any compact subset $K$ of  $\RR^{2n}$,
there exists a test function $\phi^\eps$ such that
\begin{align}
\label{pertf0}
&\sup_{(\bx,\bx^0) \in K } |\phi^\eps(\bx,\bx^0,\boeta) - \phi(\bx,\bx^0)|  \leq C \eps (1+\|\boeta\|^2)
\, , \\
&\sup_{(\bx,\bx^0) \in K } |{{\cal L}}^\eps \phi^\eps(\bx,\bx^0,\boeta) - {{\cal L}} \phi(\bx,\bx^0)| \leq C \eps (1+\|\boeta\|^3)
\, ,
\end{align}
for any $\eps \in (0,1)$,  
where ${\cal L}$  is the generator  defined by
\begin{align}
\nonumber
{\cal L} =& \sum_{j=1}^n   \widetilde{b}_j(\bx) \partial_{x_j}  
 + \sum_{j=1}^n   \widetilde{b}_j(\bx^0) \partial_{x^0_j} + \frac{1}{2}
 \sum_{i,j=1}^n [\bGamma(\bx)\bGamma(\bx)^T]_{ij} \partial^2_{x_ix_j}\\
&
+\frac{1}{2} \sum_{i,j=1}^n [\bGamma(\bx^0)\bGamma(\bx^0)^T]_{ij}\partial^2_{x^0_ix^0_j} +
\sum_{i,j=1}^n  [\bGamma(\bx)\bGamma(\bx^0)^T]_{ij}\partial^2_{x_i x^0_j}    .
\label{def:calL}
\end{align}
}
\end{lem}

\noindent
{\it Proof.}
Let $\phi(\bx,\bx^0)$ be a smooth and bounded test function.
We look for a perturbed test function $\phi^\eps$ of the form
\begin{equation}
\label{pert-test-funct}
\phi^\eps(\bx,\bx^0,\boeta) = \phi(\bx,\bx^0) + \eps \phi_{1}(\bx,\bx^0,\boeta) + \eps^2 \phi_{2}(\bx,\bx^0,\boeta) \, .
\end{equation}
Applying ${{\cal L}}^\eps$ (given by (\ref{def:Leps})) to  this $\phi^\eps$ we get
\begin{align}
\nonumber
{{\cal L}}^\eps \phi^\eps =&
\frac{1}{\eps} \Big[ Q \phi_1(\bx,\bx^0,\boeta) +  \sum_{j=1}^n \sum_{i=1}^d \sigma_{ji}(\bx) \eta_i \partial_{x_j} \phi(\bx,\bx^0)  \Big]\\
\nonumber
&+ \Big[ Q \phi_2(\bx,\bx^0,\boeta) + \sum_{j=1}^n \sum_{i=1}^d \sigma_{ji}(\bx) \eta_i \partial_{x_j} \phi_1(\bx,\bx^0,\boeta) +
 \sum_{j=1}^n \sum_{i=1}^d \sum_{k=1}^{d'} \Gamma_{jk}(\bx^0) K_{ik} \partial^2_{\eta_i x^0_j}  \phi_1(\bx,\bx^0,\boeta)\Big]
\\
\nonumber
&+
\Big[ \sum_{j=1}^d b_j(\bx)\partial_{x_j} \phi(\bx,\bx^0)+
\sum_{j=1}^d \widetilde{b}_j(\bx^0)\partial_{x^0_j} \phi(\bx,\bx^0) +\frac{1}{2} 
\sum_{i,j=1}^n \sum_{k=1}^{d'} \Gamma_{ik}(\bx^0)\Gamma_{jk}(\bx^0)\partial^2_{x^0_i x^0_j} 
\phi(\bx,\bx^0)
\Big] \\
 &  + O(\eps) \, .
\label{l-eps-phi-eps}
\end{align}
The term $O(\eps)$ depends on the first-order derivatives of $\phi_1$ and $\phi_2$ with respect to $\bx$ and on the 
first- and second-order derivatives of $\phi_1$ and $\phi_2$ with respect to $\bx^0$.

We define the first corrector $\phi_1$ to cancel the $\eps^{-1}$ term in
(\ref{l-eps-phi-eps}). 
This gives a Poisson equation for $\phi_1$
as a function of $\boeta$ with $(\bx,\bx^0)\in \RR^{2n}$ a frozen parameter.
The Poisson equation 
\begin{equation}
\label{eq:poisson1}
Q\bg=- \boeta
\end{equation}
can be solved by Fredholm alternative
because the process $\boeta^1$ has mean zero  (with respect to its invariant probability measure) \cite[Chapter 6]{book}.
We can write a solution in the form:
$$
\bg(\boeta) = \int_0^\infty \EE [ \boeta^1_s|\boeta^1_0=\boeta] {\rm d}s ,
$$
which is here linear in $\boeta$:
\begin{equation}
\label{eq:expressg}
\bg(\boeta) = {\bf A}^{-1} \boeta.
\end{equation}
Therefore we set
\begin{equation}
\label{eq:phi1}
\phi_1(\bx,\bx^0 ,\boeta) =  \sum_{j=1}^n \sum_{i=1}^d \sigma_{ji}(\bx) g_i(\boeta) \partial_{x_j} \phi(\bx,\bx^0) .
\end{equation}
We cannot define the second corrector $\phi_2$  so as to 
cancel the order-one terms in (\ref{l-eps-phi-eps}) because
that would require solving a Poisson equation with a right-hand side that is not centered. 
To center this term we subtract its mean
relative to the invariant distribution of $\boeta^1$.
This gives the Poisson equation
\begin{align*}
&
Q \phi_{2}(\bx,\bx^0,\boeta) +  
\sum_{j=1}^n \sum_{i=1}^d \sigma_{ji}(\bx) \eta_i \partial_{x_j} \phi_1(\bx,\bx^0,\boeta) +
 \sum_{j=1}^n \sum_{i=1}^d\sum_{k=1}^{d'} \Gamma_{jk}(\bx^0) K_{ik} \partial^2_{\eta_i x^0_j}  \phi_1(\bx,\bx^0,\boeta) 
 \\
&- \EE \Big[ \sum_{j=1}^n \sum_{i=1}^d \sigma_{ji}(\bx) \eta^1_{i,0} \partial_{x_j} \phi_1(\bx,\bx^0,\boeta^1_0) +
 \sum_{j=1}^n \sum_{i=1}^d\sum_{k=1}^{d'} \Gamma_{jk}(\bx^0) K_{ik} \partial^2_{\eta_i x^0_j}  \phi_1(\bx,\bx^0,\boeta^1_0)\Big] =0\, ,
\end{align*}
where the expectation $\EE$ is taken over $\boeta^1_0 $ with respect to its invariant probability measure.
Since $\phi_1$ is linear in $\boeta$, the third term of the left-hand side is independent of $\boeta$ and is equal to its expectation,
so the Poisson equation can be reduced to
\begin{align*}
Q \phi_{2}(\bx,\bx^0,\boeta) +  
\sum_{j=1}^n \sum_{i=1}^d \sigma_{ji}(\bx) \eta_i \partial_{x_j} \phi_1(\bx,\bx^0,\boeta)
- \EE \Big[ \sum_{j=1}^n \sum_{i=1}^d \sigma_{ji}(\bx) \eta^1_{i,0} \partial_{x_j} \phi_1(\bx,\bx^0,\boeta^1_0) \Big] =0\, .
\end{align*}
This equation has a solution $\phi_2$
that is a smooth function in $(\bx,\bx^0)$ and that is a quadratic form in $\boeta$.
Note that $\phi_1(\bx,\bx^0,\boeta)$ and $\phi_2(\bx,\bx^0,\boeta)$ depend only on $\bsigma(\bx)$ and its first-order derivatives, and not on $\bb$.
By assuming that $\bsigma$ belongs to ${\cal C}^2$ with bounded derivatives, we get the control of the $O(\eps)$ term in (\ref{l-eps-phi-eps}).
It   follows that
\begin{align*}
{{\cal L}}^\eps \phi^\eps =& \EE \Big[ \sum_{j=1}^n \sum_{i=1}^d \sigma_{ji}(\bx) \eta^1_{i,0} \partial_{x_j} \phi_1(\bx,\bx^0,\boeta^1_0) +
 \sum_{j=1}^n  \sum_{i=1}^d\sum_{k=1}^{d'}  \Gamma_{jk}(\bx^0) K_{ik} \partial^2_{\eta_i x^0_j}  \phi_1(\bx,\bx^0,\boeta^1_0)\Big] \\
\nonumber
&+
\Big[ \sum_{j=1}^n b_j(\bx)\partial_{x_j} \phi(\bx,\bx^0)+
\sum_{j=1}^n \widetilde{b}_j(\bx^0)\partial_{x^0_j} \phi(\bx,\bx^0) +\frac{1}{2} 
\sum_{i,j=1}^n \sum_{k=1}^{d'} \Gamma_{ik}(\bx^0)\Gamma_{jk}(\bx^0)\partial^2_{x^0_i x^0_j} 
\phi(\bx,\bx^0)
\Big]\\
 & + O(\eps)\, .
\end{align*} 
Using \eqref{eq:phi1}, the expectation takes the form
\begin{align*}
& \EE \Big[ \sum_{j=1}^n \sum_{i=1}^d \sigma_{ji}(\bx) \eta^1_{i,0} \partial_{x_j} \phi_1(\bx,\bx^0,\boeta^1_0) +
 \sum_{j=1}^n  \sum_{i=1}^d\sum_{k=1}^{d'}  \Gamma_{jk}(\bx^0) K_{ik} \partial^2_{\eta_i x^0_j}  \phi_1(\bx,\bx^0,\boeta^1_0)\Big] \\
& =
\sum_{j,j'=1}^n \sum_{i,i'=1}^d  \EE\big[ \eta^1_{i,0}g_{i'}(\boeta^1_0)\big] \sigma_{ji}(\bx) \partial_{x_j} 
\big( \sigma_{j'i'}(\bx) \partial_{x_{j'}} \phi(\bx,\bx^0) \big)
 \\
&\quad +
 \sum_{j,j'=1}^n \sum_{i,i'=1}^d  \sum_{k=1}^{d'} 
  \EE \big[ \partial_{\eta_i}  g_{i'}(\boeta^1_0) \big]\Gamma_{jk}(\bx^0) K_{ik} \partial_{x^0_j} \big(  \sigma_{j'i'}(\bx)  \partial_{x_{j'}} \phi(\bx,\bx^0)
\big) .
\end{align*}
From the explicit form (\ref{eq:expressg}) of $\bg$ we get
$$
 \EE\big[ \eta^1_{i,0}g_{i'}(\boeta^1_0)\big] =   \sum_{k=1}^d ({\bf A}^{-1})_{i'k} C_{ki} =  ({\bf A}^{-1} {\bf C})_{i'i}
\mbox{ and }
 \EE \big[ \partial_{\eta_i}  g_{i'}(\boeta^1_0) \big] =  ({\bf A}^{-1})_{i'i} .
$$
Therefore
\begin{align*}
& \EE \Big[ \sum_{j=1}^n \sum_{i=1}^d \sigma_{ji}(\bx) \eta^1_{i,0} \partial_{x_j} \phi_1(\bx,\bx^0,\boeta^1_0) +
 \sum_{j=1}^n  \sum_{i=1}^d\sum_{k=1}^{d'} \Gamma_{jk}(\bx^0) K_{ik} \partial^2_{\eta_i x^0_j}  \phi_1(\bx,\bx^0,\boeta^1_0)\Big] \\
& = 
\sum_{j,j'=1}^n \sum_{i,i'=1}^d  ({\bf A}^{-1} {\bf C})_{i'i} \sigma_{ji}(\bx) \partial_{x_j} 
\big( \sigma_{j'i'}(\bx) \partial_{x_{j'}} \phi(\bx,\bx^0) \big)
\\
&
\quad +
 \sum_{j,j'=1}^n \sum_{i,i'=1}^d  \sum_{k=1}^{d'} 
\Gamma_{jk}(\bx^0) K_{ik}   \sigma_{j'i'}(\bx) ({\bf A}^{-1})_{i'i} \partial_{x^0_j x_{j'}} \phi(\bx,\bx^0)
 \\
 & =
\sum_{j,j'=1}^n  
\big( \bsigma (\bx) {\bf A}^{-1} {\bf C} \bsigma(\bx)^T  \big)_{j'j} \partial_{x_{j'} x_j} \phi(\bx,\bx^0)  
\\
&\quad  +\sum_{j'=1}^n \Big(  \sum_{j=1}^n ( \partial_{x_j}  \bsigma(\bx) )
{\bf A}^{-1} {\bf C} \bsigma(\bx)^T )_{j'j}
\Big)
 \partial_{x_{j'}} \phi(\bx,\bx^0)  
\\
&
\quad +
 \sum_{j,j'=1}^n
  \big( \bsigma (\bx) {\bf A}^{-1} {\bf K} \bGamma(\bx^0)^T)_{j'j} \partial_{x^0_j x_{j'}} \phi(\bx,\bx^0)
\\
& = \frac{1}{2}
\sum_{j,j'=1}^n  
\big( \bsigma (\bx) ( {\bf A}^{-1} {\bf C}+{\bf C} {{\bf A}^T}^{-1}) \bsigma(\bx)^T  \big)_{jj'} \partial_{x_{j'} x_j} \phi(\bx,\bx^0)  
\\
&\quad  +\sum_{j'=1}^n\big( \widetilde{b}_{j'}(\bx) - b_{j'}(\bx)\big) \partial_{x_{j'}} \phi(\bx,\bx^0)  
+
 \sum_{j,j'=1}^n
  \big( \bGamma (\bx)   \bGamma(\bx^0)^T)_{j'j'} \partial_{x^0_j x_{j'}} \phi(\bx,\bx^0)
 .
\end{align*}
We have from ${\rm d}\boeta^1 = {\bf K} {\rm d}\bW_t -{\bf A} \boeta^1 {\rm d}t$ and It\^o's formula:
$$
{\rm d} ( \eta^1_i \eta^1_j ) = \sum_{k=1}^{d'} ( K_{ik}\eta^1_j  + K_{jk} \eta^1_i ){\rm d}W_{kt} +  \sum_{k=1}^{d'}  K_{ik} K_{jk} {\rm d}t
 - \sum_{k=1}^d  (A_{ik}\eta^1_k \eta^1_j +
A_{jk} \eta^1_k \eta^1_i) {\rm d}t .
$$
Taking the expectation (under the invariant probability measure) gives the identity:
$$
{\bf K}{\bf K}^T  - 
{\bf C} {\bf A}^T -{\bf A} {\bf C}  ={\bf 0} .
$$
By left-multiplying by ${\bf A}^{-1}$ and by right-multiplying by $ {{\bf A}^T}^{-1} $ we find
$$
{\bf A}^{-1} {\bf K}{\bf K}^T {{\bf A}^T}^{-1} - 
{\bf A}^{-1} {\bf C} -{\bf C} {{\bf A}^T}^{-1}  ={\bf 0}  ,
$$
which gives
$$
 \bsigma (\bx) ( {\bf A}^{-1} {\bf C}+{\bf C} {{\bf A}^T}^{-1}) \bsigma(\bx)^T  =  \bGamma (\bx)   \bGamma(\bx)^T ,
 $$
 and  we obtain the desired result:
$$
{\cal L}^\eps \phi^\eps = {\cal L}\phi +O(\eps).
$$
\qed

We can then 
prove Proposition \ref{prop:0} as follows.

\noindent
{\it Proof of Proposition \ref{prop:0}.}
By the perturbed test function method  \cite[Section 6.3]{book},
Lemma \ref{lem:ptf} establishes that the continuous process $(\bX^\eps,\bU)$ converges in distribution 
to the Markov process with infinitesimal generator ${\cal L}$ defined by (\ref{def:calL}).
The infinitesimal generator ${\cal L}$ can be associated to a diffusion process $(\tilde{\bX},\tilde{\bU})$ that is solution of 
the coupled SDEs:
\begin{align*}
& {\rm d} \tilde{\bX} = \widetilde{\bb}(\tilde{\bX}) {\rm d}t+   \bGamma(\tilde{\bX}){\rm d}\tilde{\bW}_t ,\\
& {\rm d} \tilde{\bX}^0 = \widetilde{\bb}(\tilde{\bX}^0) {\rm d}t+  \bGamma(\tilde{\bX}^0){\rm d}\tilde{\bW}_t ,
\end{align*} 
where $\tilde{\bW}$ is a $d'$-dimensional Brownian motion.
This shows that, if $\tilde{\bX}_0=\tilde{\bX}^0_0$ almost surely, then $\tilde{\bX}_t-\tilde{\bX}^0_t=0$ for all $t$ almost surely.
Therefore, if $\bX^\eps_0=\bX^0_0$, then the continuous process $(\bX^\eps-\bX^0)$ converges in distribution to $0$, which implies convergence in probability.
\qed


\section{Proof of Lemma \ref{prop:3}}
\label{app:proofprop3}

%
%
We define
\begin{equation}
\textcolor{black}{
\phi(\bx,\bx^0) = {g}(\bx^0)\big( {f}(\bx)-{f}(\bx^0)  \big)
\mbox{ or } 
}
\phi(\bx,\bx^0) = ({f}(\bx)-{f}(\bx^0))^2 .
\end{equation}
Lemma \ref{lem:ptf} applied to $\phi$ gives an estimate for (\ref{eq:estimeffeps2}) of order $\eps$, but the particular form of $\phi$ 
makes it possible to get $\eps^2$, as we show in the following.
We prove Lemma \ref{prop:3} in four steps.

{\it Step 1. 
There exist smooth functions $\phi_{1i},\phi_{20},\phi_{2ij},\Lambda_{1i},\Lambda_{1ijk}$  with bounded derivatives
such that
\begin{align}
\phi^\eps(\bx,\bx^0,\boeta) =& \phi(\bx,\bx^0)+\eps \phi_1(\bx,\bx^0,\boeta) +\eps^2  \phi_2(\bx,\bx^0,\boeta) ,\\
\label{eq:expandphi1}
\phi_1(\bx,\bx^0,\boeta)=& \sum_{i=1}^d \phi_{1i}(\bx,\bx^0)\eta_i , \\
 \phi_2(\bx,\bx^0,\boeta) =&  \phi_{20}(\bx,\bx^0) +\sum_{i,j=1}^d \phi_{2ij}(\bx,\bx^0)\eta_i \eta_j  ,\\
{\cal L}^\eps \phi^\eps(\bx,\bx^0,\boeta) =& \eps \Lambda_1(\bx,\bx^0,\boeta)  +O(\eps^2) ,\\
\label{eq:expandlambda1}
 \Lambda_1(\bx,\bx^0,\boeta)=& \sum_{i=1}^d \Lambda_{1i}(\bx,\bx^0) \eta_i +\sum_{i,j,k=1}^d \Lambda_{1ijk}(\bx,\bx^0)\eta_i \eta_j\eta_k .
\end{align}
}
{\it  Proof.}
We apply the perturbed test function method as described in the proof of Lemma \ref{lem:ptf} and we get the result
by keeping track of the $\boeta$-dependence of the perturbed functions $\phi_1$ and $\phi_2$.
\qed

{\it Step 2. For $s\leq t$, the conditional distribution of $\boeta^\eps_t$ given ${\cal F}_s =\sigma(\bW_u, u \leq s)$ is
\begin{equation}
\label{eq:conddistetaeps}
{\cal N}\Big( \exp\big(- \frac{{\bf A}(t-s)}{\eps^2}\big) \boeta^\eps_s, \int_0^{(t-s)/\eps^2} e^{-{\bf A} u} {\bf K}{\bf K}^T e^{-{\bf A}^T u} {\rm d} u\Big) .
\end{equation}
There exists $\lambda, C >0$ such that 
\begin{align}
\label{eq:step21}
\big| \EE \big[ \eta_{i,t}^\eps |{\cal F}_s  \big]  \big| &\leq e^{-\lambda (t-s)/\eps^2 } \|\boeta^\eps_s\|  ,\\
\big| \EE \big[ \eta_{i,t}^\eps \eta_{j,t}^\eps \eta_{k,t}^\eps |{\cal F}_s  \big]  \big| &\leq C e^{-\lambda (t-s)/\eps^2 }  \|\boeta^\eps_s\| (1+ \|\boeta^\eps_s\|^2) .
\label{eq:step22}
\end{align}
}
{\it Proof.}
We can integrate (\ref{eq:sde2}) from $s$ to $t$:
$$
\boeta^\eps_t =  \exp\big(- \frac{{\bf A}(t-s)}{\eps^2}\big) \boeta^\eps_s
+
\int_s^t  \exp\big(- \frac{{\bf A}(t-u)}{\eps^2}\big) \frac{{\bf K}}{\eps}  {\rm d} \bW_u ,
$$
which gives (\ref{eq:conddistetaeps}) and
$$
\EE\big[ \boeta^\eps_t | {\cal F}_s \big] =  \exp\big(- \frac{{\bf A}(t-s)}{\eps^2}\big) \boeta^\eps_s.
$$
Eq.~(\ref{eq:step21}) is a straightforward consequence. 
Eq.~(\ref{eq:step22}) follows from  (\ref{eq:conddistetaeps}) and Isserlis theorem for multivariate normal random vectors.
\qed

{\it Step 3. If $(\bx,\bx^0)\mapsto \psi(\bx,\bx^0)$ is a smooth function  with bounded derivatives, then there exists $C>0$ such that,
for all $i,j,k=1,\ldots,d$ and $t\in [\eps,T]$:
\begin{align}
\label{eq:step31}
\big| \EE \big[ \eta_{i,t}^\eps \psi (\bX^\eps_t,\bU_t) \big]  \big|  &\leq C \eps  ,\\
\big| \EE \big[ \eta_{i,t}^\eps \eta_{j,t}^\eps \eta_{k,t}^\eps  \psi (\bX^\eps_t,\bU_t) \big]  \big| &\leq C \eps .
\label{eq:step32}
\end{align}
}
{\it Proof.}
We have for any $t \in[ \eps,T]$ and $\delta <\eps$:
\begin{align*}
 \psi (\bX^\eps_t,\bU_t)=& \psi (\bX^\eps_{t-\delta},\bU_{t-\delta}) 
+
\int_{t-\delta}^\delta 
\frac{1}{\eps}   \psi_1(\bX^\eps_s,\bU_s,\boeta^\eps_s) {\rm d} s
\\
&+
\int_{t-\delta}^\delta \psi_2(\bX^\eps_s,\bU_s) {\rm d} s + 
\sum_{j=1}^d \int_{t-\delta}^\delta\psi_{3j}(\bX^\eps_s,\bU_s) {\rm d} W_{js} ,
\end{align*}
with
\begin{align*}
\psi_1(\bx,\bx^0,\boeta) =&  \bsigma(\bx) \boeta  \cdot \nabla_\bx \psi(\bx,\bx^0) , \\
\psi_2(\bx,\bx^0)=& 
[ {\itbf b}(\bx)\cdot \nabla_\bx\psi +\tilde{\itbf b}(\bx^0) \cdot \nabla_{\bx^0} \psi + \nabla_{\bx^0} \psi^T \bGamma(\bx^0)^T  \bGamma(\bx^0)  \nabla_{\bx^0} \psi  ](\bx, \bx^0) , \\
\bpsi_{3j}(\bx,\bx^0)=& \sum_{i=1}^n \partial_{x^0_i} \psi (\bx,\bx^0) \Gamma_{ij}(\bx^0) .
\end{align*}
We have, by (\ref{eq:step21}),
\begin{align*}
\big|
\EE\big[ \eta_{i,t}^\eps \psi (\bX^\eps_{t-\delta},\bU_{t-\delta})   \big]
\big|
&=
\big|
\EE\big[ \EE[ \eta_{i,t}^\eps|{\cal F}_{t-\delta}] \psi (\bX^\eps_{t-\delta},\bU_{t-\delta})  \big]
\big|
\leq
C \exp( - \lambda \delta/\eps^2) .
\end{align*}
Similarly, for any $s\in [t-\delta,t]$
\begin{align*}
\big|
\EE\big[ \eta_{i,t}^\eps \psi_1 (\bX^\eps_s,\bU_s,\boeta^\eps_s)  \big]
\big|
&=
\big|
\EE\big[ \EE[ \eta_{i,t}^\eps|{\cal F}_{s}] \psi_1 (\bX^\eps_s,\bU_s,\boeta^\eps_s)  \big]
\big|
\leq
C \exp( - \lambda (t-s)/\eps^2) ,\\
\big|
\EE\big[ \eta_{i,t}^\eps \psi_2 (\bX^\eps_s,\bU_s) \big]
\big|
&\leq 
C \exp( - \lambda (t-s)/\eps^2) ,
\end{align*}
and  for $j=1,\ldots,{d'}$ and for any positive $q$, 
\begin{align*}
&\big|
\EE\big[ \eta_{i,t}^\eps \int_{t-\delta}^t \psi_{3j} (\bX^\eps_s,\bU_s) {\rm d} W_{js} \big]
\big| \\
&
\leq
\sum_{k=0}^{q-1}
\big|  \EE\big[ \EE[ \eta_{i,t}^\eps |{\cal F}_{t-k\delta/q}] \int_{t-(k+1)\delta/q}^{t-k\delta/q} \psi_{3j} (\bX^\eps_s,\bU_s)  {\rm d} W_{js} \big]
\big|\\
&
\leq
\sum_{k=0}^{q-1}
\big|  \EE\big[ \EE[ \eta_{i,t}^\eps |{\cal F}_{t-k\delta/q}]^2 \big]^{1/2}
 \Big[\int_{t-(k+1)\delta/q}^{t-k\delta/q}\EE [ \psi_{3j} (\bX^\eps_s,\bU_s)^2]  {\rm d} s 
\Big]^{1/2}\\
&\leq C \sum_{k=0}^{q-1}\exp( - \lambda k \delta /(q\eps^2)) \sqrt{\delta/q} .
\end{align*}
Consequently
\begin{align*}
&\big|
\EE\big[ \eta_{i,t}^\eps \psi (\bX^\eps_t,\bU_t)  \big]
\big|
\\
&\leq 
C \exp\big( - \frac{\lambda \delta}{\eps^2}\big) + \frac{C}{\eps} \int_{t-\delta}^t \exp \big( - \frac{\lambda  (t-s)}{\eps^2}\big){\rm d} s
+
C \frac{\sqrt{\delta}}{\sqrt{q}} \sum_{k=0}^{q-1}\exp\big( - \frac{\lambda k \delta }{q\eps^2}\big)\\
&\leq 
C \exp\big( - \frac{\lambda \delta}{\eps^2}\big)  + \frac{C\eps}{\lambda}
+
C \frac{\sqrt{\delta}}{\sqrt{q}} \sum_{k=0}^{q-1}\exp\big( - \frac{\lambda k \delta }{q\eps^2}\big)  .
\end{align*}
By taking $\delta=\eps^2|\ln \eps| / \lambda$ and $q=[|\ln \eps|]$ we finally get
\begin{align*}
\big|
\EE\big[ \eta_{i,t}^\eps \psi (\bX^\eps_t,\bU_t)  \big]
\big|
\leq 
C' \eps + C'\eps  \sum_{k=0}^\infty e^{k} \leq C'' \eps  ,
\end{align*}
which gives the first desired result. The calculations with the third-order product of coefficients $\boeta$ are similar and use (\ref{eq:step22}).
\qed

{\it Step 4. Proof of Lemma \ref{prop:3}.}\\
For any $t \in [\eps,T]$, 
we have
\begin{align*}
\EE\big[ \phi(\bX^\eps_t,\bU_t)  \big] &= \EE\big[ \phi^\eps (\bX^\eps_t,\bU_t,\boeta^\eps_t) \big] 
-\eps \EE\big[  \phi_1 (\bX^\eps_t,\bU_t,\boeta^\eps_t)  \big] +O(\eps^2) \\
&
=  \EE\big[ \phi^\eps (\bX^\eps_t,\bU_t,\boeta^\eps_t)  \big] + O(\eps^2)  ,
\end{align*}
because (\ref{eq:expandphi1}) and (\ref{eq:step31}) give $\EE[\phi_1 (\bX^\eps_t,\bU_t,\boeta^\eps_t)  ] =O(\eps)$.\\
We have
\begin{align*}
\EE\big[ \phi^\eps (\bx_0,\bx_0,\boeta^\eps_0)  \big]  &= 
\eps \EE\big[ \phi_1(\bx_0,\bx_0,\boeta^\eps_0)\big]+O(\eps^2)  =O(\eps^2)  ,
\end{align*}
because $\EE[\phi_{1j}(\bx_0,\bx_0)\eta_{j,0}^\eps]=\phi_{1j}(\bx_0,\bx_0)\EE[\eta_{j,0}^\eps]=0$ for all $j=1,\ldots,d$.\\
Therefore
\begin{align*}
\EE\big[ \phi(\bX^\eps_t,\bU_t)  \big] &=
\int_0^t \EE\big[{\cal L}^\eps \phi^\eps(\bX^\eps_s,\bU_s,\boeta^\eps_s)\big] {\rm d} s +O(\eps^2)\\
&=\eps \int_0^t \EE\big[\Lambda_1 (\bX^\eps_s,\bU_s,\boeta^\eps_s)\big] {\rm d} s +O(\eps^2) \\
&=O(\eps^2),
\end{align*}
because (\ref{eq:expandlambda1}), (\ref{eq:step31})  and (\ref{eq:step32}) give $\EE[\Lambda_1 (\bX^\eps_s,\bU_s,\boeta^\eps_s)  ] =O(\eps)$ for any $s\in [\eps,t]$ and 
$\EE[\Lambda_1 (\bX^\eps_s,\bU_s,\boeta^\eps_s)  ] =O(1)$ for any $s \in [0,\eps]$.
This completes the proof of Lemma~\ref{prop:3}.
\qed

\section{Proofs of existence and uniqueness of \eqref{eq:mode1} and \eqref{eq:mode2}}
\label{app:A}
\begin{prop}
\label{thmA1}
Fix $T>0, n \in \mathbb{N}^\star$. Suppose that ${\itbf f} \in \mathcal{C}([0,T];\mathbb{R}^n)$, ${\itbf b}$ is Lipschitz, and $\varphi$ is a l.s.c. convex function satisfying \eqref{eq:condfriction}.
Then there exists a unique solution $\bx \in \mathcal{C}([0,T];\mathbb{R}^n)$ to the following differential inclusion
\begin{equation}
\label{multiode}
\bx(0) =\bx_0 \in \mathbb{R}^n ,\quad \quad
\dot \bx(t) + \partial \varphi (\bx(t)) \ni {\itbf b}(\bx(t)) + {\itbf f}(t), \: t > 0.
\end{equation}
 \end{prop}
\begin{proof}
Let $\varphi_p$   be the Moreau-Yosida regularization of $\varphi$.
For each $p \geq 1$, we consider the penalized problem
$$
\bx^p(0) = \bx_0 \in \mathbb{R}^n, 
\quad \dot \bx^p(t) + \nabla \varphi_p (\bx^p(t)) = {\itbf b}(\bx^p(t)) + {\itbf f}(t), \: t > 0. 
$$ 
This is a standard ODE with Lipschitz coefficients, so $\bx^p \in \mathcal{C}([0,T];\mathbb{R}^n)$ is well-defined. Now, we show that $\bx^p$ is a Cauchy sequence in $\mathcal{C}([0,T];\mathbb{R}^n)$.
Fix $p,q\in \mathbb{N}^\star$ and $t \in [0,T]$. We have the following expansion
\begin{align*}
\frac{1}{2} \| \bx^p(t) - \bx^q(t)\|^2 = & \int_0^t \big( \bx^p(s) - \bx^q(s) \big) \cdot \big( {\itbf b}(\bx^p(s)) - {\itbf b}(\bx^q(s)) \big) \textup{d} s\\ 
& - \int_0^t \big( \bx^p(s) - \bx^q(s) \big) \cdot \big( \nabla \varphi_p(\bx^p(s))-\nabla \varphi_q(\bx^q(s)) \big) \textup{d} s
\end{align*}
which, using the properties of ${\itbf b},\varphi_p, \varphi_q$, leads to the following inequality
\begin{equation}
\label{inegalite4Cauchy1}
\frac{1}{2} \| \bx^p(t) - \bx^q(t)\|^2 \leq C \int_0^t \| \bx^p(s) - \bx^q(s)\|^2 \textup{d} s + \left ( \frac{1}{p} + \frac{1}{q} \right ) \int_0^t  \nabla \varphi_p(\bx^p(s)) \cdot \nabla \varphi_q(\bx^q(s))  \textup{d} s.
\end{equation}
Under the assumption $\sup \limits_{p \geq 1} \sup \limits_{\bx \in \RR^n} \| \nabla \varphi_p(\bx) \| < \infty$, we deduce from the inequality above that 
$$
\frac{1}{2} \| \bx^p(t) - \bx^q(t)\|^2 \leq C \int_0^t \| \bx^p(s) - \bx^q(s)\|^2 \textup{d} s + \left ( \frac{1}{p} + \frac{1}{q} \right ) C t.
$$
Thus, we can apply Gronwall's inequality to obtain
\begin{equation}
\label{inegalite4Cauchy2}
\sup_{0 \leq t \leq T} \| \bx^p(t) - \bx^q(t) \|^2 \leq \left ( \frac{1}{p} + \frac{1}{q} \right ) C_T.
\end{equation}
Therefore $\bx^p$ is a Cauchy sequence and there exists a function $\bx \in \mathcal{C}([0,T];\mathbb{R}^n)$ such that $\bx^p \to \bx$, as $p \to \infty$ in $\mathcal{C}([0,T];\mathbb{R}^n)$. 
Next we verify that $\bx$ satisfies the differential inclusion. 
Define $\forall t \in [0,T], \: {\boldsymbol \Delta}^p(t) \triangleq \int_0^t \nabla \varphi_p(\bx^p(s)) \textup{d} s$ and denote ${\boldsymbol \Delta}(t)  \triangleq \lim \limits_{p \to \infty} {\boldsymbol \Delta}^p(t)$.
We then have 
$$
\forall t \in [0,T], \quad \bx(t) + {\boldsymbol \Delta}(t) = \bx_0 + \int_0^t {\itbf b}(\bx(s)) \textup{d} s + \int_0^t {\itbf f}(s) \textup{d} s.
$$
Moreover, since $\sup \limits_{p} \int_0^T \| \nabla \varphi_p(\bx^p(s)) \|^2 \textup{d} s < \infty$, there exists a function ${\boldsymbol \delta} \in L^2(0,T)$ such that 
$$
\forall {\itbf h}  \in L^2(0,T), \quad
\lim \limits_{p \to \infty} 
\int_0^T   \nabla \varphi_p(\bx^p(s)) \cdot {\itbf h}(s)  \textup{d} s 
= 
\int_0^T   {\boldsymbol \delta}(s) \cdot {\itbf h}(s)  \textup{d} s.
$$ 
As a consequence, we must have ${\boldsymbol \Delta}(t) = \int_0^t {\boldsymbol \delta}(s) \textup{d} s$.
Now, to finally check the differential inclusion, we want to show that $\forall {\itbf v} \in \mathcal{C}([0,T];\mathbb{R}^n)$,
$$
\forall 0 \leq t < t+h \leq T,  \quad
\int_t^{t+h}  {\boldsymbol \delta}(s) \cdot \big(  {\itbf v}(s) - \bx(s) \big) + \varphi(\bx(s)) \textup{d} s  \leq \int_t^{t+h} \varphi({\itbf v}(s)) \textup{d} s.
$$
We exploit property 2 as listed in section \ref{sec:inc} 
$$
\forall t \in [0,T], \quad \nabla \varphi_p(\bx^p(t)) \in \partial \varphi (J_p \bx^p(t))
$$
which implies $\forall {\itbf v} \in \mathcal{C}([0,T];\mathbb{R}^n)$,
$$
\forall t \in [0,T],  \quad
 \nabla \varphi_p(\bx^p(t)) \cdot \big( {\itbf v}(t) - J_p \bx^p(t) \big) + \varphi(J_p \bx^p(t)) \leq \varphi({\itbf v}(t)).
$$
We integrate the inequality above on $[t,t+h]$ and obtain $\forall {\itbf v} \in \mathcal{C}([0,T];\mathbb{R}^n)$,
$$
\forall 0 \leq t < t+h \leq T, \quad
\int_t^{t+h}  \nabla \varphi_p(\bx^p(s)) \cdot\big({\itbf v}(s) - J_p \bx^p(s) \big) + \varphi(J_p \bx^p(s)) \textup{d} s \leq \int_t^{t+h} \varphi({\itbf v}(s)) \textup{d} s.
$$
On the one hand,
by using Fatou's lemma ($\varphi\geq 0$),
the convergence $J_p  \bx^p \to \bx$ in $\mathcal{C}([0,T];\mathbb{R}^n)$ as $p \to +\infty$, and
the fact that  $\varphi$ is l.s.c, we get
$$
\forall 0 \leq t < t+h \leq T, \quad 
\liminf \limits_{p \to \infty }
\int_t^{t+h} \varphi(J_p \bx^p(s)) \textup{d} s \geq 
\int_t^{t+h} \liminf \limits_{p \to \infty }
\varphi(J_p \bx^p(s)) \textup{d} s 
\geq 
\int_t^{t+h} \varphi(\bx(s)) \textup{d} s .
$$
On the other hand, we have
\begin{align*}
\int_t^{t+h}  \nabla \varphi_p(\bx^p(s)) \cdot \big(  {\itbf v}(s) - J_p \bx^p(s) \big) \textup{d} s 
=&
\int_t^{t+h} \nabla \varphi_p(\bx^p(s)) \cdot \big( {\itbf v}(s) - \bx(s) \big) \textup{d} s  \\
&+
\int_t^{t+h}  \nabla \varphi_p(\bx^p(s)) \cdot \big( \bx(s)- J_p \bx^p(s) \big)  \textup{d} s 
.
\end{align*}
By using $\nabla \varphi_p(\bx^p(s)) \to {\boldsymbol \delta}$ in $L^2(0,T)$ weak, the first term of the right-hand side goes to zero.
By using $\|\nabla \varphi_p(\bx^p(s))\| $ is uniformly bounded in $L^2(0,T)$ and $\|J_p  \bx^p -\bx \|\to 0$ in $L^2(0,T)$,
the second term of the right-hand side goes to zero.
Therefore, we have
$$
\liminf \limits_{p \to \infty } \int_t^{t+h}  \nabla \varphi_p(\bx^p(s)) \cdot \big( {\itbf v}(s) - J_p \bx^p(s) \big) + \varphi(J_p \bx^p(s)) \textup{d} s 
\geq
\int_t^{t+h} {\boldsymbol \delta}(s) \cdot \big( {\itbf v}(s) - \bx(s) \big) + \varphi(\bx(s)) \textup{d} s ,
$$
which proves that $\bx$ is solution of \eqref{multiode}.\\
We now show uniqueness. Assume that $\bx$ and $\by$ satisfy \eqref{multiode} with $\bx(0) = \bx_0$ and $\by(0) = \by_0$.
Then 
$$
\bx(t) + {\boldsymbol \Delta}_x(t) = \bx_0 + \int_0^t {\itbf b}(\bx(s)) \textup{d} s+ \int_0^t {\itbf f}(s) \textup{d} s
$$
and
$$
\by(t) + {\boldsymbol \Delta}_y(t) = \by_0 + \int_0^t {\itbf b}(\by(s)) \textup{d} s+ \int_0^t {\itbf f}(s) \textup{d} s.
$$
With $\int_0^t {\boldsymbol \delta}_x(s) \textup{d} s =  {\boldsymbol \Delta}_x(t)$ and $\int_0^t {\boldsymbol \delta}_y(s) \textup{d} s =  {\boldsymbol \Delta}_y(t)$, we have the following inequalities
$$
\forall 0 \leq t < t+h \leq T, \quad \int_t^{t+h} 
 {\boldsymbol \delta}_x(s) \cdot \big( \by(s) - \bx(s) \big) + \varphi(\bx(s)) \textup{d} s \leq \int_t^{t+h} \varphi(\by(s)) \textup{d} s
$$ 
and
$$
\forall 0 \leq t < t+h \leq T, \quad \int_t^{t+h} 
 {\boldsymbol \delta}_y(s) \cdot \big( \bx(s) - \by(s) \big) + \varphi(\by(s)) \textup{d} s  \leq \int_t^{t+h} \varphi(\bx(s)) \textup{d} s,
$$ 
which give
$$
\forall 0 \leq t < t+h \leq T, \quad \int_t^{t+h} 
\big(  {\boldsymbol \delta}_x(s) - {\boldsymbol \delta}_y(s) \big) \cdot \big( \by(s) - \bx(s) \big) \textup{d} s  \leq 0.
$$ 
That means 
$$
\forall 0 \leq t < t+h \leq T, \quad \int_t^{t+h} \big( {\itbf b}(\bx(s))-{\itbf b}(\by(s)) - ( \dot \bx (s) - \dot \by (s) ) \big)\cdot \big( \by(s) - \bx(s) \big)
  \textup{d} s  \leq 0 ,
$$ 
which turns into
$$
\forall 0 \leq t < t+h \leq T, \quad  \| \bx(t+h) - \by(t+h) \|^2  \leq  \| \bx(t) - \by(t) \|^2 + C \int_t^{t+h} \| \bx(s) - \by(s) \|^2  \textup{d} s.
$$ 
Gronwall inequality gives the desired result and it is clear that if $\bx_0 = \by_0$ then $\forall t \in [0,T], \bx(t) = \by(t)$.
\end{proof}

\begin{prop}
\label{thmA2}
Fix $T>0, n,m \in \mathbb{N}^\star$. Suppose that ${\itbf f} \in \mathcal{C}([0,T];\mathbb{R}^n)$, ${\itbf b}^x: \RR^n \to \RR$ and ${\itbf b}^z: \RR^m \to \RR$ are Lipschitz, 
and $\varphi : \RR^n \to \RR, \: \psi : \RR^m \to \RR$ are l.s.c. convex functions, with $\varphi$ satisfying \eqref{eq:condfriction}
and $\psi$ satisfying \eqref{eq:condphipz0}.
Then there exists a unique solution $(\bx,\bz) \in \mathcal{C}([0,T];\mathbb{R}^n \times \mathbb{R}^m )$ to the following differential inclusion
$(\bx(0),\bz(0)) = (\bx_0,\bz_0) \in \mathbb{R}^n \times \mathbb{R}^m$, 
\begin{equation}
\label{multiode2}
\dot \bx(t) + \partial \varphi (\bx(t)) \ni {\itbf b}^x(\bx(t),\bz(t)) + {\itbf f}(t),
\quad \dot \bz(t) + \partial \psi (\bz(t)) \ni {\itbf b}^z(\bx(t),\bz(t)), 
\: t > 0. 
\end{equation}
 \end{prop}
 \begin{proof}
Let $\varphi_p$ and $\psi_p$ be the Moreau-Yosida regularization of $\varphi$ and $\psi$.
We consider the penalized problems $(\bx^p(0),\bz^p(0))=(\bx_0,\bz_0) \in \RR^n \times \RR^m$,
\begin{equation}
\begin{rcases}
& \dot \bx^p(t) + \nabla \varphi_p (\bx^p(t)) = \bb^x(\bx^p(t),\bz^p(t)) + {\itbf f}(t),\\
& \dot \bz^p(t) + \nabla \psi_p (\bz^p(t)) = \bb^z(\bx^p(t),\bz^p(t)),
\end{rcases}
\quad t > 0. 
\end{equation}
It can be shown than $(\bx^p,\bz^p)$ is a Cauchy sequence in $\mathcal{C}([0,T]; \RR^n \times \RR^m)$.
The proof follows similar steps as in the proof of Theorem \ref{thmA1}, 
except that we need a bound of the form 
\begin{equation}
\label{psipbound}
\int_0^t  \nabla \psi_p(\bz^p(s)) \cdot  \nabla \psi_q(\bz^q(s))  \leq C t ,
\end{equation}
where $C$ does not depend on $p,q$, without using $\sup \limits_{p \geq 1} \sup \limits_{\bz \in \RR^m} \| \nabla \psi_p(\bz) \| < \infty$,
which we do not assume.
We proceed first with the following expansion
$$
\psi_p(\bz^p(t)) = \psi_p(\bz_0) + \int_0^t\nabla \psi_p (\bz^p(s)) \cdot \big( \bb^z (\bx^p(s),\bz^p(s)) - \nabla \psi_p(\bz^p(s)) \big) \textup{d} s  ,
$$
which implies 
$$
\psi_p(\bz^p(t)) + \int_0^t \| \nabla \psi_p (\bz^p(s)) \|^2 \textup{d} s = \psi_p(\bz_0) + \int_0^t 
 \nabla \psi_p (\bz^p(s)) \cdot  \bb^z (\bx^p(s),\bz^p(s)) \textup{d} s .
$$
We get
$$
\psi_p(\bz^p(t)) + \frac{1}{2} \int_0^t \| \nabla \psi_p (\bz^p(s)) \|^2 \textup{d} s 
\leq  \psi_p(\bz_0) + \frac{1}{2} \int_0^t \| \bb^z (\bx^p(s),\bz^p(s)) \|^2 \textup{d} s
$$
which 
 implies
\begin{equation}
\label{eq:proofA1:2}
\psi_p(\bz^p(t)) + \frac{1}{2} \int_0^t \| \nabla \psi_p (\bz^p(s)) \|^2 \textup{d} s 
\leq C \left ( 1 + \int_0^t \left \{ \| \bx^p(s) \|^2 + \| \bz^p(s) \|^2 \right \} \textup{d} s \right ).
\end{equation}
Besides, 
\begin{align*}
\| \bx^p(t)\|^2 + \| \bz^p(t)\|^2 = \| \bx^p(0)\|^2 + \| \bz^p(0)\|^2 & + 2 \int_0^t 
 \bx^p(s) \cdot \big( \bb^x(\bx^p(s),\bz^p(s)) - \nabla \varphi_p(\bx^p(s)) + {\itbf f}(s) \big) \textup{d} s\\
& + 2 \int_0^t  \bz^p(s)\cdot\big(  \bb^z(\bx^p(s),\bz^p(s)) - \nabla \psi_p(\bz^p(s)) \big) \textup{d} s.
\end{align*}
Using  $\int_0^T \| {\itbf f}(s) \|^2 \textup{d} s < \infty$, $\forall \bx \in \RR^n,  \, \bx \cdot  \nabla \varphi_p(\bx) \geq 0$,
 and $\forall \bz \in \RR^m, \,  \bz \cdot \nabla \psi_p(\bz) \geq 0$, 
we get
\begin{align*}
\| \bx^p(t)\|^2 + \| \bz^p(t)\|^2 \leq C \left ( 1  +  \int_0^t  \| \bx^p(s) \|^2 + \| \bz^p(s) \|^2 \textup{d} s \right ).
\end{align*}
Thus, Gronwall inequality yields 
$$
\sup \limits_{s \leq t} \left [  \| \bx^p(s) \|^2 + \| \bz^p(s) \|^2 \right ] \leq C \exp (C t) .
$$
Substituting into (\ref{eq:proofA1:2}) and using the fact that $\psi_p\geq0$ gives \eqref{psipbound}. Next we can use the same arguments as 
in the proof of Theorem \ref{thmA1} to show that $\bx$ satisfies the differential inclusion involving $\varphi$.
We finally discuss the case of $\bz$. 
Proceeding as in the proof of Theorem \ref{thmA1}, we have $\forall {\itbf v} \in \mathcal{C}([0,T];\mathbb{R}^n)$,
$$
\forall 0 \leq t < t+h \leq T, \: 
\int_t^{t+h}   \nabla \psi_p(\bz^p(s)) \cdot \big( {\itbf v}(s) - J_p \bz^p(s) \big) + 
\psi(J_p \bz^p(s)) \textup{d} s \leq \int_t^{t+h} \psi({\itbf v}(s)) \textup{d} s.
$$
We do not assume $\sup \limits_{p \geq 1} \sup \limits_{\bz \in \RR^m} \| \nabla \psi_p(\bz) \| < \infty$,
so we cannot claim that  $J_p  \bz^p \to \bz$ in $\mathcal{C}([0,T];\mathbb{R}^m)$, in contrast with the proof
of Theorem \ref{thmA1}.
However, 
since $\nabla \psi_p (\bz^p(s))  =  p \left ( \bz^p(s) - J_p^z \bz^p(s) \right )$ and $\sup \limits_p \int_0^T \| \nabla \psi_p (\bz^p(s)) \|^2 \textup{d} s < \infty$,
we get 
$$
\lim \limits_{p \to \infty} \int_0^T \| \bz^p(s) - J_p^z \bz^p(s) \|^2 \textup{d} s = 0.
$$
We can extract a subsequence $p_k$ to get 
$$
\bz^{p_k} - J_{p_k}^z \bz^{p_k} \to 0, \: a.e \: \mbox{in} \: (0,T).
$$
 Therefore 
 $$
 J_{p_k}^z \bz^{p_k} \to \bz, \: a.e \: \mbox{in} \: (0,T).
 $$
 Observe that by the l.s.c. property of $\psi$, we have  
 $$
 \psi(\bz(t)) \leq \liminf\limits_{k\rightarrow \infty} \psi (J_{p_k}^z \bz^{p_k}(t) ), \: a.e \: \mbox{in} \: (0,T).
 $$
 Thus using Fatou's lemma, we obtain
$$
\int_t^{t+h} \psi(\bz(s)) \textup{d} s 
\leq 
\int_t^{t+h}  \liminf\limits_{k\rightarrow \infty} \psi (J_{p_k}^z \bz^{p_k}(s)) \textup{d} s
\leq \liminf\limits_{k\rightarrow \infty} \int_t^{t+h}   \psi (J_{p_k}^z \bz^{p_k}(s)) \textup{d} s . 
$$ 
Also, we have 
\begin{align*}
\int_t^{t+h}  \nabla \psi_{p_k} (\bz^{p_k}(s)) \cdot \big( {\itbf v}(s) - J_{p_k}^z \bz^{p_k}(s) \big) \textup{d} s
= & \int_t^{t+h} \nabla \psi_{p_k} (\bz^{p_k}(s)) \cdot \big( {\itbf v}(s) - \bz(s) \big) \textup{d} s\\
& + \int_t^{t+h}  \nabla \psi_{p_k} (\bz^{p_k}(s)) \cdot \big( \bz(s) - J_{p_k}^z \bz^{p_k}(s) \big) \textup{d} s.
\end{align*}
In the right-hand side, as $k \uparrow \infty$, 
the first term goes to $\int_t^{t+h} {\boldsymbol \delta}_z(s) \cdot \big({\itbf v}(s) - \bz(s) \big) \textup{d} s$ 
because $ \nabla \psi_{p_k} (\bz^{p_k})$ weakly converges to ${\boldsymbol \delta}_z$ and the second term goes to $0$ 
because 
$\| \nabla \psi_{p_k} (\bz^{p_k}(s)) \|$ is uniformly bounded in $L^2(0,T)$ with respect to $k$ 
and 
$\bz - J_{p_k}^z(\bz^{p_k})$ goes to $0$ in   $L^2(0,T;\RR^m)$. 
Therefore
$$
\forall 0 \leq t < t+h \leq T,  \quad
\int_t^{t+h}  {\boldsymbol \delta}_z(s) \cdot \big( {\itbf v}(s) - \bz(s) \big) + \psi(\bz(s)) \textup{d} s  \leq 
\int_t^{t+h} \psi({\itbf v}(s)) \textup{d} s.
$$
The proof is complete.
\end{proof}

\section{Existence and uniqueness for \eqref{eq:msde_limit1}}
\label{app:D}
Let us first take a look at the case where we remove the multivalued operator $\partial \varphi$ from the drift in \eqref{eq:msde_limit1}. The problem becomes the same as \eqref{eq:sde3} where $\bsigma$ is constant and in particular it does not involve a stochastic integral. Thus, as pointed out in page 294 of \cite{MR1121940}, the proof of existence and uniqueness of a solution (still in \cite{MR1121940}, Theorem 2.9 page 289) can be simplified in a way that makes no use of probabilistic tools. We consider the Wiener space 
$\left ( \Omega \triangleq C \left ( [0,T];\RR^d \right ),  
\mathcal{F} \triangleq \mathcal{B} \left ( C \left ( [0,T];\RR^d \right ) \right ), \PP \right )$
here 
$\Omega$ is the space of $\RR^d$-valued continuous functions on $[0,T]$ endowed with the norm 
$\forall \boldsymbol{\omega} \in \Omega, \: \| \boldsymbol{\omega} \| \triangleq \sup \limits_{0 \leq t \leq T} \| \boldsymbol{\omega}(t) \|$,
$\mathcal{B} \left ( C \left ( [0,T];\RR^d \right ) \right )$ is the Borel $\sigma$-algebra on $\Omega$,
$\PP$ is the Wiener measure;
the mappings 
indexed by $ t \in [0,T]$, $\bW_t( \cdot) \colon \Omega \to \RR^d$, 
$\boldsymbol{\omega} \mapsto \: \bW_t(\boldsymbol{\omega}) \triangleq \boldsymbol{\omega}(t)$,
the sequence of $\sigma$-algebras $\mathcal{F}^t \triangleq \sigma \{ \bW_s, \: 0 \leq s \leq t \}$
and the map $\bX \colon \Omega \to C \left ( [0,T] ; \RR^n \right )$, $\boldsymbol{\omega} \mapsto \bX(\boldsymbol{\omega}) \triangleq \bx$
where $\forall 0 \leq t \leq T, \: \bx(t) = \bx(0) + \int_0^t {\itbf b}(\bx(s)) \textup{d} s + \bGamma \boldsymbol{\omega}(t).$ 
Under $\PP$, $\bW$ is a Wiener process and $\bX(\bW)$ solves \eqref{eq:sde3} where $\bsigma$ is constant. 
In this approach, the key ingredient is the mapping $\bX$. 
For obtaining the existence and uniqueness of the solution to \eqref{eq:msde_limit1}
with the multivalued operator $\partial \varphi$, we discuss below the properties of a similar mapping to $\bX$ which involves the multivalued operator.
This is done via the so-called ``Generalized Skorokhod Problem". The discussion follows \cite{MR3308895} 
from page 245 to page 252.
We use the notation $BV[0,T]$ for the space of functions with bounded variation on $[0,T]$.
\begin{definition}[Generalized Convex Skorokhod Problem]
If a pair of functions $(\bx,\boldsymbol{\Delta})$ satisfies the following conditions
\begin{enumerate}
\item
$\bx,\boldsymbol{\Delta} : [0,T] \to \mathbb{R}^n$ are continuous, $\bx(0) = \bx_0$ and $\boldsymbol{\Delta}(0) = 0$,
\item
$\forall 0 \leq t \leq T, \: \bx(t) \in \overline{\textup{Dom} (\partial \varphi)}, \: \boldsymbol{\Delta} \in  BV([0,T]; \RR^n)$,
\item
$\forall 0 \leq t \leq T, \: \bx(t) + \boldsymbol{\Delta}(t) = x_0 + \int_0^t {\itbf b}(\bx(s)) \textup{d} s + \bGamma \boldsymbol{\omega}(t)$,
\item
$\forall 0 \leq s \leq t \leq T, \: \forall \mathfrak{z} \in \RR^n, \: \int_s^t \big( \bx(r) - \mathfrak{z} \big)\cdot \textup{d} \boldsymbol{\Delta} (r) 
+ \int_s^t \varphi (\bx(r)) \textup{d} r \leq (t-s) \varphi(\mathfrak{z})$,
\end{enumerate}
then we say that $\bx$ solves the  generalized Skorokhod problem with parameters $\partial \varphi, \bx_0, {\itbf b}$ and $\boldsymbol{\omega}$ and we use the notation $\bx = \mathcal{GSP}(\partial \varphi, \bx_0, {\itbf b}, \boldsymbol{\omega})$.
\end{definition}
\noindent Existence and uniqueness of a solution for the Generalized Skorokhod Problem can be found in Theorem 4.17 page 252. This is obtained under the following conditions : $\varphi$ is a l.s.c convex function and $\textup{int} \left ( \textup{Dom}(\varphi) \right ) \neq \emptyset$; ${\itbf b}$ is Lipschitz, $\bx_0 \in \overline{\textup{Dom} (\partial \varphi)}$ and $\boldsymbol{\omega} : [0,T] \mapsto \RR^n$ is continuous with $\boldsymbol{\omega}(0)=0$. The continuity of the mapping $\bX \colon \Omega \to C \left ( [0,T] ; \RR^n \right )$, $\boldsymbol{\omega} \mapsto \bX(\boldsymbol{\omega}) \triangleq \bx$
where $\bx = \mathcal{GSP}(\partial \varphi, \bx_0, {\itbf b}, \boldsymbol{\omega})$ is shown in proposition 4.16 page 247.
We use the notation $\cS_n^0[0,T]$ for the space of progressively measurable continuous stochastic processes (p.m.c.s.p.) from $\Omega \times [0,T]$ to $\RR^n$,
$$
\cS_n^2[0,T] \triangleq \left \{ {\itbf Z} \in \cS_n^0[0,T], \: \EE \Big[ \sup \limits_{0 \leq t \leq T} \| {\itbf Z}(t) \|^2 \Big] < \infty \right \}  .
$$ 
Within the framework of the aforementioned Wiener space, $\bX({\itbf W}) \in \cS_n^0[0,T]$ solves \eqref{eq:msde_limit1}.
Furthermore, it can be shown that $\bX({\itbf W}) \in \cS_n^2[0,T]$.

\section{Proof of Proposition \ref{prop:1}}
\label{app:proofprop1}%
{\it Proof of the first item of Proposition \ref{prop:1}.}
Using 
\eqref{eq:prop1phi}, we get
\begin{align*}
\frac{1}{2} \| \bX_t^{\eps,p} - \bX_t^{\eps,q} \|^2 
& =  \int_0^t (\bX_s^{\eps,p}-\bX_s^{\eps,q} )\cdot ( {\itbf b}(\bX_s^{\eps,p}) - {\itbf b}(\bX_s^{\eps,q} ) ) \textup{d} s\\
& \quad - \int_0^t (\bX_s^{\eps,p}-\bX_s^{\eps,q} )\cdot( \nabla \varphi_p(\bX_s^{\eps,p}) - \nabla \varphi_q(\bX_s^{\eps,q}) ) \textup{d} s\\
& \leq C \int_0^t \| \bX_s^{\eps,p} - \bX_s^{\eps,q} \|^2 \textup{d} s\\
& \quad + \left ( \frac{1}{p} + \frac{1}{q} \right ) \int_0^t  \nabla \varphi_p (\bX_s^{\eps,p}) \cdot  \nabla \varphi_q (\bX_s^{\eps,q})  \textup{d} s.    
\end{align*}
Then under the condition \eqref{eq:condfriction} and from an application of Gronwall inequality, we obtain 
$$
\| \bX_t^{\eps,p} - \bX_t^{\eps,q} \|^2 \leq C_{b,\varphi,t} \left ( \frac{1}{p} + \frac{1}{q} \right ).
$$ 
Here the constant $C_{b,\varphi,t}$ depends on $b, \varphi$ and $t$. This implies the result.
\qed

{\it Proof of the second item of Proposition \ref{prop:1}.}
Using the same arguments as in the proof of the first item, we get that
\begin{equation}
\label{eq:lem1:bis}
\forall p \geq 1, \:  \mathbb{E} \left[ \sup \limits_{ t \leq T} \| \bX^{0,p}_t - \bX^0_t \|^2\right] \leq \frac{C_T}{p} ,
\end{equation}
where $\bX^{0,p}$ is an approximation of $\bX^0$ in the following sense:
$$
 {\rm d}\bX^{0,p} + \nabla \varphi_p (\bX^{0,p}) {\rm d}t =  \bb(\bX^{0,p}) {\rm d}t +  \bGamma  {\rm d} \bW_t . 
$$
Therefore, for any $\delta >0$ and $p$,
\begin{align*}
\PP\left( \sup_{t \leq T} \|\bX^\eps_t - \bX^0_t\|\geq \delta \right)
&\leq
\PP\left( \sup_{t \leq T} \|\bX^\eps_t - \bX^{\eps,p}_t\|\geq \frac{\delta}{3} \right)
+
 \PP\left( \sup_{t \leq T} \|\bX^{\eps,p}_t - \bX^{0,p}_t\|\geq \frac{\delta}{3}  \right) \\
&\quad +
\PP\left( \sup_{t \leq T} \|\bX^{0,p}_t - \bX^0_t\|\geq \frac{\delta}{3}  \right) \\
&\leq  \frac{18 C_T}{p \delta^2} +  \PP\left( \sup_{t \leq T} \|\bX^{\eps,p}_t - \bX^{0,p}_t\|\geq \frac{\delta}{3}  \right),
\end{align*}
by Markov inequality.
From Proposition \ref{prop:0}, we have
$$
\limsup_{\eps \to 0} \PP\left( \sup_{t\in [0,T]} \|\bX^\eps_t - \bX^0_t\|\geq \delta \right)
\leq  \frac{18 C_T}{p \delta^2} ,
$$
which holds for any $p$, hence the desired result.
\qed

\section{Proof of Proposition \ref{prop:2}}
\label{app:proofprop2}%
{\it Proof of the first item of Propostion \ref{prop:2}.}
We first note that  $\psi$ satisfies (\ref{eq:condplastic}).
Using 
\eqref{eq:prop1phi} for $\varphi$ and $\psi$, we get
\begin{align*}
& \frac{1}{2} \| \bX_t^{\eps,p} - \bX_t^{\eps,q} \|^2 + \frac{1}{2} \| \bZ_t^{\eps,p} - \bZ_t^{\eps,q} \|^2  \\
& =  \int_0^t ( \bX_s^{\eps,p}-\bX_s^{\eps,q} ) \cdot ( {\itbf b}^X(\bX_s^{\eps,p},\bZ_s^{\eps,p}) - {\itbf b}^X(\bX_s^{\eps,q},\bZ_s^{\eps,q}) ) \textup{d} s\\
& \quad + \int_0^t (\bZ_s^{\eps,p}-\bZ_s^{\eps,q} ) \cdot ( {\itbf b}^Z(\bX_s^{\eps,p},\bZ_s^{\eps,p}) - {\itbf b}^Z(\bX_s^{\eps,q},\bZ_s^{\eps,q}) ) 
\textup{d} s\\
& \quad - \int_0^t ( \bX_s^{\eps,p}-\bX_s^{\eps,q} ) \cdot ( \nabla \varphi_p(\bX_s^{\eps,p}) - \nabla \varphi_q(\bX_s^{\eps,q}) ) \textup{d} s\\
& \quad - \int_0^t ( \bZ_s^{\eps,p}-\bZ_s^{\eps,q} ) \cdot ( \nabla \psi_p(\bZ_s^{\eps,p}) - \nabla \psi_q(\bZ_s^{\eps,q}) ) \textup{d} s\\
& \leq C \int_0^t \left \{ \| \bX_s^{\eps,p} - \bX_s^{\eps,q} \|^2 +\| \bZ_s^{\eps,p} - \bZ_s^{\eps,q} \|^2 \right \} \textup{d} s\\
& \quad + \left ( \frac{1}{p} + \frac{1}{q} \right ) \left \{  \int_0^t \| \nabla \varphi_p (\bX_s^{\eps,p}) \| \| \nabla \varphi_q (\bX_s^{\eps,q}) \| \textup{d} s
+ \int_0^t \| \nabla \psi_p (\bZ_s^{\eps,p}) \| \| \nabla \psi_q (\bZ_s^{\eps,q}) \| \textup{d} s  \right \}.   
\end{align*}
From \eqref{eq:condfriction}, 
$$
\int_0^t \| \nabla \varphi_p (\bX_s^{\eps,p}) \| \| \nabla \varphi_q (\bX_s^{\eps,q}) \| \textup{d} s \leq C t  .
$$
From now on we focus on $\int_0^t \| \nabla \psi_p (\bZ_s^{\eps,p}) \| \| \nabla \psi_q (\bZ_s^{\eps,q}) \| \textup{d} s$.
Let us proceed with the following expansion 
$$
\psi_p(\bZ_T^{\eps,p}) = \psi_p(\bz_0) + \int_0^T  \nabla \psi_p (\bZ_s^{\eps,p}) \cdot ( {\itbf b}^Z(\bX_s^{\eps,p},\bZ_s^{\eps,p}) - \nabla \psi_p(\bZ_s^{\eps,p}) ) \textup{d} s     ,
$$
which implies 
$$
\psi_p(\bZ_T^{\eps,p}) + \frac{1}{2} \int_0^T \| \nabla \psi_p (\bZ_s^{\eps,p}) \|^2 \textup{d} s
\leq \psi_p(\bz_0) + \frac{1}{2} \int_0^T \| {\itbf b}^Z(\bX_s^{\eps,p} ,\bZ_s^{\eps,p}) \|^2 \textup{d} s     .
$$
Therefore 
\begin{equation}
\label{eq:psibound1}
\psi_p(\bZ_T^{\eps,p}) + \frac{1}{2} \int_0^T \| \nabla \psi_p (\bZ_s^{\eps,p}) \|^2 \textup{d} s
\leq \psi_p(\bz_0) + \frac{C}{2} \int_0^T \{ 1 + \| \bX_s^{\eps,p} \|^2 + \| \bZ_s^{\eps,p} \|^2 \} \textup{d} s.     
\end{equation}
We study the term in the integral of the right-hand side of the inequality (\ref{eq:psibound1}):
\begin{align}
\nonumber
\| \bX_s^{\eps,p} \|^2 + \| \bZ_s^{\eps,p} \|^2 & =  \| \bX_0^{\eps,p} \|^2 + \| \bZ_0^{\eps,p} \|^2 + 2 \int_0^s  \bX_r^{\eps,p} \cdot \Big( {\itbf b}^X(\bX_r^{\eps,p} , \bZ_r^{\eps,p}) - \nabla \phi_p(\bX_r^{\eps,p}) + \bsigma \frac{\boeta_r^\eps}{\eps} \Big) \textup{d} r\\
\nonumber
& \quad + 2 \int_0^s \bZ_r^{\eps,p} \cdot \big( {\itbf b}^Z(\bX_r^{\eps,p} ,\bZ_r^{\eps,p}) - \nabla \psi_p(\bZ_r^{\eps,p}) \big) \textup{d} r\\ 
\nonumber
& \leq  \| \bX_0^{\eps,p} \|^2 + \| \bZ_0^{\eps,p} \|^2 + 2 C \int_0^s \{ 1 + \| \bX_r^{\eps,p} \|^2 + \| \bZ_r^{\eps,p} \|^2 \}  \textup{d} r\\ 
& \quad + \int_0^s  \bX_r^{\eps,p} \cdot \Big( \bsigma \frac{\boeta_r^\eps}{\eps} \Big) \textup{d} r ,
\label{eq:psibound1b}
\end{align}
where we have used \eqref{eq:prop2phi} 
to get the last inequality. 
We want to estimate the last term of the right-hand side of (\ref{eq:psibound1b}).
If we introduce the function $\phi(\bx,\boeta) = \bx \cdot \big( \bsigma {\bf A}^{-1}\boeta\big)$, then we get by (\ref{eq:poisson1}-\ref{eq:expressg})
\begin{align*}
\EE\big[
\eps\phi(\bX_s^{\eps,p},\boeta^\eps_s))-\eps \phi( \bX_0^{\eps,p},\boeta^\eps_0) \big] = &
- \EE\Big[  \int_0^s  \bX_r^{\eps,p} \cdot \Big( \bsigma \frac{\boeta_r^\eps}{\eps} \Big) \textup{d} r\Big] \\
&+\EE\Big[ \int_0^s \big(  \eps \bb^X(\bX_r^{\eps,p},\bZ_r^{\eps,p})  -  \eps \nabla \varphi_p (\bX_r^{\eps,p}) + \bsigma \boeta^\eps_r \big)
\cdot  \big(  \bsigma {\bf A}^{-1}\boeta^\eps_r \big) \textup{d} r\Big] .
\end{align*}
As $\boeta^\eps$ is stationary:
$$
\EE\Big[ \int_0^s \big(    \bsigma \boeta^\eps_r  \big) \cdot 
\big(  \bsigma {\bf A}^{-1}\boeta^\eps_r \big) \textup{d} r\Big] =C_0 s , \quad C_0 = \EE\big[  \big(    \bsigma \boeta^1_0 \big) \cdot 
\big(  \bsigma {\bf A}^{-1}\boeta^1_0 \big) \big].
$$
As $\nabla \varphi_p$ is bounded and $\bb^X$ is Lipschitz, we get
$$
\left|\EE\Big[   \int_0^s  \bX_r^{\eps,p} \cdot \Big( \bsigma \frac{\boeta_r^\eps}{\eps} \Big) \textup{d} r\Big] \right|
\leq
C s + C \eps \Big(\|\bX_0^{\eps,p}\|+\EE[\|\bX_s^{\eps,p}\|^2]^{1/2} +\int_0^s \{ \EE[ \|\bX_r^{\eps,p}\|^2]+ \EE[ \|\bZ_r^{\eps,p}\|^2] \}^{1/2} \textup{d} r \Big).
$$
Therefore, by substituting into (\ref{eq:psibound1b}), we can deduce that 
$$
\sup \limits_{s \leq t} \mathbb{E}   \left[ 
  \| \bX_s^{\eps,p} \|^2 + \| \bZ_s^{\eps,p} \|^2   \right]
\leq C \left ( 1 +  \int_0^t \sup \limits_{s \leq r} \mathbb{E}  
\left[ \| \bX_s^{\eps,p} \|^2 + \| \bZ_s^{\eps,p} \|^2 \right]  \textup{d} r \right) ,
$$
which yields by Gronwall's inequality
$$
\sup \limits_{s \leq t} \mathbb{E}   \left[ 
  \| \bX_s^{\eps,p} \|^2 + \| \bZ_s^{\eps,p} \|^2 \right]
\leq C \exp(C t).
$$
The constant $C$ does not depend on $\eps, p$. 
Substituting into (\ref{eq:psibound1}) and using \eqref{eq:condphipz0} and $\psi_p\geq 0$ (by (\ref{def:yosida})) gives 
\begin{equation}
\label{eq:psibound2}
\sup \limits_{p} \sup \limits_{\eps} \mathbb{E} \left[ \int_0^T \| \nabla \psi_p(\bZ_s^{\eps,p}) \|^2 \textup{d} s \right] < \infty.
\end{equation}
Finally, combining inequalities above, we obtain
\begin{align*}
\mathbb{E}\left[ \sup \limits_{s \leq t} \left \{ \| \bX_s^{\eps,p} - \bX_s^{\eps,q} \|^2 + \| \bZ_s^{\eps,p} - \bZ_s^{\eps,q} \|^2 \right \} \right]
\leq & C \int_0^t \mathbb{E} \left[ \sup \limits_{s \leq r} \left \{ \| \bX_s^{\eps,p} - \bX_s^{\eps,q} \|^2 + \| \bZ_s^{\eps,p} - \bZ_s^{\eps,q} \|^2 \right \} \right]\textup{d} r\\ 
& + C \left ( \frac{1}{p}+\frac{1}{q} \right) ,
\end{align*}
which in turn provides 
$$
\mathbb{E} \left[ \sup \limits_{s \leq t} \left \{ \| \bX_s^{\eps,p} - \bX_s^{\eps,q} \|^2 + \| \bZ_s^{\eps,p} - \bZ_s^{\eps,q} \|^2 \right \}  \right]
\leq C \left ( \frac{1}{p}+\frac{1}{q} \right).
$$
The constant $C$ does not depend on $\eps,p,q$. The proof is complete.
\qed

{\it Proof of the second
 item of Proposition \ref{prop:2}.}
 
 Using the same arguments as in the proof of Proposition \ref{prop:1}
 (note that $\EE\big[  \int_0^s  \bX_r^{0,p} \cdot \big( \bGamma\textup{d} {\itbf W}_r \big) \big]=0$), we get that
\begin{equation}
\label{eq:lem2:bis}
\forall p \geq 1, \:  \mathbb{E} \left[ \sup \limits_{ t \leq T} \left \{ \| \bZ_t^{0,p} - \bX^0_t  \|^2 + \| \bZ_t^{0,p} - \bZ^0_t \|^2\right  \} \right]
\leq \frac{C_T}{p} ,
\end{equation}
where $(\bX^{0,p},\bZ^{0,p})$ is an approximation of $(\bX^0,\bZ^0)$ in the following sense
\begin{equation}
\begin{dcases}
 {\rm d}\bX^{0,p} + \nabla \varphi_p (\bX^{0,p}){\rm d}t = \bb^X(\bX^{0,p},\bZ^{0,p}) {\rm d}t+  \bGamma {\rm d} {\itbf W}_t,\\
 {\rm d}\bZ^{0,p} + \nabla \psi_p (\bZ^{0,p}) {\rm d}t= \bb^Z(\bX^{0,p},\bZ^{0,p}) {\rm d} t. 
\end{dcases}
\end{equation}
The proof is then similar as the one of Proposition \ref{prop:1}.

\section{Proof of Lemma \ref{prop:5}}
\label{app:proofprop5}
Let $\bX^{\eps}$ satisfy (\ref{eq:mode1}).
\textcolor{black}{
We address the evaluation of $\EE[(f(\bX^\eps_t)-f(\bX^0_t))^2]$.}
For any $p>0$, by using Eqs.~(\ref{eq:lem1}) and (\ref{eq:lem1:bis}), we have
\begin{align*}
&
\EE \big[ \big(  f (\bX^\eps_t   ) - f (\bX^0_t   )\big)^2 \big] \\
&
\leq 
4 \EE \big[ \big( f (\bX_t^\eps   ) - f (\bX_t^{\eps,p}   )\big)^2 \big] 
+
4 \EE \big[ \big( f (\bX_t^{\eps,p}   ) - f (\bX_t^{0,p}   )\big)^2 \big] 
+
4 \EE \big[ \big( f (\bX^{0,p}_t   ) - f (\bX^0_t    )\big)^2 \big]  \\
&\leq \frac{8C_T \|\nabla f \|_\infty^2}{p} +
4 \EE \big[ \big( f (\bX_t^{\eps,p}   ) - f (\bX^{0,p}_t   )\big)^2 \big] .
\end{align*}
In order to get an estimate of the last term, we can follow the steps of the proof of Lemma \ref{prop:3} in the same way,
because $\nabla\phi_p$, that appears only in $\Lambda_1$, is bounded uniformly in $p$.
We get 
\begin{align*}
&\EE \big[ \big(  f (\bX^\eps_t   ) - f (\bX^0_t   )\big)^2 \big] \\
&
\leq 
4 \EE \big[ \big( f (\bX_t^\eps   ) - f (\bX_t^{\eps,p}   )\big)^2 \big] 
+
4 \EE \big[ \big( f (\bX_t^{\eps,p}   ) - f (\bX^{0,p}_t   )\big)^2 \big] 
+
4 \EE \big[ \big( f (\bX^{0,p}_t  ) - f (\bX^0_t    )\big)^2 \big]  \\
&\leq \frac{8C_T \|\nabla f \|_\infty^2}{p} +
4 C \eps^2 .
\end{align*}
As this holds true for any $p$, 
this gives the desired result.

Let $(\bX^{\eps},\bZ^{\eps})$ satisfy (\ref{eq:mode2}). 
\textcolor{black}{
We address the evaluation of $\EE[(f(\bX^\eps_t, \bZ_t^\eps )-f(\bX^0_t, \bZ_t^0 ))^2]$.}
For any $p$, by using (\ref{eq:lem2}) and (\ref{eq:lem2:bis}), we have
\begin{align*}
& \EE \big[ \big( f (\bX_t^\eps , \bZ_t^\eps  ) - f (\bX^0_t  ,\bZ^0_t  )\big)^2\big] 
\leq 
4 \EE \big[ \big( f (\bX_t^\eps , \bZ_t^\eps  ) - f (\bX_t^{\eps,p}  ,\bZ_t^{\eps,p}  )\big)^2\big] \\
& +
4 \EE \big[ \big( f (\bX_t^{\eps,p}  ,\bZ_t^{\eps,p}  )-  f (\bX^{0,p}_t  ,\bZ^{0,p}_t  )\big)^2 \big] 
+
4 \EE \big[ \big(  f (\bX^{0,p}_t  ,\bZ^{0,p}_t  )- f (\bX^0_t  ,\bZ^0_t  )\big)^2 \big]  \\
&\leq \frac{8C_T \|\nabla f \|_\infty^2}{p} +
4 \EE \big[ \big( f (\bX_t^{\eps,p}  ,\bZ_t^{\eps,p}  )-  f (\bX^{0,p}_t  ,\bZ^{0,p}_t  )\big)^2 \big] .
\end{align*}
In order to get an estimate of the last term, we can follow the steps of the proof of Lemma \ref{prop:3} 
by keeping track of the bound (\ref{eq:condplastic})
on $\nabla\psi_p$ (that appears only in $\Lambda_1$), and we get
$$
 \EE \big[ \big( f (\bX_t^\eps , \bZ_t^\eps  ) - f (\bX^0_t  ,\bZ^0_t  )\big)^2\big] 
\leq \frac{8C_T \|\nabla f \|_\infty^2}{p} +
4 C(1+p) \eps^2  .
$$
By optimizing in $p$ we get the desired result.

\section{Proof of (\ref{eq:estimadhoc})}
\label{app:C}
We consider the case of an OU noise that satisfies the equation ${\rm d} \eta^\eps = -\eps^{-2} \eta^\eps {\rm d} t + \eps^{-1} {\rm d} W$ and $\eta_0 \sim \mathcal{N}(0,1/2)$. The idea remains the same for the case of a Langevin noise.
Define $X_\star^\eps$ and $X^0_\star$ as follows:
\begin{equation}
X_{\star,t}^\eps \triangleq x_0 + \frac{1}{\eps} \int_0^t \eta_s^\eps {\rm d} s \: \: \mbox{and} \: \: X^0_{\star,t} \triangleq x_0 + W_t.
\end{equation}
We first show that 
\begin{equation}
\label{eq:adhoc1}
\mathbb{E} \Big[ \sup \limits_{0 \leq t \leq T} \left | X_{\star,t}^\eps - X^0_{\star,t} \right |^2  \Big]  = O(\eps^2 |\log \eps|).
\end{equation}
From the equation for $\eta^\eps$, it can be seen that  $X_{\star,t}^\eps = X^0_{\star,t} + \eps (\eta_0-\eta^\eps_t)$.
Thus, we have
\[
\mathbb{E}  \Big[ \sup \limits_{0 \leq t \leq T} \left | X_{\star,t}^\eps - X^0_{\star,t} \right |^2 \Big]
= \eps^2 \mathbb{E} \Big[ \sup \limits_{0 \leq t \leq T} \left | \eta_0-\eta_t^\eps \right |^2 \Big]
\leq \eps^2 +2 \eps^2 \EE \Big[ \sup \limits_{0 \leq t \leq T} \left |\eta_t^\eps \right |^2 \Big] .
\]
The process $\eta_t^\eps$ is a stationary centered Gaussian process with covariance function $\EE[ \eta_t^\eps \eta_{t'}^\eps ] 
=(1/2) \exp(- |t-t'|/\eps^2)$.
By the maximal inequality for the OU process \cite{MR1664394} we get 
$ \EE \Big[ \sup \limits_{0 \leq t \leq T} \left |\eta_t^\eps \right | \Big] \leq C \sqrt{\log( 1+T/\eps^2)}$ and by \cite[Proposition 3.19]{massart} we obtain
$$
\mathbb{E}  \Big[ \sup \limits_{0 \leq t \leq T} \left | X_{\star,t}^\eps - X^0_{\star,t} \right |^2 \Big]
\leq 
\eps^2 + \eps^2 + C^2 \eps^2 \log( 1+T/\eps^2) \leq C' \eps^2 (1+|\log \eps|) .
$$
Next, we show that 
\begin{equation}
\label{eq:adhoc2}
\mathbb{E} \Big[\sup \limits_{0 \leq t \leq T} \left | X_t^\eps - X^0_t \right |^2 \Big] = O(\eps^2 |\log \eps|).
\end{equation}
We can use an explicit formula for $X^\eps$ (resp. $X^0$) that involves $X_\star^\eps$ (resp. $X^0_\star$).
Indeed, $X_t^\eps  = \mathcal{M}_t (X_\star^\eps)$ and $X^0_t  = \mathcal{M}_t (X^0_\star)$ where $\mathcal{M}$ is the self map on the set of continuous functions defined by
$\mathcal{M}_t(f) \triangleq f(t) - \min \limits_{0 \leq s \leq t} \min (0,f(s))$. 
This leads to
$$
\vert X_t^\eps - X^0_t \vert \leq \vert X_{\star,t}^\eps - X^0_{\star,t} \vert 
+ \left \vert \min \limits_{0 \leq s \leq t} \min (0,X_{\star,s}^\eps)  -  \min \limits_{0 \leq s \leq t} \min (0,X^0_{\star,s}) \right \vert.
$$
The second term in the right-hand side can be bounded by using the following inequalites i) $| \min(0,a) - \min(0,b)| \leq |a-b|$ for all $a,b$  
and ii) 
$
\left \vert \min \limits_{0 \leq s \leq t}  f(s)  -  \min \limits_{0 \leq s \leq t} g(s) \right \vert
 \leq  \max \limits_{0 \leq s \leq t} |f(s) - g(s)|,
$ 
for all continuous funtions $f,g$.
 Therefore, 
 \[
\mathbb{E} \Big[ \sup \limits_{0 \leq t \leq T} \left | X_t^\eps - X^0_t \right |^2 \Big]
\leq
4 \mathbb{E} \Big[ \sup \limits_{0 \leq t \leq T} \left | X_{\star,t}^\eps - X^0_{\star,t} \right |^2 \Big] = O(\eps^2|\log \eps|), 
\]
which gives \eqref{eq:estimadhoc}.

\section*{Acknowledgments.} 
LM expresses his sincere gratitude to Prof. Jean Michel Coron and the Sino-French International Associated Laboratory for Applied Mathematics for being supported for travels and housing at Ecole Polytechnique.

\end{document}